\documentclass{article}
\usepackage[left=0.75in, right=0.75in, top=0.75in, bottom=0.75in]{geometry}
\usepackage{amsmath,amssymb,mathtools}
\usepackage{graphicx,color,subfig}
\usepackage{cite}

\usepackage{amsthm}

\usepackage{tikz,calc,bbding}
\usetikzlibrary{shapes, arrows, decorations, decorations.pathreplacing, decorations.pathmorphing, fit, matrix}
 \usepackage[font=small,labelfont=bf]{caption}
\usepackage{pifont}	
\usepackage[labelfont=bf]{caption}
\usepackage{authblk}
\usepackage{cite}
\usepackage{appendix}

\usepackage{hyperref}
\hypersetup{backref,bookmarks,hidelinks,breaklinks=true}

\graphicspath{{epsfiles/}}

\newtheorem{theorem}{Theorem}[section]
\newtheorem{proposition}[theorem]{Proposition}
\newtheorem{lemma}[theorem]{Lemma}

\newtheorem{definition}[theorem]{Definition}
\newtheorem{example}[theorem]{Example}
\newtheorem{problem}{Problem}

\newcommand{\longthmtitle}[1]{\mbox{}{\bf \textit{(#1).}}}

\newcommand{\cplx}{\ensuremath{\mathbb{C}}}
\newcommand{\real}{\ensuremath{\mathbb{R}}}
\newcommand{\realpos}{\ensuremath{\mathbb{R}_{>0}}}
\newcommand{\realnonneg}{\ensuremath{\mathbb{R}_{\ge 0}}}
\newcommand{\integers}{{\mathbb{Z}}}
\newcommand{\intpos}{{\mathbb{N}}}

\DeclareMathOperator*{\argmax}{\arg\max}
\DeclarePairedDelimiter{\ceil}{\lceil}{\rceil}

\newcommand{\setdef}[2]{\{#1 \; | \; #2\}}

\newcommand{\notdivides}{\hspace{-4pt}\not \hspace*{-0.5pt} |\hspace{2pt}}
\newcommand{\divides}{\hspace{1pt}|\hspace{1pt}}

\newcommand{\Ec}{\mathcal{E}}

\newcommand{\Gc}{\mathcal{G}}

\newcommand{\Ic}{\mathcal{I}}
\newcommand{\Jc}{\mathcal{J}}

\newcommand{\Nc}{\mathcal{N}}

\newcommand{\Vc}{\mathcal{V}}
\newcommand{\Wc}{\mathcal{W}}

\newcommand\iotab{{\boldsymbol{\iota}}}

\newcommand{\ones}{\mathbf{1}}

\newcommand{\diag}{\text{diag}}
\newcommand{\until}[1]{\{1,\dots,#1\}}

\newcommand\tr{\text{tr}}

\newcommand\Kok{K\!\!-\!1\!-\!k}

\usepackage{enumitem}

\renewcommand{\footnoterule}{%
  \hspace{3pt} \hrule width 0.4\textwidth height 0.5pt
  \kern 2pt
}

\newcommand{\oprocendsymbol}{\hbox{$\square$}}
\newcommand{\oprocend}{\relax\ifmmode\else\unskip\hfill\fi\oprocendsymbol}
\def\eqoprocend{\tag*{$\square$}}

\newcommand{\new}[1]{{\color{blue} #1}}
\renewcommand\new[1]{#1}

\graphicspath{{epsfiles/}}

\title{\Large \bf Heterogeneity \new{of Central Nodes} Explains the Benefits of Time-Varying Control Scheduling in Complex Dynamical Networks}

\author[1,*]{Erfan Nozari}
\author[2]{Fabio Pasqualetti}
\author[1]{Jorge Cort\'es}
\affil[1]{Department of Mechanical and Aerospace Engineering, University of California, San Diego}
\affil[2]{Department of Mechanical Engineering, University of California, Riverside}
\affil[*]{Corresponding author (email: enozari@ucsd.edu)}
\date{}

\begin{document}

\maketitle
\textbf{Despite extensive research and remarkable advancements in the
  control of complex dynamical networks, most studies and practical
  control methods limit their focus to time-invariant control
  schedules (TICS).  This is both due to their simplicity and the fact
  that the
  benefits of time-varying control schedules (TVCS) have remained
  largely uncharacterized.  In this paper we study networks with
  linear and discrete-time dynamics and analyze the role of network
  structure in TVCS.  First, we show that TVCS can significantly
  enhance network controllability over TICS, especially when applied
  to large networks.  Through the analysis of a
  scale-dependent notion of nodal centrality,
  we then show that optimal TVCS involves the actuation of the most
  central nodes at appropriate spatial scales at all
  times. Consequently, it is the scale-heterogeneity of the
  central-nodes in a network that determine whether, and to what
  extent, TVCS outperforms conventional policies based on TICS.
  \new{Here, scale-heterogeneity of a network refers to how diverse
    the central nodes of the network are at different spatial (local
    vs. global) scales.}  Several analytical results and case studies
  support and illustrate this relationship.}
  

\vskip 20pt 

\section{Introduction}

Many natural and man-made systems, ranging from the
nervous system to power and transportation grids to societies, exhibit
dynamic behaviors that evolve over a sparse and complex network. The
ability to control such network dynamics is not only a theoretically
challenging problem but also a barrier to fundamental breakthroughs
across science and engineering. While multiple studies have addressed
various aspects of this problem, several fundamental questions remain
unanswered, including to what extent the capability of controlling a
different set of nodes over time can improve the controllability of
large-scale, complex networked systems.

Controllability of a dynamical network (i.e., a network that supports
the temporal evolution of a well-defined set of nodal \emph{states})
is classically defined as the possibility of steering its state
arbitrarily around the state space through the application of external
inputs to (i.e., \emph{actuation} of) certain \emph{control
  nodes}~\cite{REK:63}. This raises a fundamental question: how does
the choice of control nodes affect network controllability?
Hereafter, we refer to this as the \emph{control scheduling
  problem}~\cite{YYL-JJS-ALB:11,NJC-EJC-DAV-JSF-CTB:12,AO:14}. Notice
that in this classical setting, attention is only paid to the
\emph{possibility} of arbitrarily steering the network state, but not
to the \emph{difficulty and energy cost} of doing so. This has
motivated the introduction of several controllability metrics to
quantify the required effort in the control scheduling
problem~\cite{GY-JR-YL-CL-BL:12,FP-SZ-FB:14,THS-JL:14,THS-FLC-JL:16,VT-MAR-GJP-AJ:16}. While
a comprehensive solution has remained elusive, these works have
collectively revealed the role of several factors in the control
scheduling problem such as the network size and
structure~\cite{FP-SZ-FB:14}, nodal
dynamics~\cite{NJC-EJC-DAV-JSF-CTB:12} and
centralities~\cite{YYL-JJS-ALB:11,THS-JL:14}, the number of control
nodes~\cite{FP-SZ-FB:14}, and the choice of controllability
metric~\cite{THS-FLC-JL:16}.
\new{This problem has also close connections with the optimal sensor scheduling problem, see, e.g.~\cite{SP-PB-GJP:15,MAB:16,HZ-RA-SS:17,VT-YX-SP-PB-GJP:18} and the references therein.}

The majority of the above literature, however, implicitly relies on
the assumption of time-invariant control schedules (TICS), namely,
that the control node(s) is fixed over time. Depending on the specific
network structure, this assumption may come at the expense of a
significant limitation on its controllability, especially for
large-scale systems where distant nodes inevitably exist relative to
any control node. Intuitively, the possibility of time-varying control
schedules (TVCS), namely, the ability to control different nodes at
different times, allows for targeted interventions at different
network locations and can ultimately decrease the control effort to
accomplish a desired task. On the other hand, from a practical
standpoint, the implementation of TVCS requires the ability to
geographically relocate actuators or the presence of actuation
mechanisms at different, ideally all, network nodes, and more
sophisticated control policies. This leads to a critical trade-off
between the benefits of TVCS and its implementation costs which has
not received enough, if any, attention in the literature.

The significant potential of time-varying schedules for control (and
also sensing, which has a dual interpretation to control) has led to
the design of (sub)optimal \new{sensor~\cite{LZ-WZ-JH-AA-CJT:14,STJ-SLS:15} and control~\cite{YZ-FP-JC:16-cdc,DH-JW-HZ-LS:17} scheduling algorithms} in recent years. 
While constituting a notable leap forward and the benchmark for the methods
developed in this paper, these works are oblivious to the fundamental
question of whether,
and to what extent, TVCS provides an improvement in network
controllability compared to TICS.
Our previous work~\cite{EN-FP-JC:17-acc} has studied the former
question (i.e., whether TVCS provides \emph{any} improvement over
TICS) in the case of undirected networks, but did not consider directed networks 
or, more importantly, addressed the latter
question of how large the relative improvement in network
controllability~is.  Given the trade-off between benefits and costs of
TVCS,
a clear answer to this question is vital for the practical application
of TVCS in real-world complex networks.

In this paper, we address these two questions in the context of
discrete-time linear dynamics evolving over directed networks.
Since the implementation costs of TVCS are greatly
  domain-specific and do not follow any common pattern of dependence
  on the control schedule, we here provide an in-depth analysis of the
  benefits of TVCS. This provides the necessary information for
  comparison with the costs of implementing TVCS in any specific
  application in order to decide between TICS and TVCS. 
  
  To this end,
we show that 
\emph{$2k$-communicability},
a new notion of nodal centrality that we define here, plays a
fundamental role in TVCS. This notion measures the centrality of each
node in the network at different spatial scales.  Throughout this
work, the \emph{spatial scale} (or simply \emph{scale}) of any notion
of centrality is defined as the maximum topological distance between
pairs of nodes that allows them to affect the centrality of each
other, where topological distance between a pair of nodes refers to
the minimum number of edges in the graph of the network that should be
traversed to go from one to the other. In particular, the spatial
scale of degree centrality is $1$, while the spatial scale of
eigenvector centrality is $\infty$.  Based on the distinction between
local and global nodal centralities (i.e., centralities with small and
large spatial scales, respectively), we show that the optimal control
node at every time instance is the node with the largest centrality at
the appropriate scale (i.e., the node with the largest
$2k$-communicability at an appropriate~$k$). Accordingly, our main
conclusion is that the benefit of TVCS
is directly related to the \emph{scale-heterogeneity of central nodes}
in the network: the most benefit is gained in networks where the
highest centrality is attained by various nodes at different spatial
scales, while this benefit starts to decay as fewer nodes dominate the
network at all scales \new{(i.e., scale-homogeneity).}

\new{ Moreover, we provide an extensive discussion of how the
  dynamical adjacency matrix of a network can (and should) be
  extracted from its static connectivity, a vital step that is often
  ignored in the literature. Indeed, our simulation results show that
  this step has a significant effect on the benefit of TVCS, with
  \emph{transmission} networks (networks with states that represent
  physical quantities transmitted over the network) benefiting
  significantly more than \emph{induction} networks (those with
  non-physical states that induce state dynamics over the network)
  from TVCS.  }

\section{Notation and Preliminaries}\label{sec:prelims}

In this section, we introduce our notation and briefly review some preliminary concepts that will be used throughout the work. We use $\real$ and $\intpos$ denote the set of reals and positive integers, respectively.
Given $x \in \real^n$,
$x_i$ and $(x)_i$ refer to its $i$th component. Similarly, $a_{i j}$
and $(A)_{i j}$ refer to the $(i, j)$th entry of $A$, and $a_i$ refers
to its $i$'th column.  
Given a matrix $M \in \real^{n \times n}$, its
trace, determinant, and eigenvalue with smallest magnitude are
denoted by $\tr(M)$, $\det(M)$, and $\lambda_\text{min}(M)$,
respectively.

\subsection{Graph theory}
A weighted undirected graph $\Gc = (\Nc, \Ec, A)$ consists of a vertex set
$\Nc = \until{n}$, an edge set
\begin{align*}
\Ec = \setdef{\{i, j\}}{i \text{ is connected to } j}, 
\end{align*}
and an adjacency matrix $A \in \realnonneg^{n
  \times n}$ where, for any $i, j \in \Nc$, $a_{i j} \ge 0$ is the
weight of the edge \new{from node $j$ to node $i$}.  A path in~$\Gc$ from
node $i$ to $j$ is a finite sequence $\ell_0, \ell_1, \ldots, \ell_p$
of nodes where $\ell_0 = i$, $\ell_p = j$, and $\{\ell_{m - 1},
\ell_m\} \in \Ec$ for $\ell \in \until{p}$.  A cycle is a path with
$\ell_0 = \ell_p$.  For $k \ge 1$, $(A^k)_{i j}$ gives the (weighted)
number of paths of length $k$ between nodes $i$ and~$j$. A regular
graph of degree $k$ is a graph where all the vertices have $k$
neighbors. A strongly regular graph with parameters $(n, k, \lambda,
\mu)$ is a regular graph of $n$ nodes with degree $k$ where any two
adjacent vertices have $\lambda$ common neighbors and any pair of
non-adjacent vertices have $\mu$ neighbors in common. Given a network
$\Gc$ with $n$ nodes, a cone on $G$ is a network with $n + 1$ nodes
where the last one is connected to all others.

\subsection{Network centrality}
We briefly review here three centrality measures with spectral
characterizations.  Consider a network of size $n$ represented by the
adjacency matrix~$A$.

\paragraph{Eigenvector centrality~\cite{PB:07-sn,PB:72}:}
Let $v_i \in \realnonneg$ denote the centrality value of node $i \in
\Nc$.
Eigenvector centrality is based on the idea that the influential nodes
are the ones that are connected to other influential nodes. In other
words, $v_i \propto \sum_{j = 0}^n a_{i j} v_j$ for all $i$. This requires the
  existence of a constant $\lambda > 0$ such that $\lambda v_i =
  \sum_{j = 0}^n a_{i j} v_j$ for all~$i$. In matrix notation, $v
  = [v_1 \ \cdots \ v_n]^T$, this becomes $A v = \lambda v$, which is
  an eigenvalue problem. Since $A$ is non-negative, by the
  Perron-Frobenius Theorem~\cite[Fact 4.11.4]{DSB:09}, there always
  exists a pair $(\lambda, v) \in \realpos \times \realnonneg^n$ such
  that $A v = \lambda v$. This vector~$v$
is thus defined as the vector of (right) eigenvector centralities.
The same argument can be repeated by reversing the direction of
influence flow in the network, leading to the vector of left
eigenvector centralities (i.e., a positive vector $u$ such that $u^T A
= \lambda u^T$).

\paragraph{Exponential and resolvent
  communicability~\cite{EE-NH:08,CK:13}:}\label{subsubsec:comm}
The communicability of a node measures its
ability to communicate with the rest of the network. Different notions
of communicability have been proposed for complex networks. For a
given node~$i$, these include exponential communicability
$(e^{\beta A})_{i i}$ and the resolvent communicability
$((I - \beta A)^{-1})_{i i}$, respectively, where $\beta > 0$.
From the power series
expansion of $e^{\beta A}$ and $(I - \beta A)^{-1}$, it follows that
the exponential and resolvent communicabilities count the total number
of cycles that pass through node $i$, weighting the ``importance'' of
cycles of length $k$ by $\beta^k / k!$ and $\beta^k$,
respectively. Thus, the role of $\beta$ is to determine how local/global 
these measures are: increasing $\beta$ increases the weights
of longer cycles. One can show~\cite{CK:13} that in the extreme cases
of $\beta \to \infty$ in the exponential case and
$\beta \to \frac{1}{\lambda_\text{max}(A)}$ in the resolvent case,
both notions result in the same rankings of nodes as eigenvector
centrality.

\paragraph{Degree centrality:}
The degree centrality of node $i$ is the sum of the $i$-th row (or
column) of $A$ and provides a measure of the immediate influence of
node $i$ on its neighbors.

\new{
\sloppy
\section{Problem Statement: Comparison of Time-Varying and Time-Invariant Control Scheduling}\label{subsec:model}
}

\new{We consider a network of $n$ nodes that communicate
over a graph $\Gc = (\Nc, \Ec, A)$ that is in general weighted and directed 
(see Appendix~\ref{sec:c2a} for methods of
obtaining $A$ from network connectivity structure). Each node
$i$ has a \emph{state} value $x_i \in \real$ that evolves over time
through the interaction of node $i$ with its neighbors in $\Gc$ and an external control $u$. Assuming that these interactions are linear and time-invariant, we have}
\begin{align}\label{eq:dyn}
  x(k + 1) = A x(k) + b(k) u(k), \quad k \in \{0, \dots, K - 1\},
\end{align}
where \new{$x = (x_1, \dots, x_n) \in \real^n$ is the network state,} $u(k) \in \real$ is the control input, $b(k) \in \real^n$ is
the \emph{time-varying} input vector, 
and $K$ is the time
horizon. For simplicity of exposition, we consider only one
  control input at a time, but the discussion is generalizable to
  multi-input networks (cf. Appendix~\ref{sec:mi}). Define 
  \begin{align}\label{eq:iota}
  \iota_k \in \Nc,
  \end{align}
  to be the index of the node to which the control signal $u(k)$ is
  applied at time $k$. Then, $b(k)$ is equal to the $\iota_k$'th
column of the identity matrix. For the sake of simplicity, we
  here assume that all the network nodes are actuatable, so $\iota_k
  \in \Nc$. \new{If a subset of nodes are \emph{latent}, (i.e., not actuatable),} further
  challenges arise and thus we postpone the analysis of this case to
  Section~\ref{sec:latent}.
  
\new{
The dynamical network~\eqref{eq:dyn} is \emph{controllable} if its state can be steered from arbitrary $x(0) = x_0$ to arbitrary $x(K) = x_f$ using the control input $\{u(k)\}_{k = 0}^{K - 1}$ or, equivalently, if the \emph{controllability Gramian}
\begin{align}\label{eq:Gramian}
  \Wc_K = \sum_{k = 0}^{K - 1} A^k b(\Kok) b(\Kok)^T (A^T)^k,
\end{align}
is nonsingular~\cite{CTC:98}. In general, the eigenvalues of $\Wc_K$ determine how large the \emph{unit-energy reachability set} (the set of states $x_f$ that can be reached from the origin $x_0 = 0$ using controls with unit energy) is (cf. Appendix~\ref{sec:measures} for derivation). Therefore, various measures of controllability based on the eigenvalues of $\Wc_K$ have been proposed, most notably $\tr(\Wc_K)$, $\tr(\Wc_K^{-1})^{-1}$, $\det(\Wc_K)$, $\lambda_{\min}(\Wc_K)$. Each metric has its own benefits and limitations, on which we elaborate more in the following.
}

\begin{center} 
\begin{figure*}
\begin{center} 
\begin{tikzpicture} 
\node (intro-ex-TV) {\parbox{0.15\textwidth}{
\subfloat[]{
\includegraphics[width=0.13\textwidth]{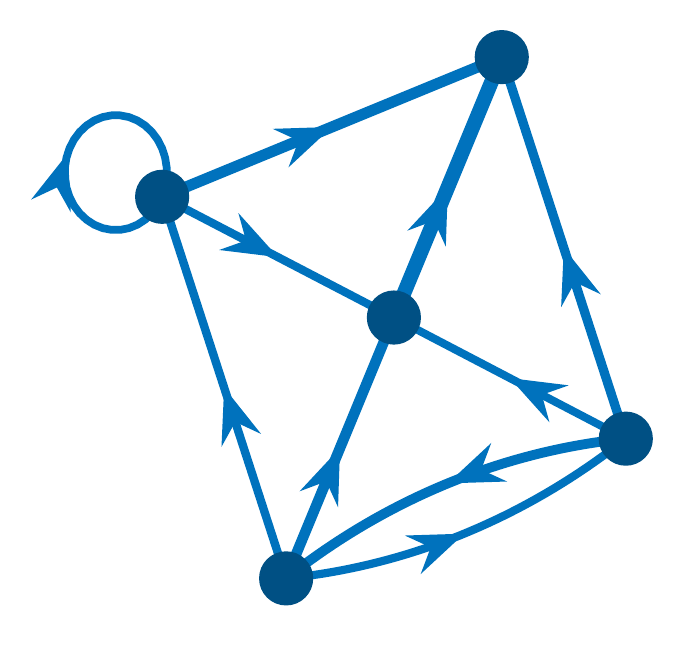}
}}};
\node[below of=intro-ex-TV, xshift=-10, yshift=-50pt] (table-TV) {\parbox{0.3\textwidth}{
\subfloat[]{\scriptsize
\addtolength\tabcolsep{-3pt}
\begin{tabular}{c|ccc}
$f(\cdot)$ & $f_{\max}^\text{TI}$ & $f_{\max}^\text{TV}$ & $\chi$
\\ \hline
$\tr(\cdot)$ & $2.00$ & $2.70$ & $0.35$
\\
$\tr((\cdot)^{-1})^{-1}$ & $1.26 \!\!\times\!\! 10^{-7}$ & $8.22 \!\!\times\!\! 10^{-4}$ & $6.5 \!\!\times\!\! 10^3$
\\
$\det(\cdot)$ & $9.90 \!\!\times\!\! 10^{-11}$ & $7.42 \!\!\times\!\! 10^{-10}$ & $6.49$
\\
$\lambda_{\min}(\cdot)$ & $1.27 \!\!\times\!\! 10^{-7}$ & $1.10 \!\!\times\!\! 10^{-4}$ & $8.7 \!\!\times\!\! 10^2$
\end{tabular}
}}};
\node[right of=table-TV, xshift=170pt, yshift=60pt] (air500) {\parbox{0.3\textwidth}{
\subfloat[]{
\includegraphics[width=0.3\textwidth]{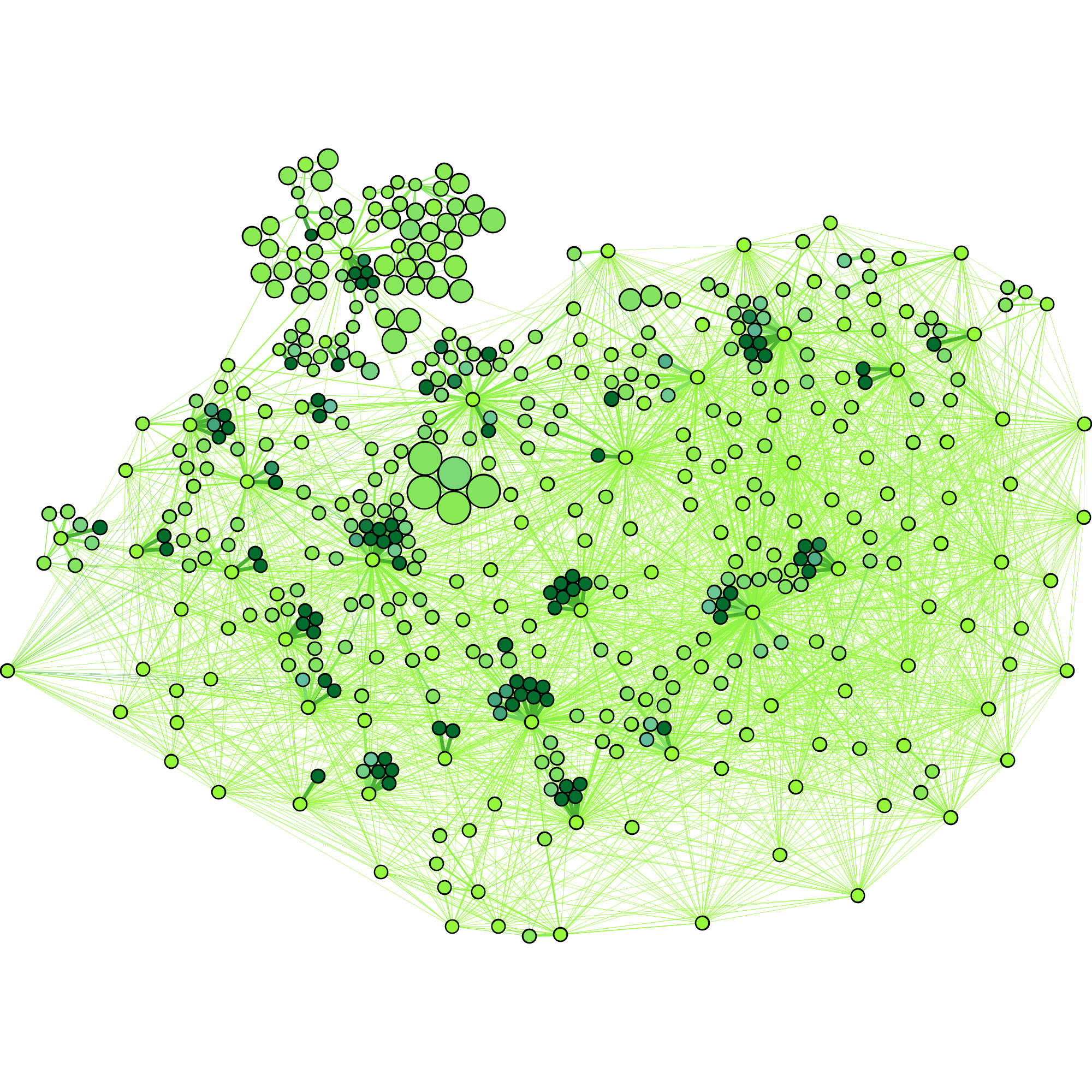}
}}};
\node[below of=intro-ex-TV, xshift=-40, yshift=-140pt] (intro-ex-TI) {\parbox{0.15\textwidth}{
\subfloat[]{
\includegraphics[width=0.3\textwidth]{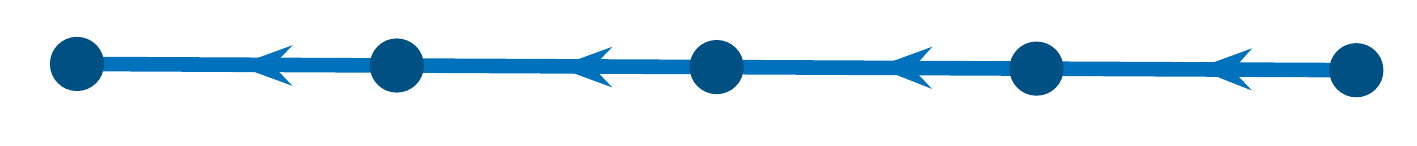}
\
}}};
\node[below of=intro-ex-TI, xshift=115pt, yshift=-35pt] (table-TI) {\parbox{0.5\textwidth}{
\subfloat[]{\scriptsize
\addtolength\tabcolsep{-3pt}
\begin{tabular}{c|ccc}
$f(\cdot)$ & $f_{\max}^\text{TI}$ & $f_{\max}^\text{TV}$ & $\chi$
\\ \hline
$\tr(\cdot)$ & $5$ & $5$ & $0$
\\
$\tr((\cdot)^{-1})^{-1}$ & $0.2$ & $0.2$ & $0$
\\
$\det(\cdot)$ & $1$ & $1$ & $0$
\\
$\lambda_{\min}(\cdot)$ & $1$ & $1$ & $0$
\end{tabular}
}}};
\node[right of=table-TI, xshift=90pt, yshift=40pt] (uci_forum) {\parbox{0.3\textwidth}{
\subfloat[]{
\includegraphics[width=0.3\textwidth]{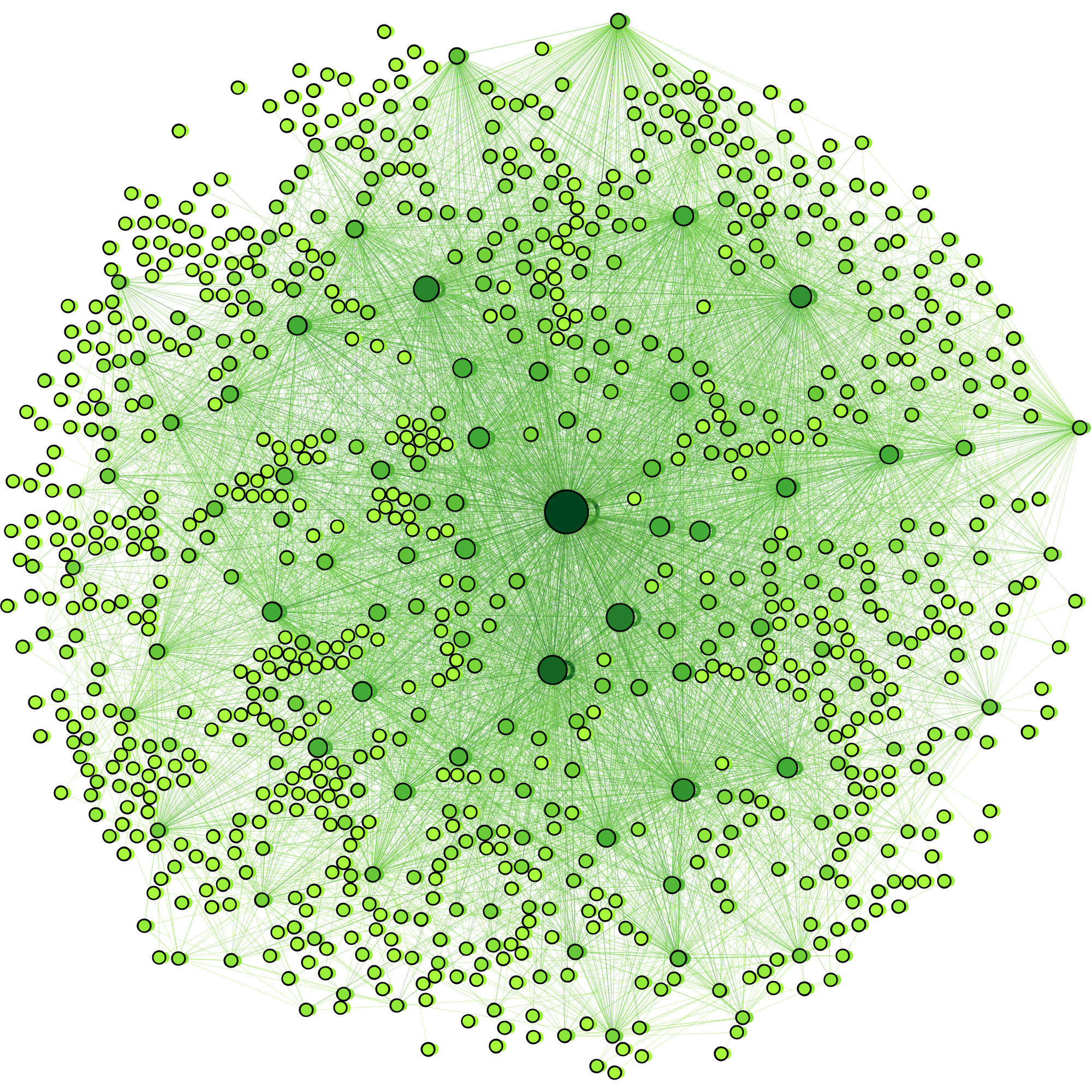}
}}};
\end{tikzpicture} 
\end{center} 
\caption{Advantage of TVCS in dynamic networks. \textbf{(a)} A small example
  network of $n = 5$ nodes. The thickness of each edge $(i, j)$
  illustrates its weight $a_{ij}$. \textbf{(b)} The optimal values of TICS and
  TVCS (equations~\eqref{eq:opt-TI} and~\eqref{eq:opt}, respectively)
  and the relative TVCS advantage (equation~\eqref{eq:chi}) for the network in
  (a). \textbf{(c)} An air transportation network among the busiest airports in
  the United States (see 'air500' in Table~\ref{tab:real} for
  details). The network is undirected, and the dynamical adjacency
  matrix $A$ is computed from static connectivity using the
  transmission method (cf. Appendix~\ref{sec:c2a}). This is an example of a
  network that significantly benefits from TVCS with $\chi \simeq 20\%$. \textbf{(d)} A small example network of
  the same size as (a) but with no benefit from TVCS. \textbf{(e)} The optimal
  values of TICS and TVCS (equations~\eqref{eq:opt-TI}
  and~\eqref{eq:opt}, respectively) and the relative TVCS
  advantage (equation~\eqref{eq:chi}) for the network in (d). We see that the
  network does not benefit from TVCS independently of the choice of
  controllability metric. \textbf{(f)} A social network of students at the
  University of California, Irvine (see 'UCI Forum' in
  Table~\ref{tab:real} for details). Similar to (c), the network is
  undirected and the adjacency matrix is computed using the
  transmission method. This network, however,
  does not benefit from TVCS ($\chi = 0$). In (c) and (f), the
  controllability measure of equation~\eqref{eq:tr} is used due to the large size
  of the network. In both cases, the color intensity and size of nodes
  represent their values of $R_i(1)$ and $R_i(K - 1)$, respectively
  ($K = 10$). While there is a close correlation between nodal size
  and color intensity in (f) (i.e., the darkest nodes are also the
  largest), this is not the case in (c), which is the root cause for
  the difference in their $\chi$-values.}
\label{fig:tv-adv}
\end{figure*}
\end{center}

Assume, \new{for now, that} $f(\Wc_K) \ge 0$ is any of the aforementioned controllability
measures. \new{In \emph{optimal control scheduling},} we seek to choose the control nodes
$\{\iota_k\}_{k = 0}^{K - 1}$ (or, equivalently,
$\{b(k)\}_{k = 0}^{K - 1}$) optimally. 
The conventional approach in the
literature~\cite{YYL-JJS-ALB:11,NJC-EJC-DAV-JSF-CTB:12,GY-JR-YL-CL-BL:12,FP-SZ-FB:14,AO:14,THS-JL:14,VT-MAR-GJP-AJ:16,THS-FLC-JL:16}
\new{is to assume a constant control node, thus called the \emph{time-invariant control scheduling (TICS) problem}:}
\begin{subequations}\label{eq:opt-TI}
  \begin{align}
    \text{\new{TICS:}} \hspace{50pt} \max_{\iota_0, \dots, \iota_{K - 1} \in \Nc} \hspace{10pt} &f(\Wc_K) \hspace{100pt}
    \\
    \text{s.t.} \hspace{25pt} &\iota_0 = \cdots = \iota_{K - 1}
  \end{align}
\end{subequations}
The main advantage of TICS is its simplicity, from theoretical,
computational, and implementation perspectives.  However,
this simplicity comes at a possibly significant cost in terms of
network controllability, compared to the case where the control nodes
$\{\iota_k\}_{k = 0}^{K - 1}$ are independently chosen, namely,
\begin{align}\label{eq:opt}
\text{\new{TVCS:}} \hspace{50pt} \max_{\iota_0, \dots, \iota_{K - 1} \in \Nc} \hspace{10pt} f(\Wc_K). \hspace{100pt}
\end{align}
This approach, namely, time-varying control scheduling (TVCS), is at
least as good as TICS, but has the potential to improve network
controllability significantly. Figure~\ref{fig:tv-adv}(a-b) illustrates a small
network of $n = 5$ nodes together with the optimal values of
equations~\eqref{eq:opt-TI} and~\eqref{eq:opt} and the relative
advantage of TVCS over TICS, defined as
\begin{align}\label{eq:chi}
\chi = \frac{f_{\max}^\text{TV} - f_{\max}^\text{TI}}{f_{\max}^\text{TI}}.
\end{align}
Three observations are worth highlighting. First, the value of $\chi$
is extremely dependent on the choice of controllability measure $f$,
and different choices lead to orders of magnitude change in
$\chi$. Second, the relative advantage of TVCS over TICS is
significant for all choices of the controllability measure, with the
minimum improvement of $\chi = 35\%$ for the choice of $f(\cdot) =
\tr(\cdot)$. The fact that $f(\cdot) = \tr(\cdot)$ results in the
smallest value of $\chi$ relative to other measures is consistently
observed in synthetic and real-world networks, and stems from the fact
that $\tr(\Wc_K)$ has the smallest sensitivity (greatest robustness)
to the choice of control schedule. Finally, even with optimal TVCS,
$\lambda_{\min}(\Wc_K)$ is orders of magnitude less than $1$,
indicating the inevitable existence of very hard-to-reach directions
in the state space. 
  This shows that efficient
  controllability cannot be maintained in all directions in the state
  space even using TVCS and even in very small networks with control over $1/5
  = 20\%$ of the nodes. Except for $\tr(\Wc_K)$, all the measures rely
  heavily on this least-controllable direction, while $\tr(\Wc_K)$
  trades this off for improved controllability in the most efficient
  directions in the state space. \new{See Appendix~\ref{sec:measures} for further discussion of this tradeoff.}

Despite the significant increase in size and complexity, the same core
principles outlined above apply to controllability of real-world
networks. The large size of these networks, however, imposes new
constraints on the choice of the controllability measure $f$ that make
the use of $f(\cdot) = \lambda_{\min}(\cdot)$,
$\tr((\cdot)^{-1})^{-1}$, and $\det(\cdot)$ numerically infeasible and
theoretically over-conservative, as discussed in detail in Appendix~\ref{sec:measures}. As a result, we resort to
the particular choice of controllability measure
\begin{align}\label{eq:tr}
  f(\Wc_K) = \tr(\Wc_K),
\end{align}
for networks beyond $n \simeq 15$.  Since this measure has the
smallest sensitivity to the choice of $\{\iota_k\}_{k = 0}^{K - 1}$
(Figure~\ref{fig:tv-adv}(b)), we expect any network that benefits from
TVCS using the choice of equation~\eqref{eq:tr} to also benefit from it using other
Gramian-based measures (while the converse is not necessarily true,
i.e., there are networks that significantly benefit from TVCS using
other measures but show no benefit in terms of
$\tr(\Wc_K)$). Figure~\ref{fig:tv-adv}(c) illustrates an air
transportation network among the busiest airports in the United
States, comprising of $n = 500$ nodes. 
\new{Using~\eqref{eq:tr}, we see $\chi \simeq 20\%$ improvement in
controllability, verifying our expectation about the
benefits of TVCS.}

In spite of this potential benefit, TVCS has usually higher
computational and implementation costs.
 These include the higher computational cost of computing the
optimal TVCS, and that of installing an actuator at several (ideally
all) nodes of the network. Further, not all networks benefit from TVCS
alike. A simple directed chain network with the same size as that of
Figure~\ref{fig:tv-adv}(a) gains absolutely no benefit from TVCS,
independently of the choice of $f$ 
(Figure~\ref{fig:tv-adv}(d-e)). \new{Similarly, $\chi = 0$ is
also observed in larger, complex networks,
indicating that the optimal TVCS and the optimal TICS are the same (Figure~\ref{fig:tv-adv}(f)).}

\new{ These observations collectively raise a fundamental question
  that constitutes the main problem studied in this paper. Before
  formally stating the problem, we need a definition for ease of
  reference.

\begin{definition}\longthmtitle{Class $\Vc$ and $\Ic$ networks}\label{def:VI}
  Consider a dynamical network described by~\eqref{eq:dyn} and the
  measure $\chi$ introduced in~\eqref{eq:chi}. We say that the network
  belongs to class $\Vc$ if it has $\chi > 0$ and we say it belongs to
  class $\Ic$ otherwise ($\chi = 0$).
\end{definition}

In words, class $\Vc$ networks are those that benefit from TVCS and
class $\Ic$ networks are those that do not. Our main problem of
interest is then as follows.

\begin{problem}\label{prob}
  Given the set of all dynamical networks described by dynamics of the
  form~\eqref{eq:dyn}, characterize the sets $\Vc$ and $\Ic$ in terms
  of the network structure $A$ and develop efficient and
  easy-to-interpret methods for distinguishing between them.
\end{problem}
}

In the following, we restrict our attention to the choice
of controllability measure in equation~\eqref{eq:tr} due to its applicability to
all network sizes and carry a thorough analysis of its properties in
order to address Problem~\ref{prob}.

\new{
\section{Main Results}\label{sec:res}

In this section, we present our main results regarding
Problem~\ref{prob}. First, we introduce a new notion of
communicability that is pivotal to the solution of
Problem~\ref{prob}. Then, we present our results regarding the
characterization of class $\Vc$ and $\Ic$ networks and, finally, study
the case of networks with latent nodes declared earlier.  }

\new{
\subsection{$2k$-Communicability and Scale-Heterogeneity}\label{sec:2k-comm}
}

Consider the TVCS problem in equation~\eqref{eq:opt} with $f(\cdot) =
\tr(\cdot)$. Using the definition of the controllability Gramian \new{in~\eqref{eq:Gramian}} and the
invariance property of trace under cyclic permutations, we can write
\begin{align*}
  \tr(\Wc_K) = \sum_{k = 0}^{K - 1} b(\Kok)^T (A^k)^T A^k b(\Kok).
\end{align*}
Therefore,
\begin{align*}
\max_{\iota_0, \dots, \iota_{K - 1}} \tr(\Wc_K) = \sum_{k = 0}^{K - 1} \max_{\iota_{\Kok}} b(\Kok)^T (A^k)^T A^k b(\Kok),
\end{align*}
where each term $b(\Kok)^T (A^k)^T A^k b(\Kok)$ is the
$\iota_{\Kok}$'th diagonal entry of $(A^k)^T A^k$ \new{(cf. equation~\eqref{eq:iota})}.  Therefore, the
optimization in~\eqref{eq:opt} boils down to finding the largest diagonal
element of $(A^k)^T A^k$ and applying $u(\Kok)$ to this node.
\new{On the other hand, for the TICS problem in~\eqref{eq:opt-TI}} we have
\begin{align*}
  \tr(\Wc_K) = b^T \Bigg(\sum_{k = 0}^{K - 1} (A^k)^T A^k\Bigg) b,
\end{align*}
so one has to instead find the largest diagonal entry of $\sum_{k =
  0}^{K - 1} (A^k)^T A^k$ and apply all the control inputs $u(0),
\dots, u(K - 1)$ to this same node, which is clearly sub-optimal with
respect to TVCS. 
\new{This discussion motivates the following definition.

\begin{definition}\longthmtitle{$2k$-communicability}\label{def:2kcomm}
Given the network dynamics~\eqref{eq:dyn}, the $2k$-communicability of a node $i \in \Nc$ is defined as
\begin{align}\label{eq:Rik}
R_i(k) = ((A^k)^T A^k)_{i i}, \qquad i \in \Nc, \quad k \ge 0.
\end{align} 
\end{definition}
}

\begin{figure*}
\begin{center}
\begin{tikzpicture} 
\node (graph) {\parbox{0.33\textwidth}{
\subfloat[]{
\includegraphics[width=0.15\textwidth]{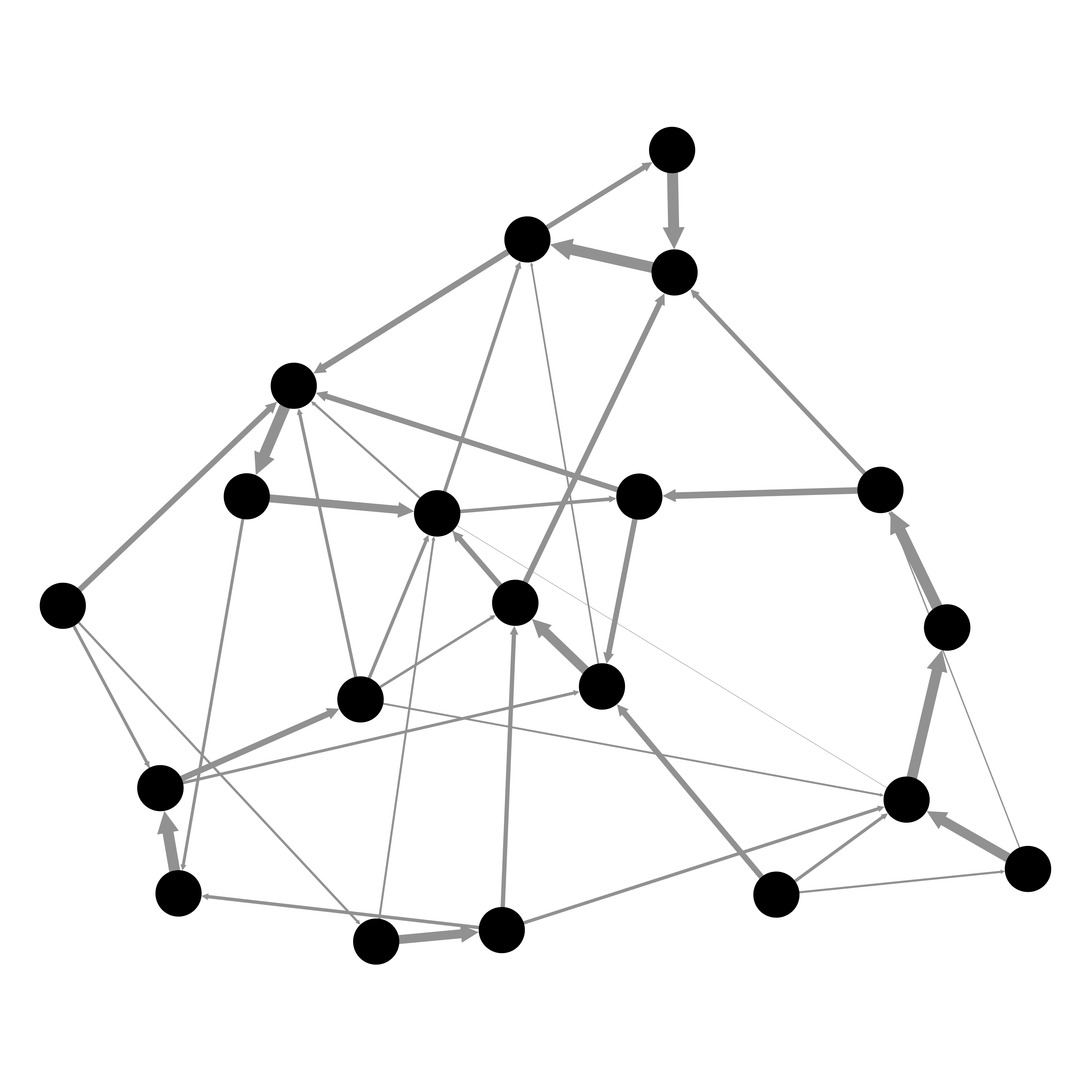}
}}};
\node[below of=graph, xshift=-30pt, yshift=-90pt] (comm) {\parbox{0.45\textwidth}{
\subfloat[]{
\includegraphics[width = 0.4\textwidth]{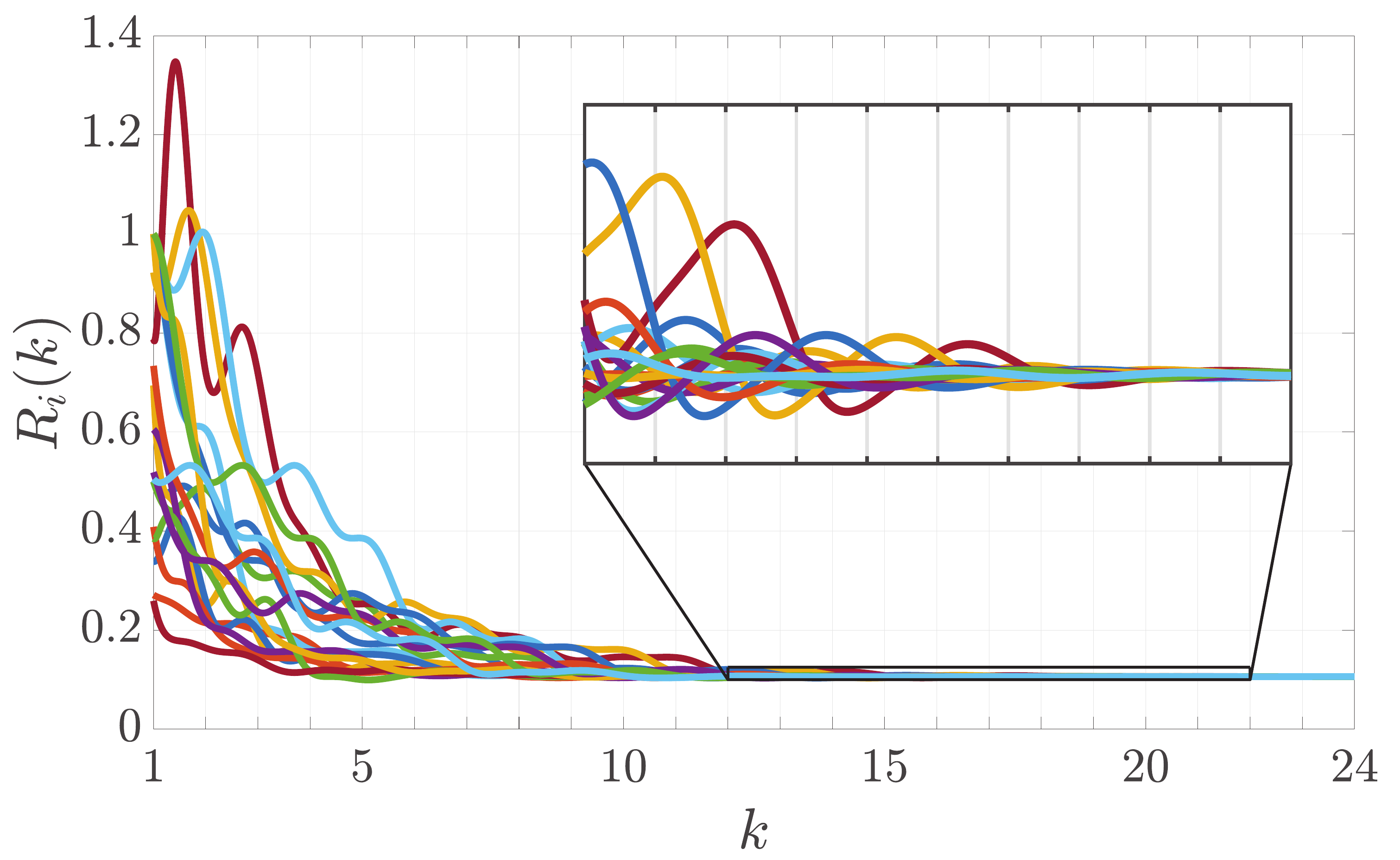}
}}};
\node[right of=graph, xshift=140pt] (0) {\parbox{0.3\textwidth}{
\subfloat[]{
\includegraphics[width=0.2\textwidth]{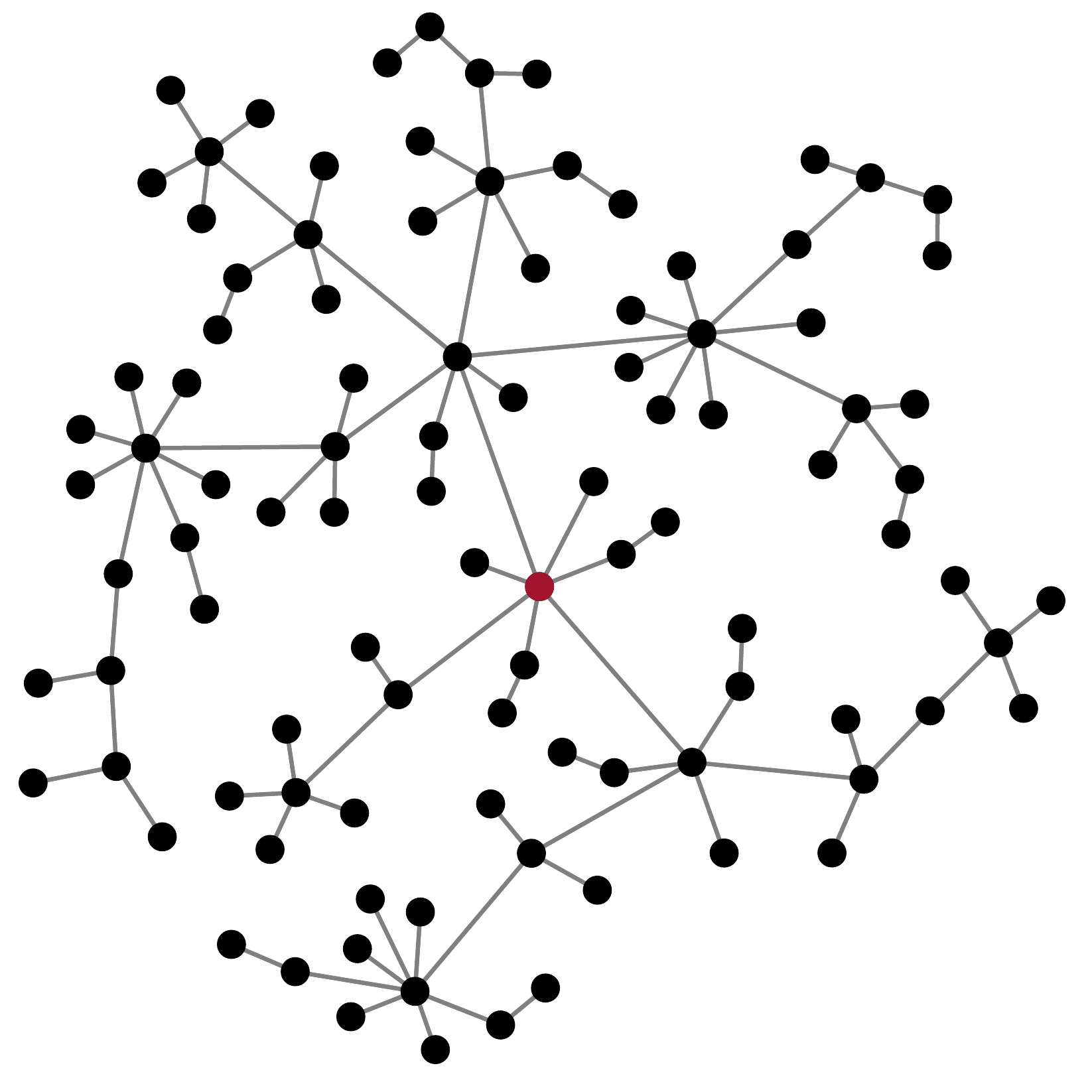}
}}};
\node[right of=0, xshift=100pt] (1) {\parbox{0.3\textwidth}{
\subfloat[]{
\includegraphics[width=0.2\textwidth]{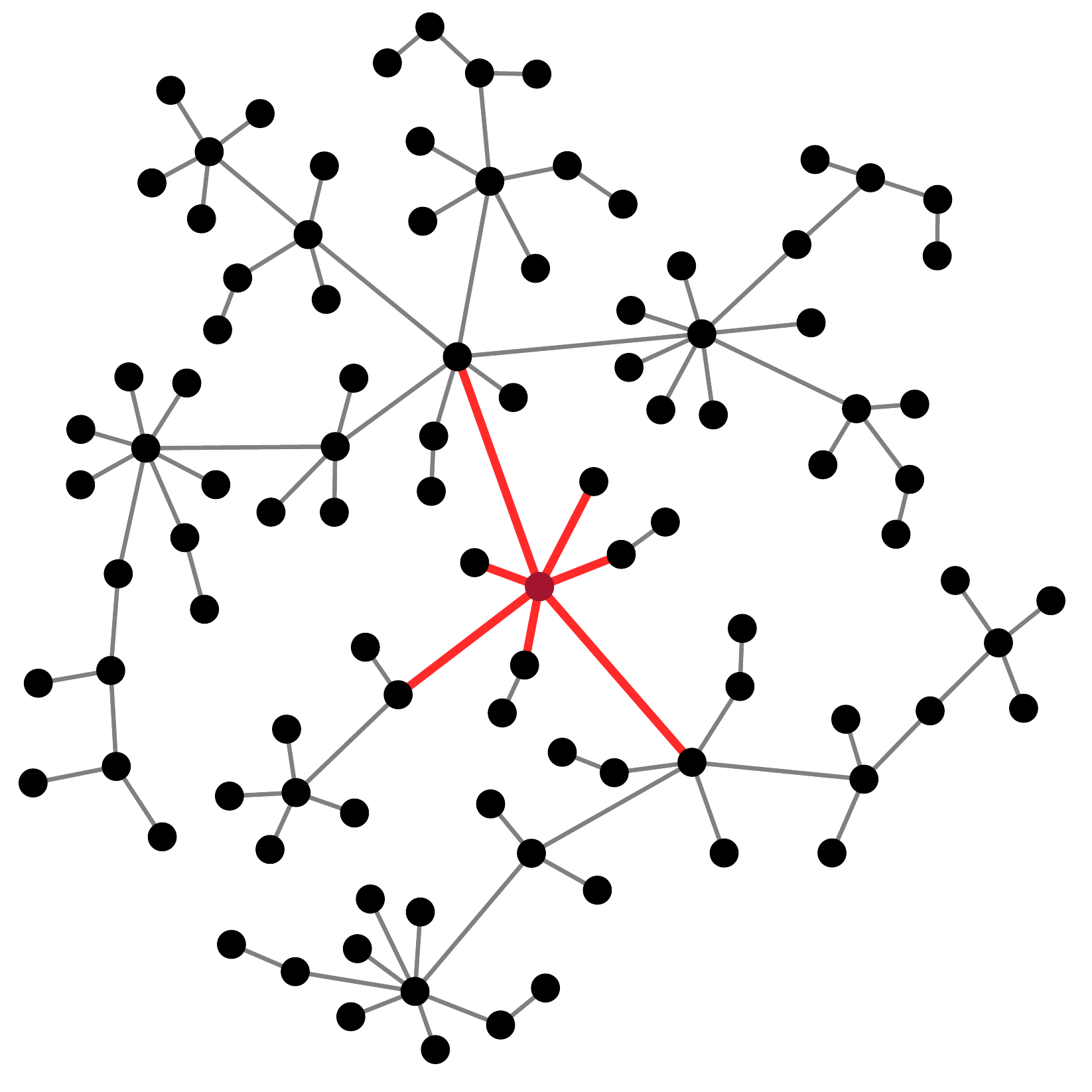}
}}};
\node[below of=0, yshift=-100pt] (2) {\parbox{0.3\textwidth}{
\subfloat[]{
\includegraphics[width=0.2\textwidth]{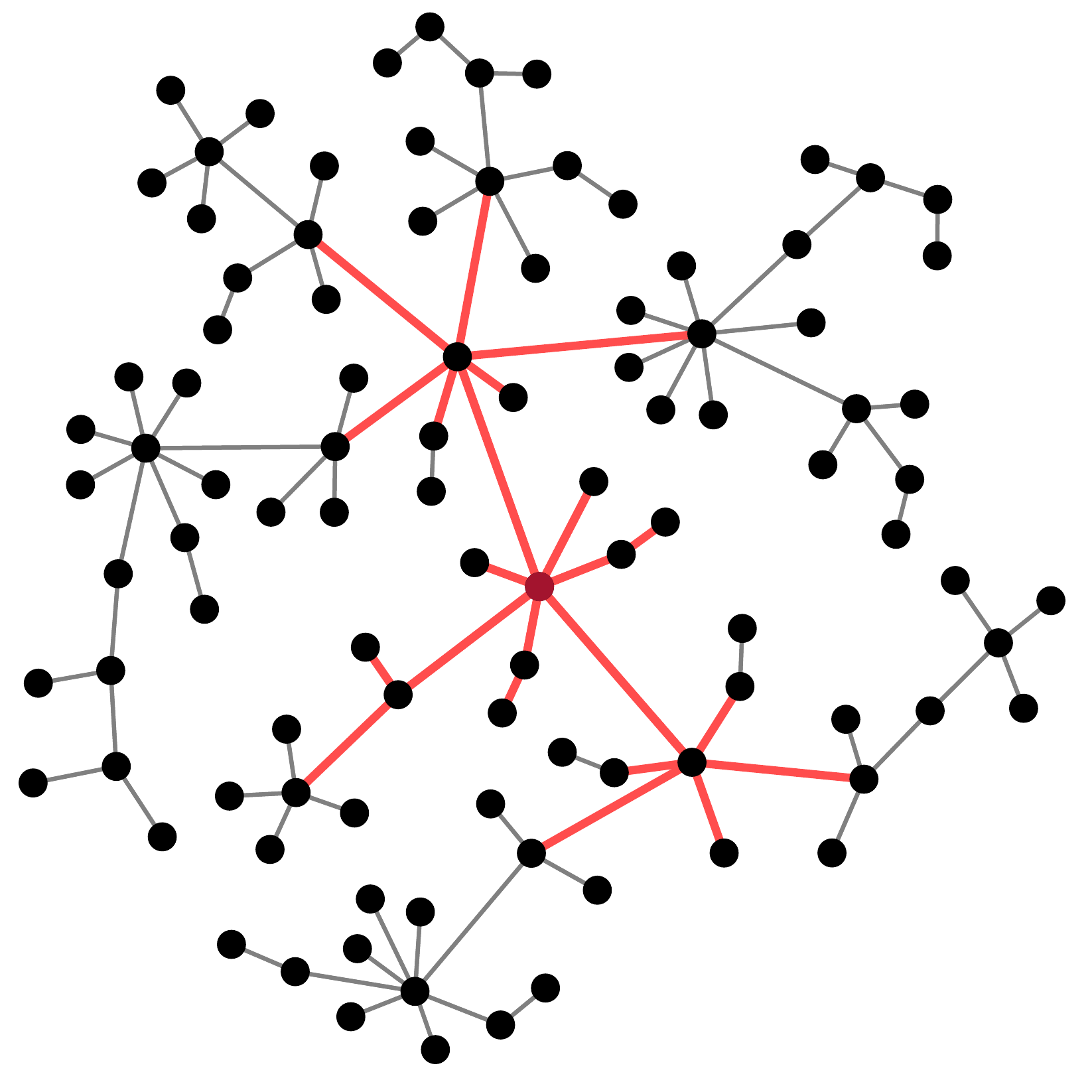}
}}};
\node[right of=2, xshift=100pt] (3) {\parbox{0.3\textwidth}{
\subfloat[]{
\includegraphics[width=0.2\textwidth]{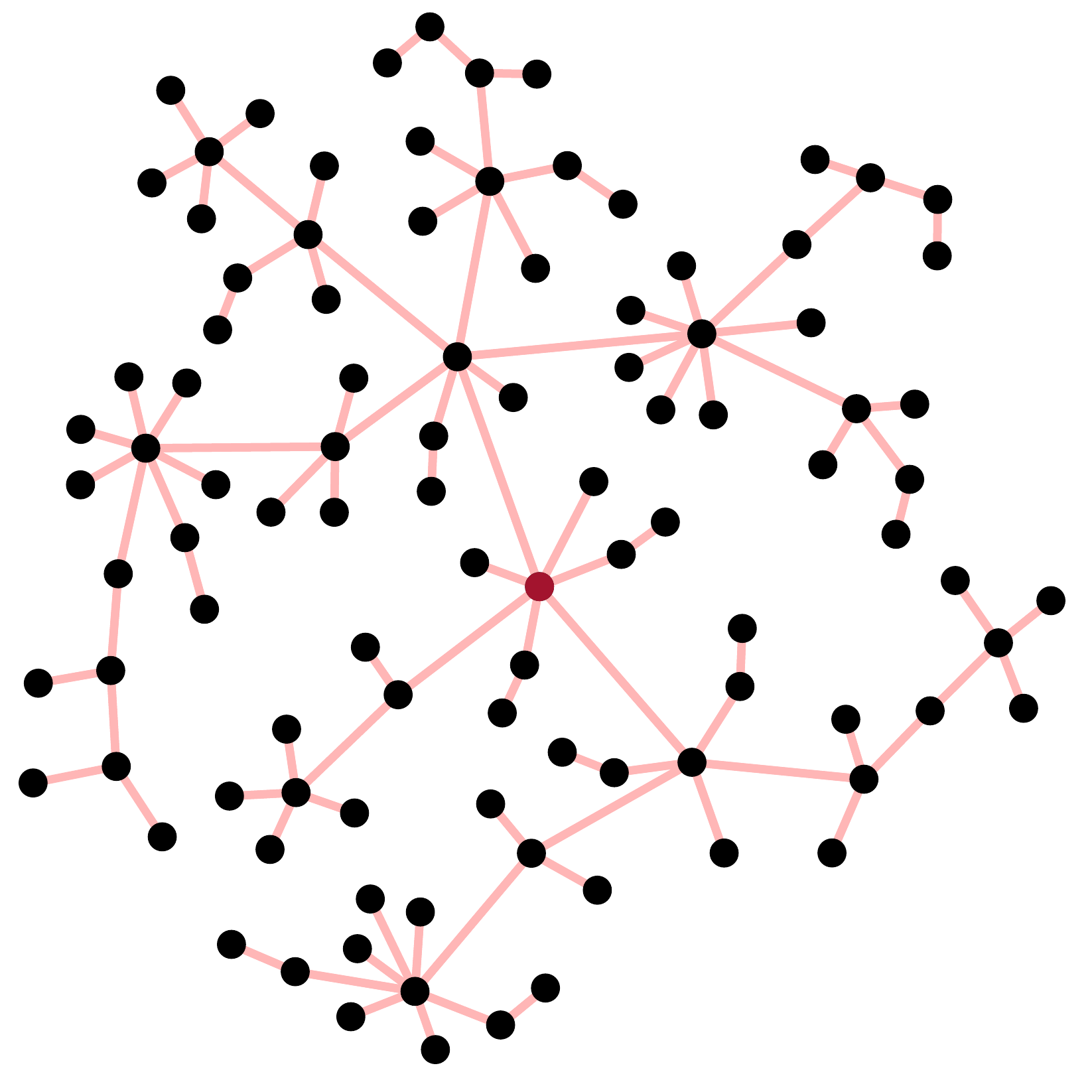}
}}};
\end{tikzpicture} 
\end{center} 
\caption{$2k$-communicability of dynamical networks. \textbf{(a)} An example
  network of $n = 20$ nodes for illustration of the dependence on $k$
  of nodal $2k$-communicabilities. The thickness of the edges is
  proportional to their weights. \textbf{(b)} The evolution of the functions
  $\{R_i(k)\}_{i = 1}^n$. Although these functions are originally only
  defined over integer values of $k$, we have extended their domain to
  real numbers for better illustration of their crossings and
  oscillatory behavior. Oscillatory behavior only arises when $A$ has
  complex-valued eigenvalues (otherwise, $R_i(k)$ is strictly convex).
  \textbf{(c)} An example network of $n = 100$ nodes for illustration of the
  scaling property of $2k$-communicability. The node whose
  $2k$-communicabilities are to be computed (i.e., ``node $i$'') is
  depicted in red. \textbf{(d-f)} The $2$-, $4$-, and $14$-communicability of
  the node depicted in red, as determined by its $1$-, $2$-, and
  $7$-hop incoming paths. We see that $R_i(1)$ only depends on the
  immediate (out-)neighbors of $i$, but as $k$ grows, $R_i(k)$ encodes
  more global information.}
\label{fig:comm}
\end{figure*}

Figure~\ref{fig:comm}(a-b) illustrates the
evolution of $R_i(k)$ as a function of $k$ for all $i \in \Nc$ for a
sample network of $n = 20$
nodes.

Perhaps the most salient property of $2k$-communicability is the
extent to which it relies on the local interactions among the
nodes. Recall, cf.~\cite{CDG-GFR:01}, that for any $k$, the $(i, j)$
entry of $A^k$ equals the total number of paths of length $k$ from
node $i$ to $j$ (if the graph is weighted, each path counts as its
weight, equal to the product of the weights of its
edges). From equation~\eqref{eq:Rik}, we see that $R_i(k)$ equals the sum
of the squares of the total (weighted) number of paths of length $k$
ending in node $i$. In other words, $R_i(k)$ only depends on
connections of node $i$ with its $k$-hop out-neighbors, and is
independent of the rest of the network. Therefore, $R_i(k)$ is a local
notion of centrality for small $k$ and it incorporates more global
information as $k$ grows. In particular, as shown in Appendix~\ref{sec:relations}, $R_i(k)$ is closely related
to
\begin{itemize}
\item the out-degree centrality of node $i$ for $k = 1$;
\item the left eigenvector centrality of node $i$ for $k \to \infty$.
\end{itemize} 
\new{This scaling property of $2k$-communicability is illustrated in Figure~\ref{fig:comm}(a-d) for an example network of $n = 100$ nodes. 
Accordingly, we take the left eigenvector centrality squared as the definition of $R_i(\infty)$ in the sequel.}

\new{
The scaling property of $2k$-communicability also plays an important role in Problem~\ref{prob}.
For ease of reference, let
\begin{align*} 
r(k) \in \Nc
\end{align*}
denote the index of the node that has the largest $R_i(k)$. 
\new{Then, according to the discussion above,
\begin{align}\label{eq:iota-r}
  \iota_k^* = r(K - 1- k),
\end{align}
which forms the core connection between $2k$-communicability and TVCS.} From
this, we see that the optimal TVCS involves the application of $u(0)$
to the node $r(K - 1)$ with the highest global centrality and
gradually moving the control node until we apply $u(K - 2)$ to the
node $r(1)$ with the highest local centrality
(the control node at time $K - 1$ is arbitrary as
  $R_i(0) = 1$ for all $i$).  The intuition behind this procedure is
simple. At $k = 0$, the control input has enough time to propagate
through the network, which is why the highest globally-central node
should be controlled. As we reach the control horizon $K$, the control
input has only a few time steps to disseminate through the network,
hence the optimality of locally-central nodes.
This further motivates our definition of \emph{scale-heterogeneity}, as follows.

\begin{definition}\longthmtitle{Scale-heterogeneity of dynamical networks}\label{def:het}
Consider the network dynamics~\eqref{eq:dyn} subject to the TVCS problem~\eqref{eq:opt} with $2k$-communicability as defined in Definition~\ref{def:2kcomm}. The network is called \emph{scale-homogeneous} if $r(1) = r(2) = \cdots = r(\infty)$ and scale-heterogeneous otherwise. Accordingly, the more varied $\{r(k)\}_{k = 1}^\infty$ and $\{R_{r(k)}(k)\}_{k = 1}^\infty$ are, the more scale-heterogeneous the network is.
\end{definition}

Based on this definition, we see that the scale-heterogeneity is the main factor in the benefit of TVCS over TICS. In fact, scale-homogeneous and scale-heterogenous networks are the same as class $\Ic$ and $\Vc$ networks, respectively, due to~\eqref{eq:iota-r}. Further, note that the degree of scale-heterogeneity provides a \emph{geometric and qualitative} characterization of the amount of benefit TVCS has over TICS and distinguishes between networks in $\Vc$ that only marginally benefit from TVCS and those the benefit significantly (while $2k$-communicability is a more \emph{quantitative} notion used for computational assignment of networks to class $\Vc$ or $\Ic$).

It follows immediately from Definition~\ref{def:het} that determining the scale-heterogeneity of a network requires computation of all $\{r(k)\}_{k = 1}^\infty$ which is infeasible. Next, we seek simple and computationally efficient conditions to be used as a proxy for scale-heterogeneity.
}

\new{
\subsection{Identifying Class $\Vc$ Networks}\label{sec:V}
}

\new{In this section we discuss a sufficient condition for scale-heterogeneity that, when satisfied, ensures that a network belongs to class $\Vc$. This condition, given next, relies on the fact that $r(1)$ and $r(\infty)$ are particularly important elements of $\{r(k)\}_{k = 1}^\infty$ in determining scale-heterogeneity. The proof of this theorem is given in Appendix~\ref{sec:proofs}.}

\begin{theorem}\longthmtitle{Class $\Vc$ networks}\label{thm:main-suf}
  \new{Consider the TVCS problem~\eqref{eq:opt} for the network dynamics~\eqref{eq:dyn}.}
  Assume that the adjacency matrix $A$ is irreducible, aperiodic, and
  diagonalizable. If
  \begin{align*}
    \argmax_{i \in \Nc} R_i(1) \cap \argmax_{i \in \Nc}
    R_i(\infty) = \emptyset, 
  \end{align*}
  then the network belongs to class $\Vc$ for sufficiently large $K$. \oprocend
\end{theorem}

The condition of $A$ being irreducible is equivalent to the network
being strongly connected, and thus not restrictive. Likewise, $A$
being aperiodic is not restrictive as it requires that there exists no
integer number greater than 1 that divides the length of every cycle
in the network (satisfied, in particular, if any self-loops
exist). Finally, $A$ is almost always diagonalizable in the Lebesgue
sense, i.e., the set of non-diagonalizable $A$ has Lebesgue measure
zero.

Consider again the
networks of Figure~\ref{fig:tv-adv}(c and f). Here, the color
intensity of each node indicates its value $R_i(1)$ while its size
corresponds to its value $R_i(K - 1)$. Clearly, the first few largest
and darkest nodes are distinct in Figure~\ref{fig:tv-adv}(c), while
there is a close correlation between nodal size and darkness in
Figure~\ref{fig:tv-adv}(f), illustrating the root cause of their
difference in benefiting from TVCS.

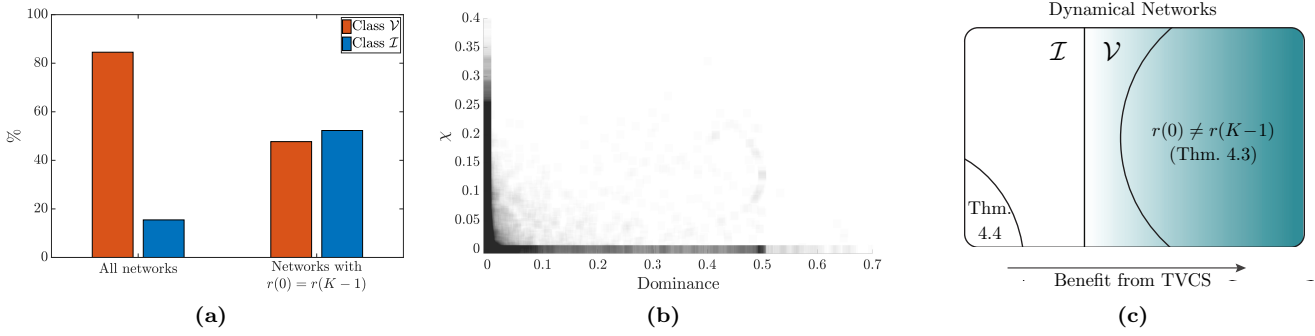
\begin{figure*}
  \begin{center} 
    \subfloat[]{
      \includegraphics[width=0.31\textwidth]{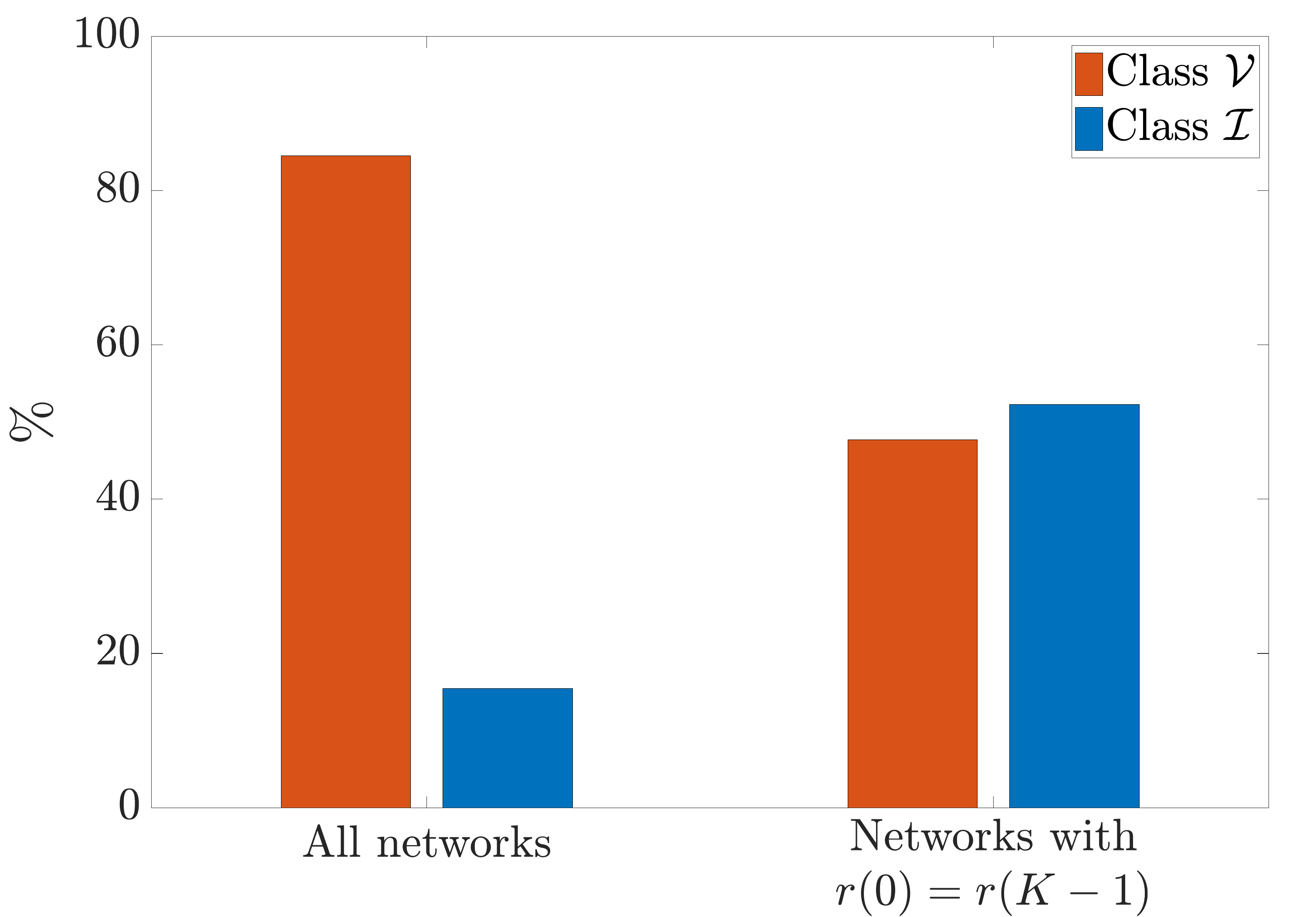}
    }
    \subfloat[]{
      \includegraphics[width=0.34\textwidth]{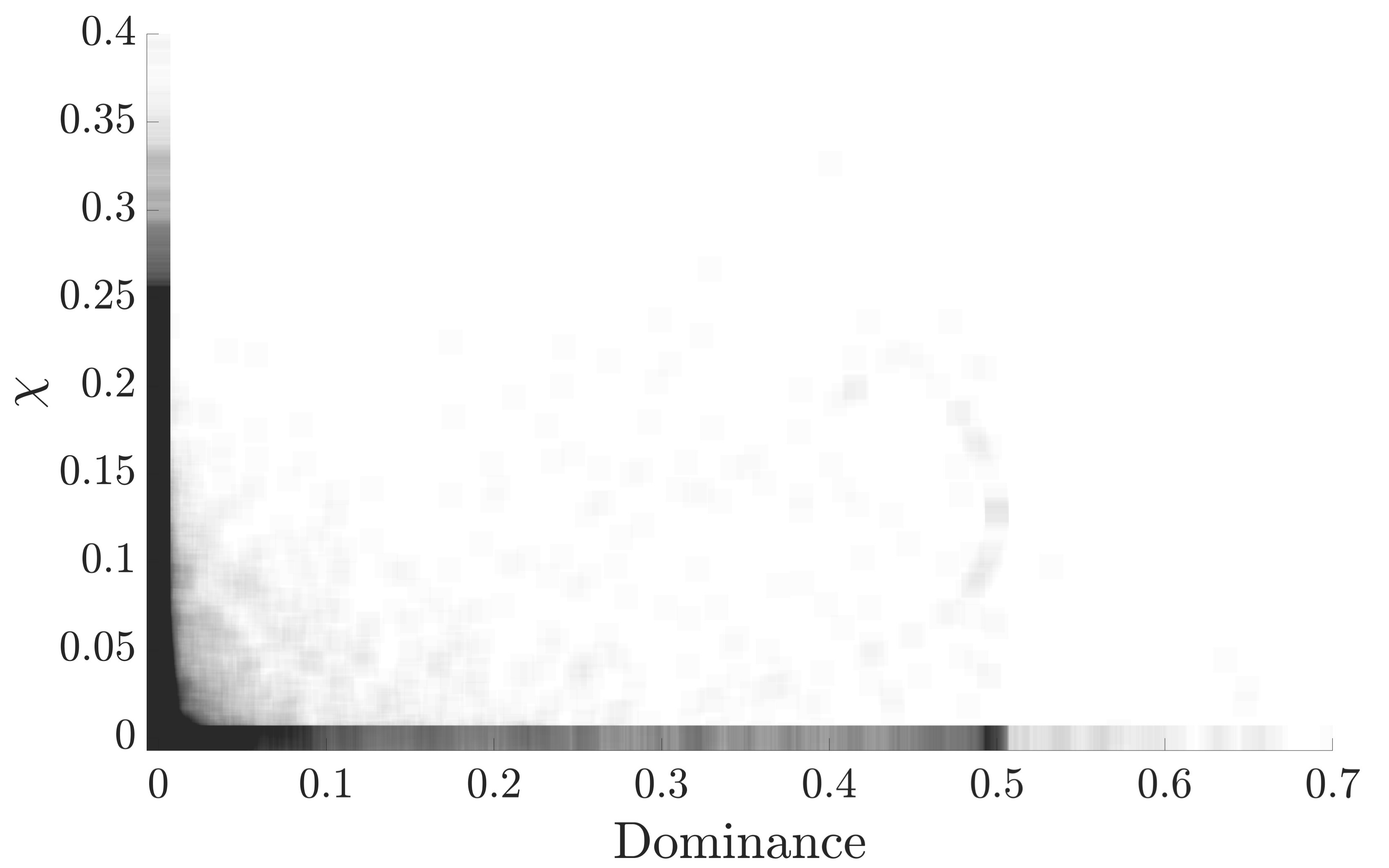}
    }
\hspace{10pt}
    \subfloat[]{
    \begin{tikzpicture}
    \node (venn) {\includegraphics[width=0.27\textwidth]{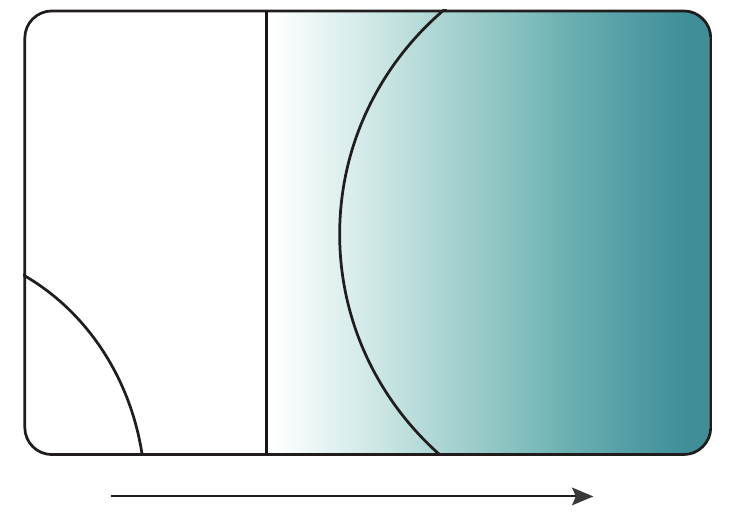}};
    \node[below of=venn, yshift=-20pt, scale=0.7] (benefit) {Benefit from TVCS};
    \node[below left of=venn, xshift=-34pt, yshift=-6pt, scale=0.7] (subI) {\parbox{30pt}{\centering Thm. \\ \ref{thm:main-nec}}};
    \node[right of=venn, xshift=2pt, yshift=3pt, scale=0.7] (subV) {\parbox{70pt}{\centering $r(0) \neq r(K \!-\! 1)$ \\ (Thm.~\ref{thm:main-suf})}};
    \node[above left of=venn, xshift=-8pt, yshift=18pt] (I) {$\Ic$};
    \node[above right of=venn, xshift=-28pt, yshift=18pt] (V) {$\Vc$};
    \node[above of=venn, yshift=25pt, scale=0.7] (nets) {Dynamical Networks};
    \end{tikzpicture} 
    }
  \end{center} 
  \caption{The role of $2k$-communicability in distinguishing between
    networks of class $\Vc$ ($\chi > 0$) and $\Ic$ ($\chi =
    0$). \textbf{(a)} The proportion of random networks in $\Vc$
      and $\Ic$. A total of $10^5$ random connectivity matrices were
      generated with logarithmically-uniform $n$ in $[10^1, 10^3]$,
      uniform sparsity $p$ in $[0, 1]$, and uniform pairwise
      connectivity weight in $[0, 1]$, and then transformed to
      adjacency matrices $A$ using the transmission method
      (cf. Appendix~\ref{sec:c2a}).  A time-horizon of $K = 10$ is
      used for all networks. While more than $80\%$ of all networks
      belong to class $\Vc$, this number drops to less than $50\%$
      among networks with $r(1) = r(K - 1)$ (i.e., networks where the
      same node has the greatest local and global
      centralities). \textbf{(b)} The $\chi$-value of the same
      networks as in (a) that have $r(1) = r(K - 1)$ as a function of
      the dominance of the node $r(0)$. For the node $r(0)$, its
      \emph{dominance} (over the rest of the network) is a measure of
      how distinctly $R_{r(0)}(1)$ and $R_{r(0)}(K - 1)$ are larger
      than $R_i(1)$ and $R_i(K - 1)$, respectively, for $i \neq r(0)$
      (cf. Appendix~\ref{sec:dom}). Each gray square represents
      one randomly generated network, so the darkness of each area
      represents the probability of observing random networks with
      that value of (dominance, $\chi$). A rapid decay of $\chi$ with
      dominance is clear, such that networks with positive dominance
      have very low probability of having $\chi > 0$. \textbf{(c)} A
    Venn diagram illustrating the decomposition of dynamical networks
    based on the extent to which they benefit from TVCS. The color
    gradient is a depiction of this extent, as measured by $\chi$
    (equation~\eqref{eq:chi}), where darker areas correspond to higher
    $\chi$. As shown in (a) and (b), the class of networks for which
    $r(0) \neq r(K - 1)$ is only a subset of $\Vc$ but provides a good
    approximation for it.}
  \label{fig:Km1eq1}
\end{figure*}

If a network has $r(0) = r(K - 1)$, it is still possible that the
network belongs to class $\Vc$. In fact, about half of the networks
with $r(0) = r(K - 1)$ still belong to $\Vc$
(Figure~\ref{fig:Km1eq1}(a)). However, these networks have a value of
$\chi$ of no more than $3\%$ on average, and in turn this value
quickly decreases with the dominance of the node $r(0)$ over the rest
of the network nodes (Figure~\ref{fig:Km1eq1}(b)). This is a
  strong indication that, for most practical purposes, the test based
  on $2k$-communicability is a valid indicator of whether a network
  benefits from TVCS. Furthermore, in the case of undirected
networks, it is possible to analytically prove that a network belongs
to class $\Ic$ ($\chi = 0$) if certain conditions based on the
eigen-decomposition of the adjacency matrix $A$ are satisfied, as shown next.

\new{
\subsection{Identifying Class $\Ic$ Networks}
}

\new{
Complementary to Section~\ref{sec:V}, here we discuss some necessary conditions for scale-heterogeneity based on the eigen-structure of the network that characterize subsets of $\Ic$.
} 
Let $A = V
\Lambda V^T$ be the eigen-decomposition of $A$, where $V = [v_{i
  j}]_{n \times n}$ 
and $\Lambda =
\diag(\lambda_1,\dots,\lambda_n)$ \new{is the diagonal matrix of eigenvalues} with $|\lambda_1| \ge |\lambda_2|
\ge \cdots \ge |\lambda_n|$. 
Further, let $W = [w_{i j}]_{n \times n}$ be the doubly stochastic matrix such that $w_{i j} = v_{i j}^2$ for
all $i, j \in \until{n}$.  
The next result, proven in Appendix~\ref{sec:proofs}, characterizes three undirected sub-classes of~$\Ic$.

\begin{theorem}\longthmtitle{Class $\Ic$
    networks}\label{thm:main-nec}
   \new{Consider the TVCS problem~\eqref{eq:opt} for the network dynamics~\eqref{eq:dyn}. Assume that the network is undirected (i.e., $A = A^T$) and that, without loss of generality, the node with the largest eigenvector centrality is labeled as node $1$.}
  If any of the following conditions holds:
  \begin{enumerate}
  \item $\frac{1 - w_{1 1}}{w_{1 1}} \le \frac{|\lambda_1| -
      |\lambda_2|}{|\lambda_1| - |\lambda_n|}$,
  \item $w_{1 1} + w_{1 2} = 1$,
  \item the network has three or fewer nonzero eigenvalues with
    different absolute values and $1 \in \argmax_i R_i(1)$,
  \end{enumerate}
  then, 
  \begin{align}\label{eq:1-in-argmax}
      1 \in \argmax_{1 \le i \le n} R_i(k), \quad \forall k \in \{0,
      \ldots, K - 1\},
    \end{align}
    i.e., selecting the node with the largest eigenvector centrality at
  every time step is the solution to~\eqref{eq:opt}. \oprocend
\end{theorem}

The conditions in Theorem~\ref{thm:main-nec} are based on the
eigen-decomposition of the network adjacency matrix $A$ and thus
abstract. However, these conditions can be interpreted as follows:
\begin{enumerate}
\item Condition~\emph{(i)} holds for networks where there is a
  sufficiently distinct central node, in the sense of eigenvector
  centrality, and the network dynamics is dominated by the largest
  eigenvalue.  An extreme case of such networks is a totally disconnected
  network where $W = I$ and the highest authority is the node with the
  largest self-loop.

\item Condition~\emph{(ii)} holds for networks where the eigenvector centrality of
  all nodes is determined by the weight of the link to the most
  eigenvector-central node. To see this, note that we have $w_{1 j} = 0$ for $j
  \ge 3$, implying $v_{1 j} = 0, j \ge 3$. Since the rows of $V$
  are orthogonal, we deduce $ v_{i 2} = \alpha v_{i 1}$ for all $i \ge
  2$, where $\alpha = -v_{1 1} / v_{1 2}$ is constant. Using $A = V
  \Lambda V^T$, we have
  \begin{align*}
    a_{1 i} = \lambda_1 v_{1 1} v_{i 1} + \lambda_2 v_{1 2} v_{i 2} =
    (v_{1 1} \lambda_1 + \alpha v_{1 2} \lambda_2) v_{i 1},
  \end{align*}
  so $v_{i 1} \propto a_{1 i}$ for all $i \ge 2$. Examples of such networks are star
  networks with no (or small-weight) self-loops (cf. Proposition~\ref{prop:star}).

\item Regarding condition~\emph{(iii)}, the most well-known families
  of networks with three distinct eigenvalues are the complete
  bipartite networks and connected strongly regular
  networks. Moreover, cones on $(n, k, \lambda, \mu)$-strongly regular
  graphs satisfying
  $\lambda_\text{min}(A) (\lambda_\text{min}(A) - k) = n$ are also
  known to have three distinct eigenvalues~\cite{ERV:96}.  The other
  condition $1 \in \argmax_i R_i(1)$ holds when the node with the
  largest eigenvector centrality (i.e., $r(\infty)$) has also the largest
  $2$-communicability. The
  simplest example of a network with these properties is the star
  network (with no or equal self-loops).
\end{enumerate}

The general abstraction from these cases is that a network belongs to
class $\Ic$ if it contains a sufficiently distinct central node, which
reinforces our main conclusion that $\Vc$ is the class of networks
with multiple scale-heterogeneous central nodes.
  The inclusion relationships between the various classes of networks introduced in this section are summarized in Figure~\ref{fig:Km1eq1}(c).

\new{
While Theorem~\ref{thm:main-nec} is only applicable to undirected networks, it has a straightforward extension to normal networks (i.e., directed networks with normal $A$). Using the same proof technique as in Theorem~\ref{thm:main-nec}, it can be shown that the exact same results hold if one replaces the eigenvalues and eigenvectors with singular values and singular vectors of $A$. Interestingly in this case, $R_i(\infty)$ coincides with HITS hub/authority centrality of node $i$ squared~\cite{JMK:99}.
}

\new{
\subsection{Networks with Latent Nodes}\label{sec:latent}
}

\new{
As mentioned in Section~\ref{subsec:model}, in many real-world applications of TVCS not all the
nodes are available/accessible for control. In this case, we call a node \emph{manifest} if it can be actuated and \emph{latent} if it cannot.
The natural solution would then be to choose the control nodes optimally among the
manifest nodes.
If the adjacency matrix $A$ of the network is fixed and given, this is the best solution. However, there are cases (often in man-made networks) where $A$ itself can be changed, at least among the manifest nodes. We call such a change of structure a \emph{manipulation}. If manipulation is possible in a network with latent nodes, another solution to TVCS is to manipulate the network among the manifest nodes such that the optimal control
nodes (when computed without any restrictions on control scheduling)
lie among the manifest nodes for all time.
The following results provides a guarantee that this is always effective, 
provided that the manipulation is sufficiently strong and not acyclic.
}

\begin{theorem}\longthmtitle{Network manipulation and TVCS in networks with latent nodes}\label{thm:sn}
  Consider the optimal node selection problem~\eqref{eq:opt} over
  a time horizon~$K$. Given a network of $n$ nodes with adjacency
  matrix $A_0 \in \real^{n \times n}$, let $E \in \real^{n \times n}$
  be a nonnegative matrix of the form \vspace{-6pt}
  \begin{center}
  \begin{tikzpicture}
  \node[draw=none] at (0, 0) (matrix) {
  \parbox{100pt}{
  \begin{align*}
  \arraycolsep=8pt
  \def\arraystretch{1.4}
  \def\bracketHeight{18pt}
  E = 
  \left[\rule{0cm}{\bracketHeight}\right.
  \hspace{-5pt}
  \text{\raisebox{1pt}{$\begin{array}{cc} \star & 0 \\ 0 & 0 \end{array}$}}
  \hspace{-5pt}
  \left]\rule{0cm}{\bracketHeight}\right.
  \end{align*}
  }
  };
  \node at (37pt, 8pt) (right-top1) {$\big\}$};
  \node[right of=right-top1, xshift=-21pt] (right-top2) {$\scriptstyle n_1$};
  \node at (37pt, -9pt) (right-bottom1) {$\big\}$};
  \node[right of=right-bottom1, xshift=-15pt] (right-bottom2) {$\scriptstyle n - n_1$};
  \node[rotate=270] at (0pt, 18pt) (top-left1) {$\big\{$};
  \node[above of=top-left1, yshift=-23pt] (top-left2) {$\scriptstyle n_1$};
  \node[rotate=270] at (21pt, 18pt) (top-right1) {$\big\{$};
  \node[above of=top-right1, yshift=-22pt] (top-right2) {$\scriptstyle n - n_1$};
  \node[right of=matrix, xshift=35pt] (comma) {$,$};
  \end{tikzpicture}
  \end{center}
  \vskip -13pt
  \noindent corresponding to the manifest subnetwork involving the first $n_1 < n$ nodes
  (this is without loss of generality, since nodes can be renumbered) and
  consider the dynamic network
  described by~\eqref{eq:dyn} with adjacency matrix $A = A_0 + \alpha
  E$, where $\alpha > 0$.
    Then, if $E$ is not acyclic, there exists $\overline \alpha > 0$
  such that for $\alpha > \overline \alpha$,
  \begin{align}\label{eq:man}
    r(k) \in \until{n_1},
  \end{align}
  for all $k \in \{0, \dots, K - 1\}$. Furthermore, if $A_0$ and $E$ are symmetric (the corresponding networks are undirected), $\bar \alpha$ can be found in closed form and~\eqref{eq:man} holds for all $k \ge 1$. \oprocend
\end{theorem}

Both requirements of Theorem~\ref{thm:sn} (that $\alpha E$ is
sufficiently strong and acyclic) have clear interpretations. First,
depending on how large the size of the manifest subnetwork is and how
central its nodes already are (pre-manipulation), larger manipulation
may be necessary to turn them into central nodes at various scales
(i.e., $r(k)$ for $k = \{0, \dots, K - 1\}$). Second, for the manifest
nodes to become central at arbitrarily global scales (i.e., $r(k)$ for
$k \sim K \to \infty$), the manipulation must contain paths of
arbitrarily long lengths, which are absent in acyclic networks.

According to Theorem~\ref{thm:sn}, manipulation of the manifest
subnetwork is effective even when the manifest nodes are among the
least central nodes of the network (before the manipulation). In this
case, as we increase $\alpha$ from $0$, the manifest nodes usually
first turn into the most locally-central nodes ($\alpha \not\ge \bar
\alpha$ yet), and then also into globally-central nodes ($\alpha >
\bar \alpha$). The following example illustrates this phenomenon in a
simple star network where the center node is latent and the peripheral
nodes are manifest.

\begin{example}\longthmtitle{Undirected star networks with varying self-loop
    weights}\label{ex:star}
   Consider an undirected uniform star network given by
   \begin{align*}
    A_0 = \left[\begin{array}{cc} l_p I_{n - 1} & a_{c p} \ones_{n - 1} \\ a_{c p} \ones_{n - 1}^T & l_c \end{array}\right],
  \end{align*}
  \new{where $\ones_{n - 1}$ denotes the $(n - 1)$-dimensional vector
    of all ones and} the positive constants $l_c$, $l_p$, and $a_{c
    p}$ are the central self-loop weight, peripheral self-loop weight,
  and the link weight between the center node and any peripheral node,
  respectively. The $2k$-communicabilities of this network are
  computed analytically in Proposition~\ref{prop:star} (nodes are
  re-labeled here for conformity with Theorem~\ref{thm:sn}). It
  follows from~\eqref{eq:rik-star} that for any $i \in \{1, \dots, n -
  1\}$,
      \begin{align}\label{eq:star-diff1} 
       R_n(1) - R_i(1) = l_c^2 - l_p^2 + (n - 2) a_{c p}^2.
    \end{align}
    Therefore, if $l_p \le
    l_c$, then $R_n(k) > R_i(k)$ for all $k \ge 1$,
    i.e., the center node is the optimal control node at all
      times.  However, when $l_c < l_p$, the network can exhibit
    different behaviors. From~\eqref{eq:rik-star}, we can also see that 
    \begin{align}\label{eq:star-diff-inf}
    \lim_{k \to \infty} R_n(k) > \lim_{k \to \infty} R_i(k) \Leftrightarrow \lambda_1 - l_p > a_{c p}.
    \end{align}
    Define $\underline l_p = \sqrt{l_c^2 + (n - 2) a_{c
          p}^2}$ and $\overline l_p = l_c + (n - 2)
      a_{c p}$. Using~\eqref{eq:star-diff1}-\eqref{eq:star-diff-inf} and after some
    computations, one can see that 
    \begin{center}
      \begin{tabular}{ll}
        $r(k) = n$ for all $k$, & 
        if $l_p \le \underline l_p$,
        \\[1.5ex]
        $r(1) = \until{n - 1}$ but $r(k) = n$ for large enough $k$, & if
        $\underline l_p < l_p < \overline l_p$, 
        \\[1.5ex]
        $r(k) =
        \until{n - 1}$ for all $k$, &  if $l_p \ge \overline l_p$.
    \end{tabular}      
  \end{center}
  In other words, when the manipulation is weak, the (latent) center
  node is the optimal control node at all times. As the manipulation
  gains strength, scale-heterogeneity emerges, making the (manifest)
  peripheral nodes the optimal control node at local scales while the
  center node remains still the optimal control node at global
  scales. Finally, when the manipulation is strong enough,
  scale-heterogeneity vanishes, leaving the (manifest) peripheral
  nodes as the optimal control nodes at all scales. Notice that with
  the terminology of Theorem~\ref{thm:sn},
   \begin{align*}
     E = \left[\begin{array}{cc} I_{n - 1} & 0 \\ 0 &
         0 \end{array}\right], \quad n_1 = n - 1, \quad \alpha = l_p,
     \quad \text{and} \quad \overline \alpha = \overline
     l_p. \eqoprocend
   \end{align*}
\end{example}

\begin{figure}
  \begin{center} 
    \includegraphics[width=0.45\textwidth]{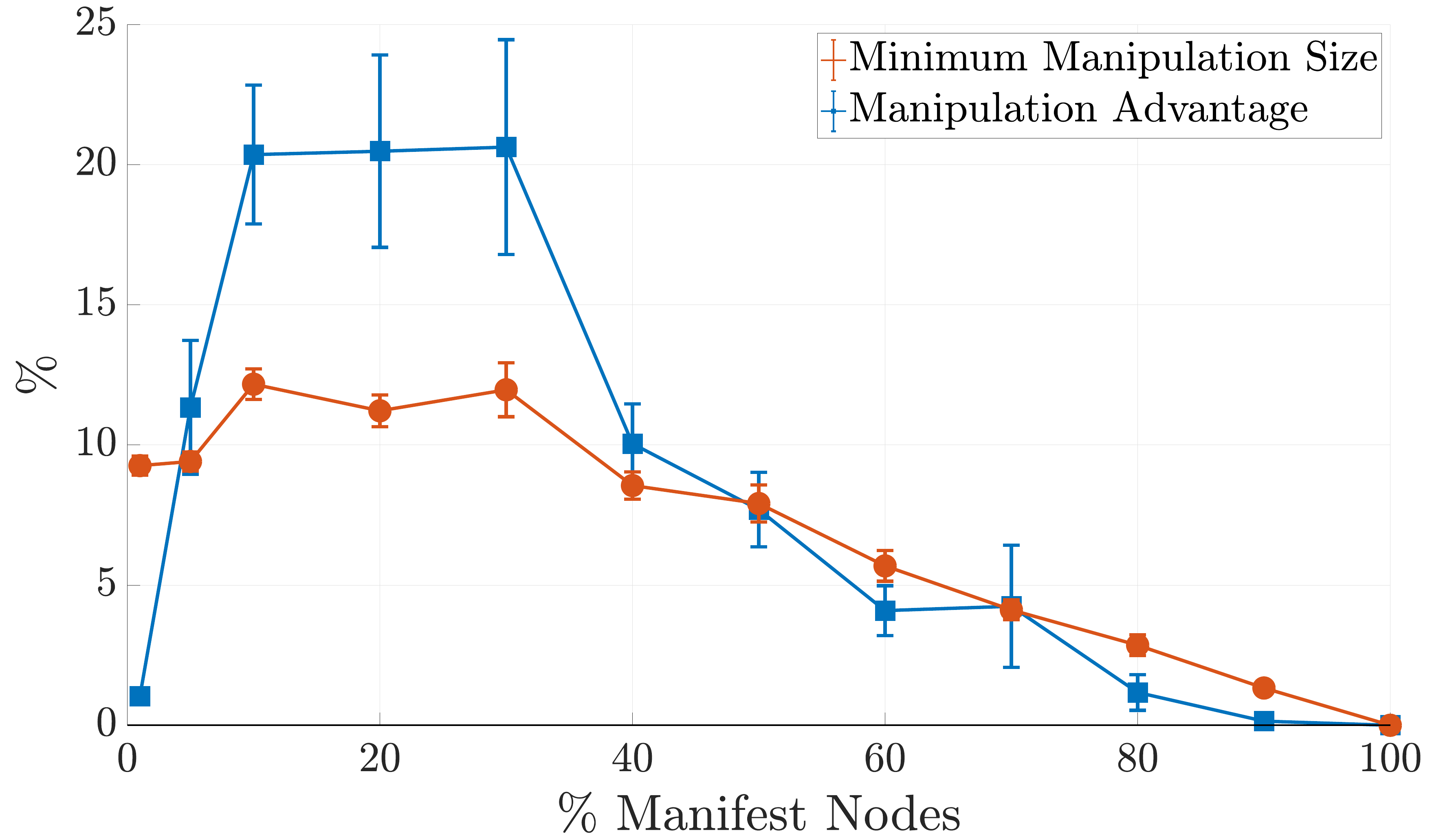}
    \caption{Manipulation of manifest subnetworks in order to obtain
      an all-manifest optimal TVCS. The horizontal axis represents the
      percentage of manifest nodes in the network. In red, we show the
      minimum size of manipulation needed for the optimal TVCS to only
      include manifest nodes, relative to the size of the initial
      adjacency matrix (both measured by induced matrix $2$-norm). In
      blue, we depict the optimal (i.e., maximal) value of
      $\tr(\Wc_K)$ for the case where the minimal manifest
      manipulation is applied, relative to the maximal value of
      $\tr(\Wc_K)$ subject to the constraint that all the control
      nodes are manifest (the former is with manipulation and without
      constraints on the control nodes, while the latter has no
      manipulation but control node constraints). Results are for
      $10^3$ random networks of logarithmically-uniform sizes in
      $[10^1, 10^3]$ but otherwise similar to
      Figure~\ref{fig:Km1eq1}. Markers (circles/squares) represent
      average values and error bars represent standard error of the
      mean (s.e.m). In both cases, the overall adjacency matrix is
      normalized by its spectral radius for fairness of comparison. We
      see that medium-sized manifest subnetworks ($5 \sim 20\%$) are
      the hardest yet most fruitful to manipulate.}
    \label{fig:man}
  \end{center} 
\end{figure}

A fair concern, however, exists regarding the minimum size of the
manipulation needed to make the TVCS all-manifest. If this is
excessively high, the prescribed approach may be infeasible in
practice. Nevertheless, among networks of various size and structure,
random manipulations with norm of about $10\%$ of the norm of $A$ are
on average sufficient (Figure~\ref{fig:man}). Here, we see that the
largest manipulations are needed for manifest subnetworks of about
$10\%$ the total size of the network. This is because when the size of
the manifest subnetwork is extremely small, manipulations are focused
on this small subset of nodes and thus more efficient, while with
extremely large manifest subnetworks, the majority of the nodes are
accessible for control and there is little restriction on the TVCS.

Finally, Figure~\ref{fig:man} also shows the comparison, in terms of
controllability, of the manipulation-based approach against the
alternative approach of selecting an optimal TVCS with the additional
constraint that control nodes must be manifest (without any
manipulation of the dynamics), which results in a sub-optimal
all-manifest TVCS. For the comparison to be fair, we normalize each
network by its spectral radius (largest magnitude of its eigenvalues),
and then compare the optimal value of their TVCS
(equation~\eqref{eq:opt}). We see that the amount of relative
advantage produced by manifest subnetwork manipulation is comparable
to the relative size of the manipulation, except for medium-sized
manifest subnetworks ($5 \sim 20\%$ of nodes), where the manipulation
advantage is about two times its size.

\new{
\section{Case Study: TVCS in Synthetic and Real Networks}
}

\begin{figure}[tb!]
\begin{center} 
  \tikzset{filled node/.style={draw, fill, text width=5pt, inner
      sep=0pt, circle}}
  \subfloat[]{
    \centering
    \def\nodesep{15pt}
    \begin{tikzpicture}[every node/.style=filled node]
      \node at (0, 0) (1) {};
      \node at (\nodesep, 0) (2) {}
      edge[very thick] (1);
      \node[draw=none, fill=none] at (\nodesep/2, 7pt) (a2) {$a$};
      \node at (2*\nodesep, 0) (3) {}
      edge[very thick] (2);
      \node[draw=none, fill=none] at (3*\nodesep/2, 7pt) (a3) {$a$};
      \node[draw=red!70!black, fill=red!70!black] at (3*\nodesep, 0) (4) {}
      edge[very thick] (3);
      \node[draw=none, fill=none] at (5*\nodesep/2, 7pt) (a4) {$a$};
      \node at (4*\nodesep, 0) (5) {}
      edge[very thick] (4);
      \node[draw=none, fill=none] at (7*\nodesep/2, 7pt) (a5) {$a$};
      \node at (5*\nodesep, 0) (6) {}
      edge[very thick] (5);
      \node[draw=none, fill=none] at (9*\nodesep/2, 7pt) (a6) {$a$};
      \node at (6*\nodesep, 0) (7) {}
      edge[very thick] (6);
      \node[draw=none, fill=none] at (11*\nodesep/2, 7pt) (a7) {$a$};
      \node[draw=none, fill=none] at (0, -25pt) (ghost) {};
    \end{tikzpicture}
  }
  \hspace{15pt}
  \subfloat[]{
    \centering
    \begin{tikzpicture}[every node/.style=filled node]
      \def \n {9}
      \def \radius {25pt}
      \def \lradius {30pt}
      \foreach \s in {1,...,\n}
      {
        \node at ({360/\n * (\s - 1)}:\radius) {};
        \draw[very thick] ({360/\n * (\s - 1)}:\radius) 
        arc ({360/\n * (\s - 1)}:{360/\n * (\s)}:\radius);
        \node[draw=none, fill=none] at ({360/\n * (\s - 1) + 360/2/\n}:\lradius) {$a$};
      }
      \node[draw=red!70!black, fill=red!70!black] at ({360/\n*2}:\radius) {}; 
    \end{tikzpicture}
  }
  \hspace{15pt}
  \subfloat[]{
    \centering
    \begin{tikzpicture}[every node/.style=filled node]
      \def \n {7}
      \def \radius {25pt}
      \def \sradius {17pt}
      \def \angleshift {15pt}
      \node[draw=red!70!black, fill=red!70!black] at (0, 0) (c) {};
      \foreach \s in {1,...,\n}
      {
        \node at ({360/\n * (\s - 1)}:\radius) {}
	edge[very thick] (c);
        \node[draw=none, fill=none] at ({360/\n * (\s - 1) + \angleshift}:\sradius) {$a$};
      }
      \node[draw=none, fill=none] at (0, -30pt) (ghost) {};
    \end{tikzpicture}
  }
  \caption{Simple networks with closed-form
    $2k$-communicabilities. \textbf{(a)} A line network, \textbf{(b)} a ring network,
    and \textbf{(c)} a star network. All networks are undirected and have
    homogeneous edge weights $a$. The $2k$-communicabilities of these
    networks are analytically computed (cf. Appendix~\ref{sec:simple-networks}),
    concluding that all networks belong to class $\Ic$, with the
    optimal control node depicted in red in each case (the optimal
    control node is arbitrary in a ring network due to its
    symmetry).}\label{fig:line-ring-star}
\end{center} 
\end{figure}
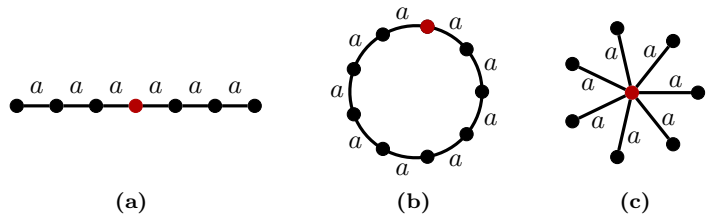

Here, we discuss the benefits of TVCS and its relation to network
structure for several examples of synthetic and real networks.
We start with the classical deterministic examples of undirected line,
ring, and star networks (Figure~\ref{fig:line-ring-star}). Due to
their simple structure, the $2k$-communicabilities of these networks
can be analytically computed in closed form
(cf. Appendix~\ref{sec:simple-networks}). Using these results,
it follows that for the line and star networks, the optimal control
node is always the center node (or any of the two center nodes if a
line has even number of nodes), while the optimal control node is
arbitrary in a ring network. Notice that in all cases, it is the
\emph{homogeneity} of these networks that results in a single node
having the greatest centrality at all scales
(cf. Example~\ref{ex:star} for non-homogeneous star networks
that have scale-heterogeneous central nodes and thus belong to class
$\Vc$).

\begin{figure*}
\subfloat[]{
\includegraphics[height=0.19\textheight]{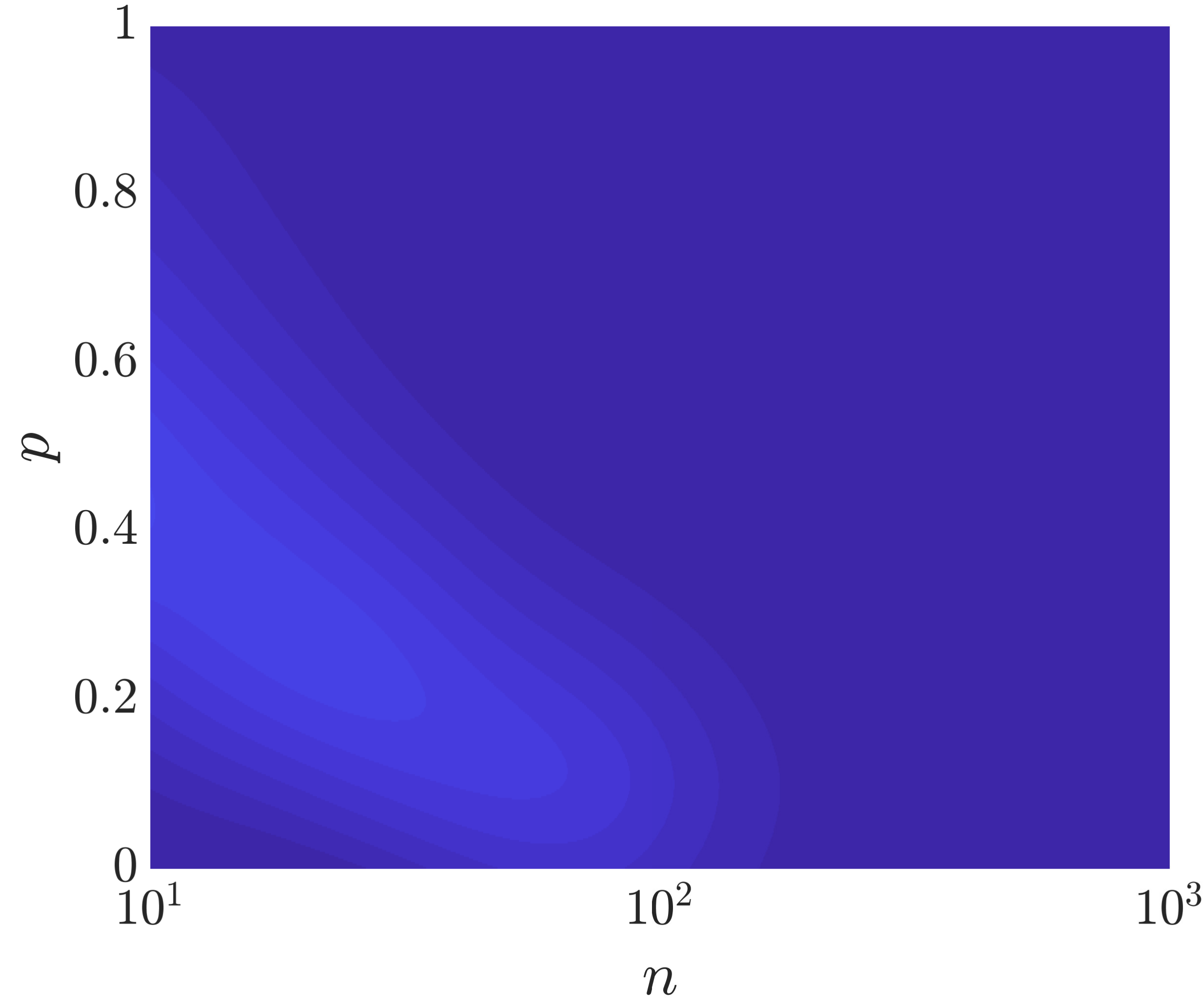}
}
\subfloat[]{
\includegraphics[height=0.19\textheight]{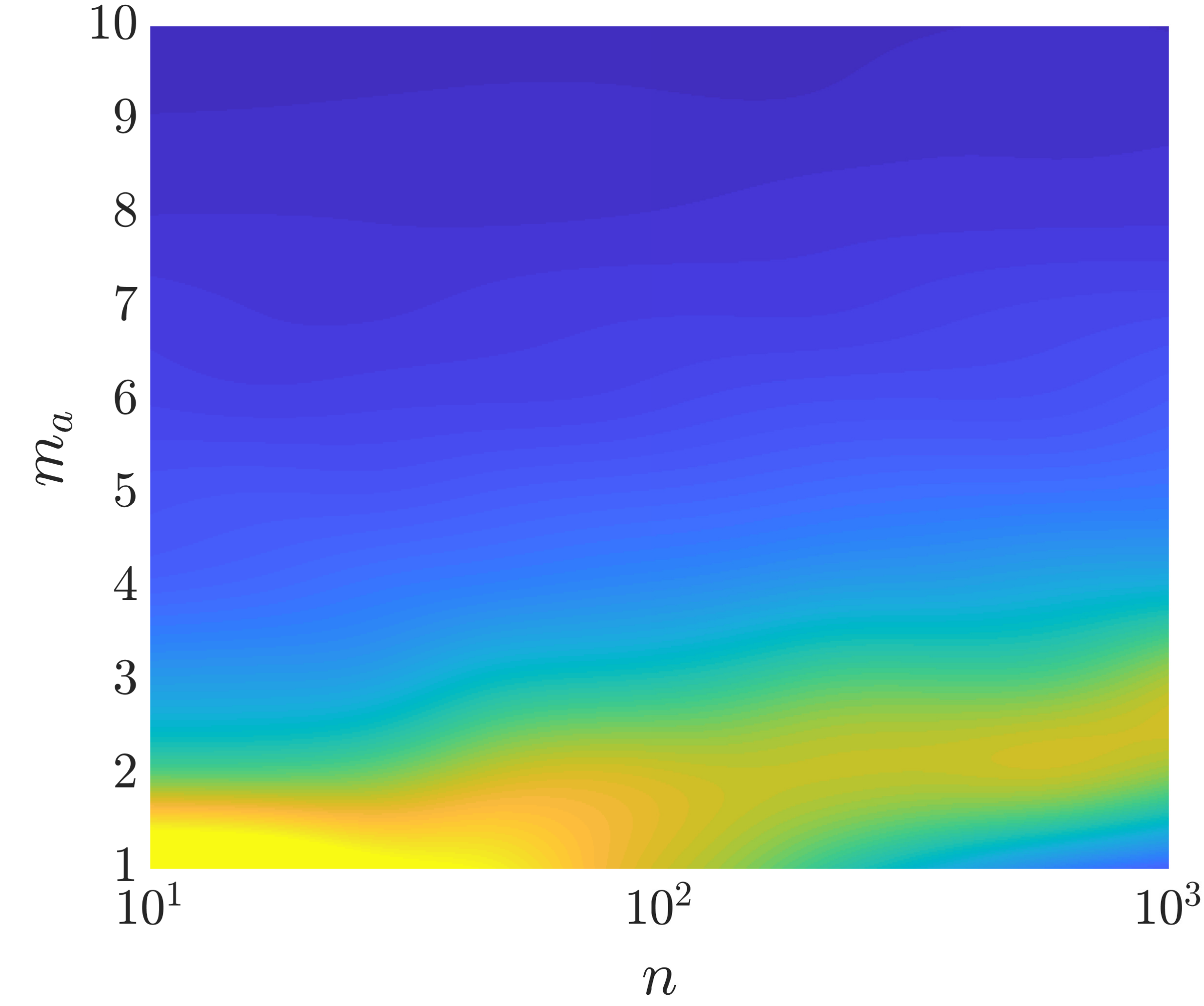}
}
\subfloat[]{
\includegraphics[height=0.19\textheight]{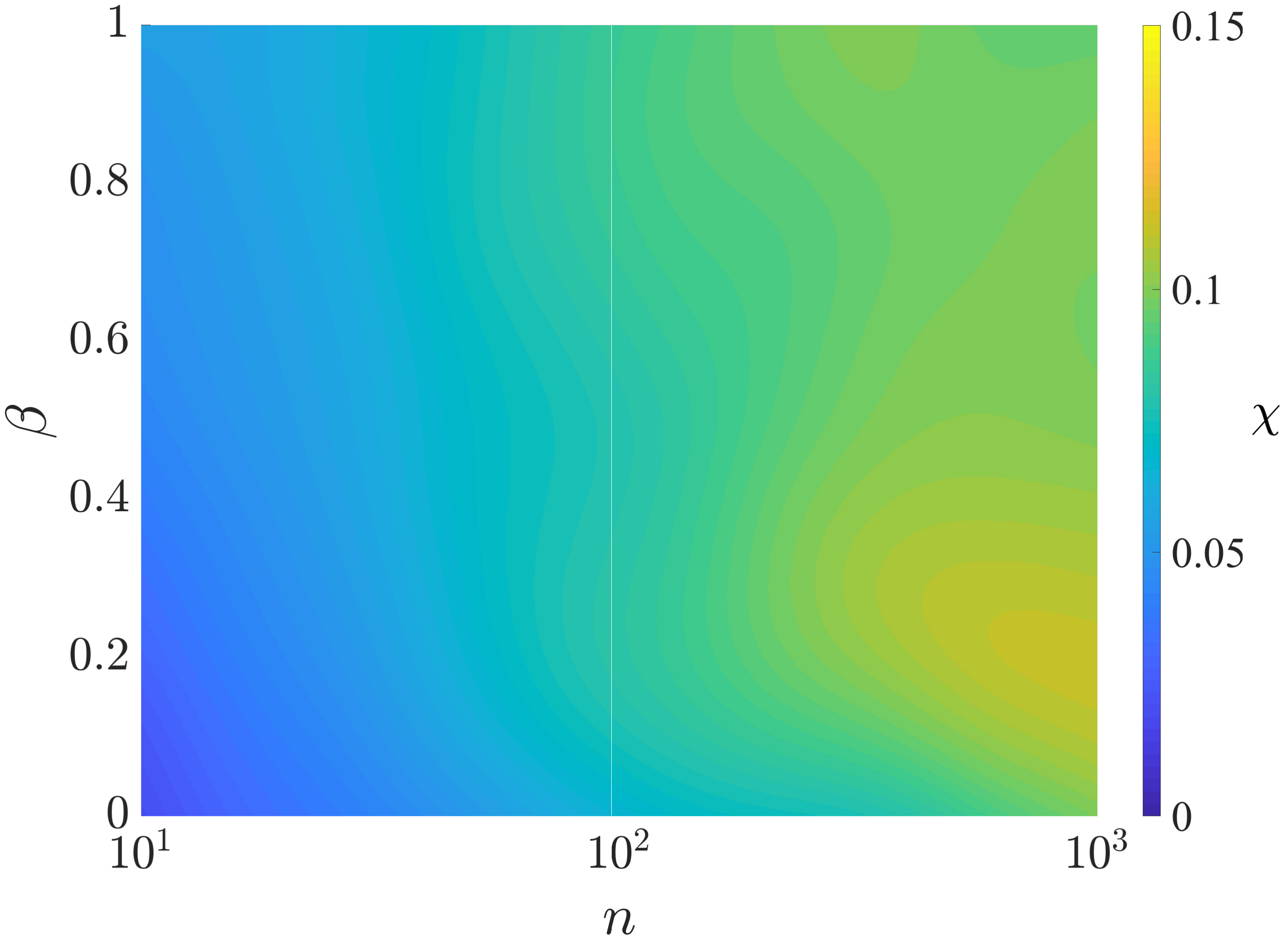}
}
\caption{The average $\chi$-value for \textbf{(a)} ER, \textbf{(b)}
  BA, and \textbf{(c)} WS probabilistic networks. The horizontal axis
  determines the size of the network $n$ in all cases, while the
  vertical axis determines the values of the corresponding parameters
  for each network: edge probability $p$ for ER, growth (link
  attachment) rate $m_a$ for BA, and rewiring probability $\beta$ for
  WS. After constructing the unweighted connectivity according to each
  algorithm (ER, BA, or WS), standard uniformly random weights are
  assigned to each edge, which is then converted to $A$ using
  transmission method (cf. Appendix~\ref{sec:c2a}). For each value of
  $n$ and network parameter over a coarse mesh ($\sim 100$ points),
  $100$ networks are generated and the average of their $\chi$-value
  is computed, which is then smoothly interpolated over a fine mesh
  (MATLAB \texttt{csaps}).}
\label{fig:rand}
\end{figure*}

Next, we analyze the role of TVCS in three classes of probabilistic
complex networks that are widely used to capture the behavior of
various dynamical networks. These include the Erd\"os-R\'enyi (ER)
random networks, Barab\'asi-Albert (BA) scale-free networks, and
Watts-Strogatz (WS) small-world networks. Each network has its own
characteristic properties, and these properties lead to different
behaviors under TVCS. The average $\chi$-values of these networks are
computed for various values of $n$ and network parameters
(Figure~\ref{fig:rand}). For ER networks, $\chi$ is in general small,
and decays with $n$. This is because ER networks, especially when $n$
is large, are extremely homogeneous. This homogeneity is further
increased during the transmission method, leading to a network matrix
$A$ that is extremely insensitive to the choice of control nodes.

The connectivity structure of BA networks, in contrast, is extremely
inhomogeneous, with one (sometimes 2) highly central nodes and a
hierarchy down to peripheral leafs. As one would expect, this implies
a small $\chi$-value since the center node has the highest centrality
at all scales (Figure~\ref{sfig:rand-all}). However, when the
connectivity matrix is transformed to $A$ using the transmission
method, the incoming links to all nodes are made uniform (adding up to
$1$). This in turns make the centrality levels of all the nodes
comparable, leading to high $\chi$-values observed (notice that the
underlying connectivity structures are still highly inhomogeneous,
distinguishing them from the homogeneous ER networks). Notice that as
the growth rate $m_a$ is increased, smaller networks tend towards
complete graphs and high $\chi$ values \emph{shift} to larger $n$.

As our last class of probabilistic networks, WS networks have the
broadest range of size-parameter values with significant $\chi$. As
one would expect, $\chi$ is low near $\beta = 0, 1$, corresponding to
regular ring lattice and ER networks, respectively. For $\beta \sim
0.2$, there is a sufficiently high probability of having multiple
nodes that are close to many rewired links (increasing their
centrality), yet there is a low probability that these nodes, and the
nodes close to them, are rewired all alike, resulting in heterogeneous
central nodes and high $\chi$-values. This heterogeneity is increased
with $n$ as larger networks have more possibilities of rewiring every
edge.

Finally, we used the tools and concepts introduced so far to analyze
TVCS in several real-world dynamical networks
(Table~\ref{tab:real}). These networks are chosen from a wide range of
application domains, from neuronal networks to transportation and
social networks. According to the type of dynamics evolving over each
network, we have used either the transmission or induction method to
obtain its dynamical adjacency matrix from its static connectivity
(\new{the ``$C \to A$" column,} cf. Appendix~\ref{sec:c2a}). 

We have computed the $\chi$-value for each network using a variable
time horizon $K \le 50$, with the results ranging from 0 to more than
$30\%$ for different networks. These large variations even within each
category signify both the potential benefits of TVCS and the
possibility of its redundancy, a contrast that has been pivotal to our
discussion.  \new{In particular, four facts about these results worth
  highlighting. (i) As measured by $\tr(\Wc_K)$, the majority of
  networks tested do not benefit from TVCS, but a few do so
  significantly. (ii) Despite coming from various domains, the
  networks that do significantly benefit from TVCS share
  \emph{scale-heterogeneity} as their common qualitative property
  (cf. Section~\ref{sec:2k-comm}). (iii) Networks with inductive $C
  \to A$ transformation benefit significantly less from TVCS than
  those with transmission $C \to A$ transformation. (iv) Significantly
  higher values of $\chi$ are expected for all networks if using
  $\lambda_{\min}(\Wc_K)$ or similar measures for controllability,
  cf. Appendix~\ref{sec:measures}.}

In the last column, we have also indicated whether the most local and
most global central nodes coincide in each network. Recall that this
is a sufficient but not necessary condition for a network to be in
class $\Vc$ (Theorem~\ref{thm:main-suf} and
Figure~\ref{fig:Km1eq1}). Though only sufficient, this simple metric
can correctly classify class members of $\Vc$ from $\Ic$ among these
networks, except for the WesternUS power network, for which $r(0) =
r(K - 1)$ only marginally holds (the dominance of $r(0)$ is $0$)
(cf. Figure~\ref{fig:Km1eq1}(b)).

\begin{table*}
\addtolength\tabcolsep{0pt}
\centering
\caption{Characteristics of the real-world networks studied in the paper}
\label{tab:real}
\resizebox{\textwidth}{!}{%
\begin{tabular}{llcccccccc}
\hline
Category & Name & $n$ & $|E|$ & Directed & $C \to A$ & $\chi (\%)$ & \parbox{5em}{\center\vspace*{-5pt}$r(0) = r(K - 1)$\vspace*{5pt}} & \parbox{6em}{\center\vspace*{-5pt}Dominance of \\ $r(0)$ ($\times 10^{-3}$) \vspace*{5pt}} & ref.
\\
\hline
Neuronal & 	BCTNet fMRI & 	638 & 	37250 & 	N & 	T & 	$1.8$ & 	N & 	N/A &	\cite{MR-OS:10}
\\
& 			Cocomac & 		58 & 		1078 & 	Y & 	T & 	$5.5$ & 	N & 	N/A & 	\cite{RB-TW-MD:12}
\\
& 			BCTNet Cat & 		95 & 		2126 & 	Y & 	T & 	$1.9$ & 	N & 	N/A & 	\cite{MR-OS:10}
\\
& 			C. elegans & 		306 & 	2345 & 	Y & 	T & 	$0$ & 	Y & 	$0$ & 	\cite{DJW-SHS:98}
\\[10pt]
Transportation & air500 & 			500 & 	5960 & 	N & 	T & 	$22.4$ & 	N & 	N/A & 	\cite{VC-RPS-AV:07}
\\
& 			airUS & 			1858 & 	28236 & 	Y & 	T & 	$0$ & 	Y & 	$0$ & 	\cite{TO:11}
\\
& 			airGlobal & 		7976 & 	30501 & 	Y & 	T & 	$0$ & 	Y & 	$0$ & 	\cite{TO:11}
\\
& 			Chicago & 		1467 & 	2596 & 	N & 	T & 	$0$ & 	Y & 	$0$ & \cite{RWE-KSC-YJL-DEB:83,DEB-KSC-MEF-YJL-KTL-RWE:85}
\\[10pt]
Gene Regulatory & E. coli & 		4053 & 	127544 & 	N & 	T & 	$0$ & 	Y & 	$0$ & 	\cite{HK-JES-JS-IL:15}
\\[10pt]
PPI & 		Yeast & 			2361 & 	13828 & 	N & 	T & 	$0$ & 	Y & 	$0$ & 	\cite{DB-YZ-LC-HX-XZ-HL-JZ-SS-LL-NZ-GL-RC:03}
\\
& 			Stelzl & 			1706 & 	6207 & 	Y & 	T & 	$0$ & 	Y & 	$0$ & 	\cite{US-UW-ML-CH-FHB-HG-MS-MZ-AS-SK-JT-SM-CA-NB-SK-AG-ET-AD-SK-BK-WB-HL-EEW:05}
\\
& 			Figeys & 			2239 & 	6452 & 	Y & 	T & 	$0$ & 	Y & 	$0$ & 	\cite{RME-PC-FE-HL-PT-SC-LM-MDR-LO-ML-RT-MD-YH-AH-LM-SZ-OO-YVB-ME-YS-JV-MA-JPPL-HSD-IIS-BK-KH-KC-KG-BM-RK-SLLA-MFM-GBM-TT-DF:07}
\\
& 			Vidal & 			3133 & 	12875 & 	N & 	T & 	$0$ & 	Y & 	$0$ & 	\cite{JR-KV-TH-THK-AD-NL-GFB-FDG-MD-NAG:05}
\\[10pt]
Power & 		WesternUS & 		4941 & 	13188 & 	N & 	T & 	$33.7$ & 	Y & 	$0$ & 	\cite{DJW-SHS:98}
\\[10pt]
Food & 		Florida & 			128 & 	2106 & 	Y & 	T & 	$34.6$ &	 N & 	N/A & 	\cite{REU-JJH-MSE:00}
\\
& 			LRL & 			183 & 	2494 & 	Y & 	T & 	$27.3$ & 	N & 	N/A & 	\cite{NDM-JJM-TK-MS:91}
\\[10pt]
Social & 		Facebook group & 	4039 & 	176468 & 	N & 	I & 	$0.4$ & 	N & 	N/A & 	\cite{JL-JJM:12}
\\
& 			E-mail & 			1005 & 	25571 & 	Y & 	I & 	$0$ & 	Y & 	$40.5$ & 	\cite{HY-ARB-JL-DFG:17,JL-JK-CF:07}
\\
& 			Southern Women & 	18 & 		278 & 	N & 	I & 	$0$ & 	Y & 	$1.6$ & 	\cite{AD-BBG-MRG:41}
\\
& 			UCI P2P & 		1899 & 	20296 & 	Y & 	I & 	$0$ & 	Y & 	$5.5$ & 	\cite{TO-PP:09}
\\
& 			UCI Forum & 		899 & 	142760 & 	N & 	I & 	$0$& 	Y & 	$2.8$ & 	\cite{TO:13}
\\
& 			Freeman's EIES & 	48 & 		830 & 	Y & 	I & 	$0$ & 	Y & 	$1.4$ & 	\cite{SCF-LCF:79}
\\
& 			Dolphins & 		62 & 		318 & 	N & 	I & 	$0$ & 	Y & 	$0.7$ & 	\cite{DL-KS-OJB-PH-ES-SMD:03}
\\[10pt]
Trust & 		Physicians & 		241 & 	1098 & 	Y & 	I & 	$8.8$ & 	N & 	N/A & 	\cite{JC-EK-HM:57}
\\
& 			Org. Consult Advice & 46 & 	879 & 	Y & 	I & 	$0$ & 	Y & 	$0.1$ & 	\cite{RLC-AP:04}
\\
& 			Org. Consult Value & 46 & 	858 & 	Y & 	I & 	$0$ & 	Y & 	$1.2$ & 	\cite{RLC-AP:04}
\\
& 			Org. R\&D Advice & 	77 & 		2228 & 	Y & 	I & 	$6 \times 10^{-3}$ & N & N/A & \cite{RLC-AP:04}
\\
& 			Org. R\&D Aware & 	77 & 		2326 & 	Y & 	I & 	$0$ & 	Y & 	$0.3$ & 	\cite{RLC-AP:04}
\\
\hline
\end{tabular}%
}
\vspace{5pt}
\caption*{For each network, we have reported the number of nodes $n$, number of edges $|\Ec|$ (with each bidirectional edge counted twice), whether the network is directed, the method used for obtaining dynamical adjacency matrix $A$ from static connectivity $C$ ($A \to C$), the $\chi$ value (equation~\eqref{eq:chi}), and whether the most local and global central nodes coincide ($r(0) = r(K - 1)$). Since the value of $\chi$ is a function of $K$, we have chosen the value of $K \le 50$ that has the largest $\chi$ for each network. Detailed descriptions of these datasets are provided in Appendix~\ref{sec:real}.}
\end{table*}

\section{Conclusions and Discussion}

Despite the breadth and depth of existing literature on the
controllability of complex networks and control scheduling, the
significant potential of TVCS has been greatly overlooked. This work
strives to explore the advantages of TVCS in linear dynamical networks
and obtaining theoretical and computational relationships between
these advantages and network structure. Using Gramian-based measures
of controllability, we showed that TVCS can significantly enhance the
controllability of \emph{many but not all} synthetic and real
networks. This motivated the pursuit of identifying properties based
on network structure that explain when, why, and by how much TVCS is
beneficial.

Using the newly introduced notion of $2k$-communicability, we showed
that the scale-heterogeneity of central nodes in a network is the main
cause and correlate of TVCS advantages. If a network has several
distinct central nodes at different scales, the optimal TVCS involves
starting the control from the most global central nodes and gradually
moving towards most local ones as the time horizon is approached. If,
on the other hand, a single node acquires the highest centrality at
all scales, optimal TVCS prescribes the sole control of this node
over the entire horizon, leading to optimality of TICS.

A striking finding that defied our expectations is the effect of
network dynamics, beyond its raw connectivity structure, on
TVCS. Here, we differentiated between the raw connectivity structure
of a network (obtained using specific field knowledge and measure the
\emph{relative} strength of nodal connections) and its dynamical
adjacency matrix which determines the evolution of network state over
time. Depending on the nature of network state, we proposed two
methods, transmission and induction, for obtaining the dynamical
adjacency matrix from static connectivity. The effects of these
methods, however, is noteworthy on the benefits of TVCS, even though
the underlying network connectivity is the same (Table~\ref{tab:real}
and Figure~\ref{sfig:rand-all}).  While the transmission method
significantly enhances the merit of TVCS, the induction method
depresses it (both compared to raw connectivity). We believe the
reason for the former is the additional \emph{homogeneity} that the
transmission method introduces among the nodes, while the latter is
due to the conversion from continuous to discrete-time dynamics, which
enables long-distance connections even over small sampling times (due
to the fact that interactions occur over infinitesimal intervals in
continuous time) \new{(cf. Section~\ref{sec:c2a} and
  Figure~\ref{fig:gamma})}. These results suggest that controllability
of network dynamics is not only a function of its structural
connectivity, but also greatly relies on the type of dynamics evolving
over the network, an aspect that has received little attention in the
existing literature and warrants future research.

\new{Our discussion so far applies to networks with and without
  self-loops alike. However, it follows from the results in
  Section~\ref{sec:res} that self-loops play an important role in
  TVCS. This is because (i) the self-loop of each node directly adds
  to its $2k$-communicability for all $k$, and (ii) the self-loop of
  each node also contributes indirectly to the $2k$-communicability of
  its neighbors less than $k-1$ hops away.  As a result, the self-loop
  of any node has the largest effect on its own $2k$-communicability
  for all $k$, but also a lesser effect on the $2k$ communicability of
  all other nodes in the network. This latter effect becomes smaller
  and limited to higher $k$ for more distant nodes. A clear
  demonstration of the effects of self-loops can be seen in
  Example~\ref{ex:star}, where as the self-loops of the peripheral
  nodes get stronger, they gradually become the central nodes in the
  network, first at local scales (small $k$) and eventually at all
  scales.}

Further, the focus of our discussion has so far been on single input networks
where one node is controlled at a time, in order to enhance the
simplicity and clarity of concepts. Nevertheless, our results have
straightforward generalizations to multiple-input networks
(cf. Appendix~\ref{sec:mi}). If $m$ denotes the number of control
inputs, the optimal TVCS involves applying these control inputs to the
$m$ nodes with the highest centralities at the appropriate scale at
every time instance (i.e., the $m$ nodes with the largest $R_i(K - 1 -
k)$ have to be controlled at every time instance $k$). It is clear
that the additional flexibility due to the additional inputs makes
$\Vc$ larger, i.e., more networks have $\chi > 0$. Nevertheless, this
additional flexibility also makes TICS significantly more
efficient. Therefore, it is not immediately clear whether this
enlargement of $\Vc$ also entails larger $\chi$ for networks with the
same size and sparsity. In fact, increasing $m$ reduces average $\chi$
for all the classes of ER, BA, and WS networks (Figure~\ref{sfig:mi}),
suggesting that the additional flexibility is more advantageous for
TICS than TVCS.

Regardless of the number of inputs ($1$ or more), an important
implicit assumption of TVCS is that this number is limited, i.e., no
more than $m$ nodes can be controlled at every time instance. This may
at first seem over-conservative since TVCS requires, by its essence,
the installation of actuators at all (or many) nodes of the
network. Therefore, one might wonder why limit the control to only $m$
nodes at every time instance when all the nodes are ready for
actuation. The answer lies within the practical limitations of
actuators. For ideal actuators, distributing the control energy over
as many nodes as possible is indeed optimal. However, this is not
  possible in many scenarios, including when (i) actuators exhibit
  nonlinear \emph{dead-zone} behaviors, so that each one requires a
  sizable activation energy. In many applications ranging from
  distributed industrial processes to opinion dynamics in social
  networks, nodes cannot be actuated with arbitrarily small amounts of
  control energy. If $E_{\min}$ is the minimum activation energy of
  any actuator, at least $m E_{\min}$ is required for actuation of $m$
  nodes at a time. Thus, when $E_{\min}$ is sizable and $n$ is large,
  simultaneous actuation of all nodes ($m = n$) requires a significant
  amount of control energy which is often infeasible (notice that the
  dead-zone behavior of actuators does not violate the linearity
  assumption in~\eqref{eq:dyn} as one can replace $u$ with $v =
  \phi(u)$, where $\phi$ denotes dead-zone nonlinearity); (ii)
  actuators are geographically disperse so that precise coordination
  becomes difficult or time-consuming. A familiar example of this is
  the social opinion dynamics in pre-election times during political
  campaigns, where rallies and speeches by candidates act as control
  inputs to the network. Even though all nodes may be actuatable, at
  most one node can be actuated at every time; (iii) simultaneous
  control of proximal nodes results in actuator interference. This is
  the case in many biological networks. In neuronal networks, for
  instance, common control technologies such as TMS do not allow for
  simultaneous actuation of all cortical areas due, in part, to
  electromagnetic interference between multiple sources of actuation
  (note that TVCS is still possible by installation and sequential
  activation of multiple coils at different locations); and when (iv)
  actuators are controlled via communication channels with limited
  capacity, so that only a small number of devices can be
  simultaneously operated. This may be the case in industrial
  applications where large numbers of geographically distributed
  actuators are remotely (and centrally) controlled over shared
  communication channels with limited bandwidth. In all these
  scenarios, TVCS has the potential to significantly enhance network
  controllability, conditioned on the scale-heterogeneity of the
  central nodes in the network.

  Although the dynamics of all real networks have some degrees of
  nonlinearity, the analysis of linear(ized) dynamics is a standard
  first step in analysis of dynamical properties of complex
  networks~\cite{YYL-JJS-ALB:11,NJC-EJC-DAV-JSF-CTB:12,GY-JR-YL-CL-BL:12,FP-SZ-FB:14,AO:14,THS-JL:14,YZ-FP-JC:16-cdc,VT-MAR-GJP-AJ:16,THS-FLC-JL:16,LZ-WZ-JH-AA-CJT:14,STJ-SLS:15,DH-JW-HZ-LS:17}. This
  is mainly due to the fact that stability and controllability of
  linearized dynamics of a nonlinear network implies the same
  properties \emph{locally} for the original nonlinear dynamics,
  making linear dynamics a powerful tool in analyzing many dynamical
  properties that are in general intractable for nonlinear
  dynamics. The local validity of linearization, however, is a main
  limitation of this work, particularly in networks where the change
  of state is significant relative to the size of the domain over
  which the linearization is valid.  \new{For these networks, whether
    the nonlinearity enhances or decreases the benefits TVCS with
    respect to its linearization is in general dependent on the type
    of nonlinearity. However, for saturation nonlinearities, being
    perhaps the most widespread, we expect TVCS to be more beneficial
    than linear counterparts. This is because in TICS all the control
    input is injected through a fixed node, requiring the state of
    that node to potentially undergo large over- and undershoots in
    order to convey sufficient input to the rest of the
    network. Saturation clearly prevents this from happening, further
    limiting the scope of TICS.  The generalization of this work to
    nonlinear dynamics with saturation and linear \emph{time-varying}
    dynamics (namely, $A(k)$ instead of $A$ in
    equation~\eqref{eq:dyn}) is a warranted next step for future
    exploration of the role of TVCS in general nonlinear networks.}

\appendices
\bigskip
\begin{center} 
{\LARGE \bfseries Appendix}
\end{center} 

\section{Obtaining Dynamical Adjacency Matrix from Static Connectivity}\label{sec:c2a}

\begin{figure}
\centering
\includegraphics[width=0.4\linewidth]{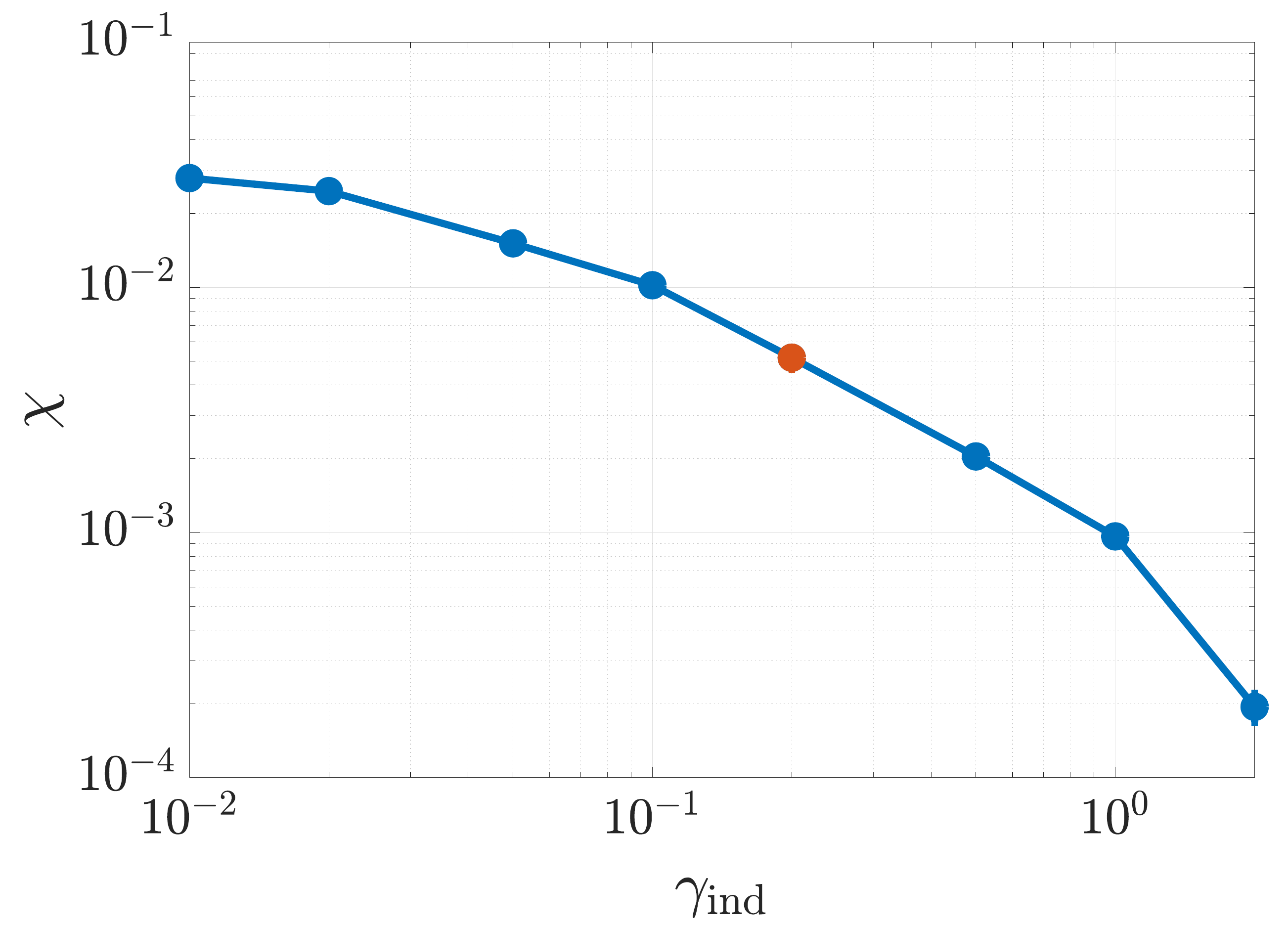}
\caption{The average value of $\chi$ for the induction method and varying values of $\gamma_{\rm ind}$ (corresponding to varying discretization step sizes $T_s$). Each point represents the average value of $\chi$ for $50$ realizations of ER networks with $n = 100$ and $p = 0.2$ and vertical bars (when visible) show one standard error of the mean (s.e.m.). For each network, the value of $K \le 10^3$ that gives the largest value of $\chi$ is chosen. The average value of $\chi$ drops with $\gamma_{\rm ind}$, showing the effect of discretization on $\chi$ and the merit of TVCS. The red point corresponds to $\gamma_{\rm ind} = 0.2$ used throughout this work.}
\label{fig:gamma}
\end{figure}

A standard starting point for the analysis of network dynamics of the form~\eqref{eq:dyn}
is the assumption that the network adjacency matrix $A$ is known. While this is a valid assumption (as the construction of $A$ is itself the subject of vast research in network identification and corresponding field sciences), care should be taken in how one interprets raw network connectivity matrices. Usually, the network structure is described not by its dynamic adjacency matrix $A$ (which determines the evolution of network \emph{state} according to~\eqref{eq:dyn}) but rather by its static connectivity matrix $C$
(our implicit assumption is that each node has a well-defined state that evolves over time through network dynamics, so our discussion is not applicable to completely static networks).
While for any $i, j \in \Nc$, $a_{ij}$ describes the impact of $x_j$ on $x_i$ over one time step (relative to $x_j$), $c_{ij}$ often describes the strength of the link $(i, j)$ in arbitrary units (e.g., number of synapses between two neurons, capacity of high-voltage lines between two generators, or number of seats on a flight). In particular, multiplying $C$ by a positive constant results in an equivalent description of the network structure, yet multiplying $A$ by a constant significantly alters network dynamics. Here, we outline two methods for obtaining $A$ from $C$, and describe example domains where each method seems more relevant. Consider an arbitrary link $(i, j) \in \Ec$.
\begin{itemize}
\item \textbf{Transmission:} This method applies to dynamical networks
  where at each time step, the value of the state of node $i$ is
  itself affected (reduced) as a result of interaction with neighbor
  node $j$. Here, the state of each node corresponds to a physical
  quantity that is \emph{transmitted} to its neighbors in order to
  affect their states. Neuronal, transportation, food, gene
  regulatory, protein-protein interaction, and power networks are all
  examples of this type of interaction. If the sampling time is chosen
  long enough such that ``current" state of a node is completely
  diffused through the network until the next time step, we can obtain
  $A$ from $C$ using
\begin{align*}
A = C D_{C, \text{in}}^{-1},
\end{align*}
where $D_{C, \text{in}}$ is the \emph{augmented in-degree matrix} of
$C$ (a diagonal matrix with the sum of the columns of $C$ on its
diagonal, except where the sum of a column of $C$ is zero, in which
case the corresponding diagonal element of $D_{C, \text{in}}$ is
$1$). This means that over each time step, $x_i$ is transmitted to the
in-neighbors of node $i$ proportionally to their connectivity strength,
if $i$ has any in-neighbors, and preserved otherwise.
\item \textbf{Induction:} This method is appropriate for networks in
  which nodal states are not physical quantities and thus do not
  reduce as a result of network interactions. Opinion or epidemic
  dynamics evolving over social and/or trust networks have such
  properties. Here, in order to compute $A$ from $C$, we start from
  the underlying continuous-time dynamics $\dot x = (-\alpha I + C) x$
  where $\alpha > 0$ is chosen such that $-\alpha I + C$ is stable
  (Hurwitz), and then discretize it to obtain~\eqref{eq:dyn}, where
\begin{align*}
A = e^{(-\alpha I + C) T_s},
\end{align*}
and $T_s$ is the sampling time~\cite[eq. (4.17)]{CTC:98}. From the
expansion of matrix exponential ($e^M = I + M + \frac{M^2}{2} +
\frac{M^3}{3!} + \cdots$), we see that $A$ does not inherit the
sparsity pattern of $C$ (and $G$) since nodes interact in continuous
time. However, if $\|(-\alpha I + C)^2 T_s^2 / 2\| \ll \|(-\alpha I +
C) T_s\|$, then the sparsity pattern of $C$ is almost preserved in
$A$. \new{Therefore, in this work we use $T_s = \gamma_{\rm ind} / \|\alpha I + C\|$ for the induction
method with $\gamma_{\rm ind} = 0.2$ unless otherwise stated.
Further, Figure~\ref{fig:gamma} shows the effect of $\gamma_{\rm ind}$ on the value of $\chi$ when using the induction method. As expected, the larger $\gamma_{\rm ind}$, the larger $T_s$, the closer $A$ gets to $\lim_{k \to \infty} A^k$, the more similar $2k$-communicabilities for different $k$ become, and the smaller $\chi$ becomes.}
\end{itemize} 

\begin{figure}
  \begin{center} 
    \setlength{\tabcolsep}{2pt}
    \begin{tabular}{cccc}
    &
      {\tiny Erd\"os-R\'enyi} & {\tiny
      Barab\'asi-Albert} &  {\tiny Watts-Strogatz}
      \\
      \rotatebox{90}{\hspace*{2ex} \tiny Transmission Method} 
      &
      \includegraphics[height=0.13\textheight]{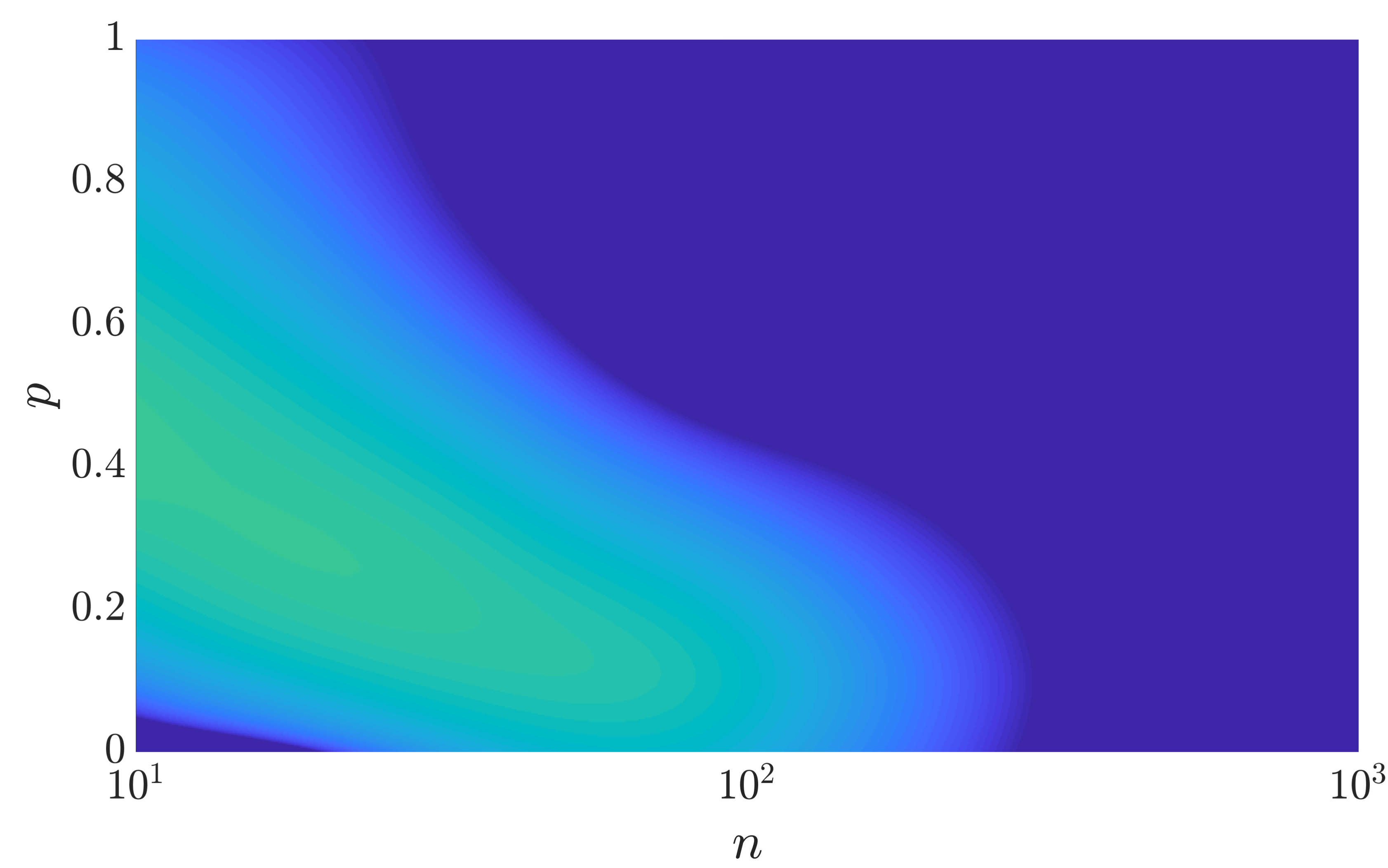}
      &
      \includegraphics[height=0.13\textheight]{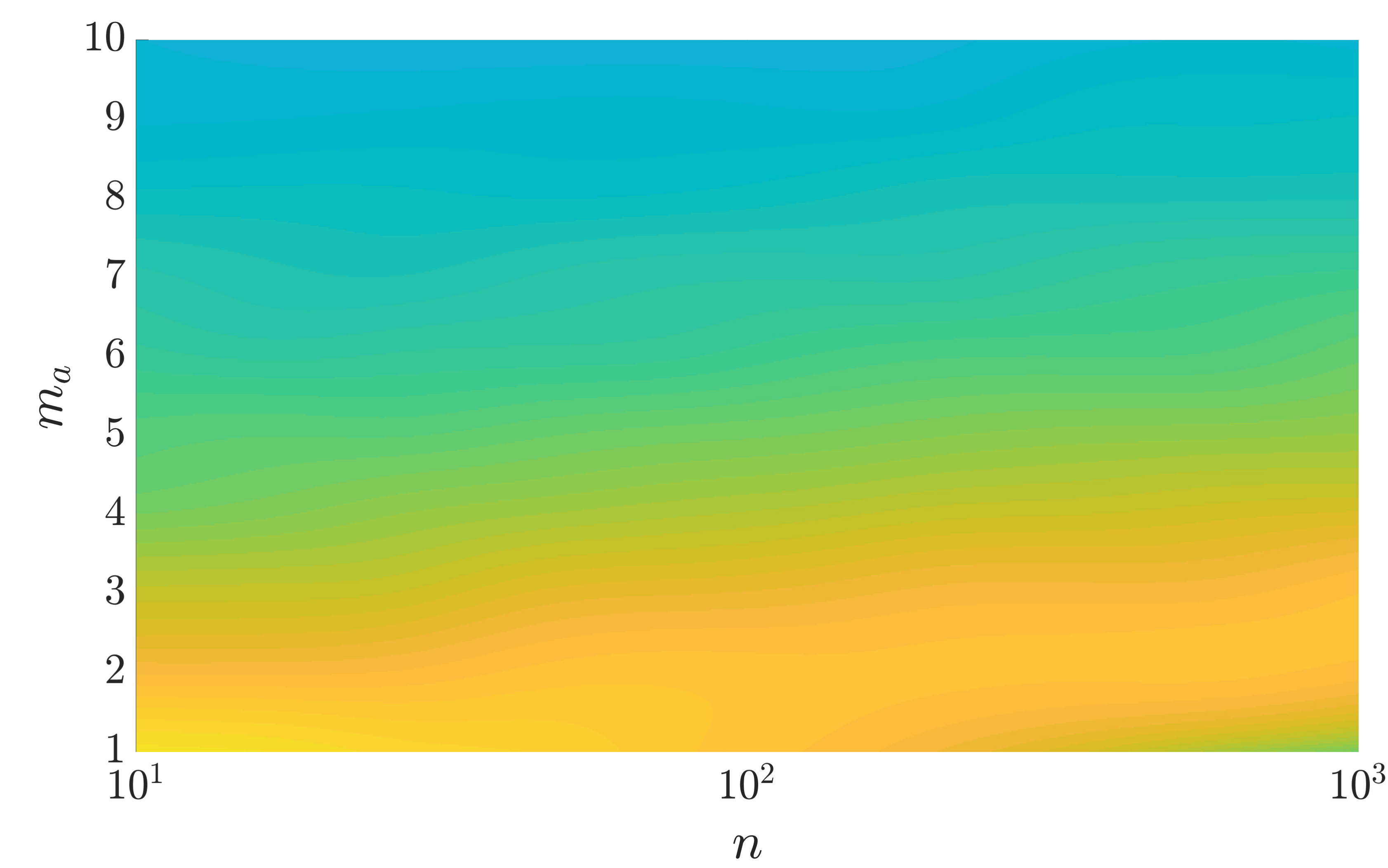}
      &
      \includegraphics[height=0.13\textheight]{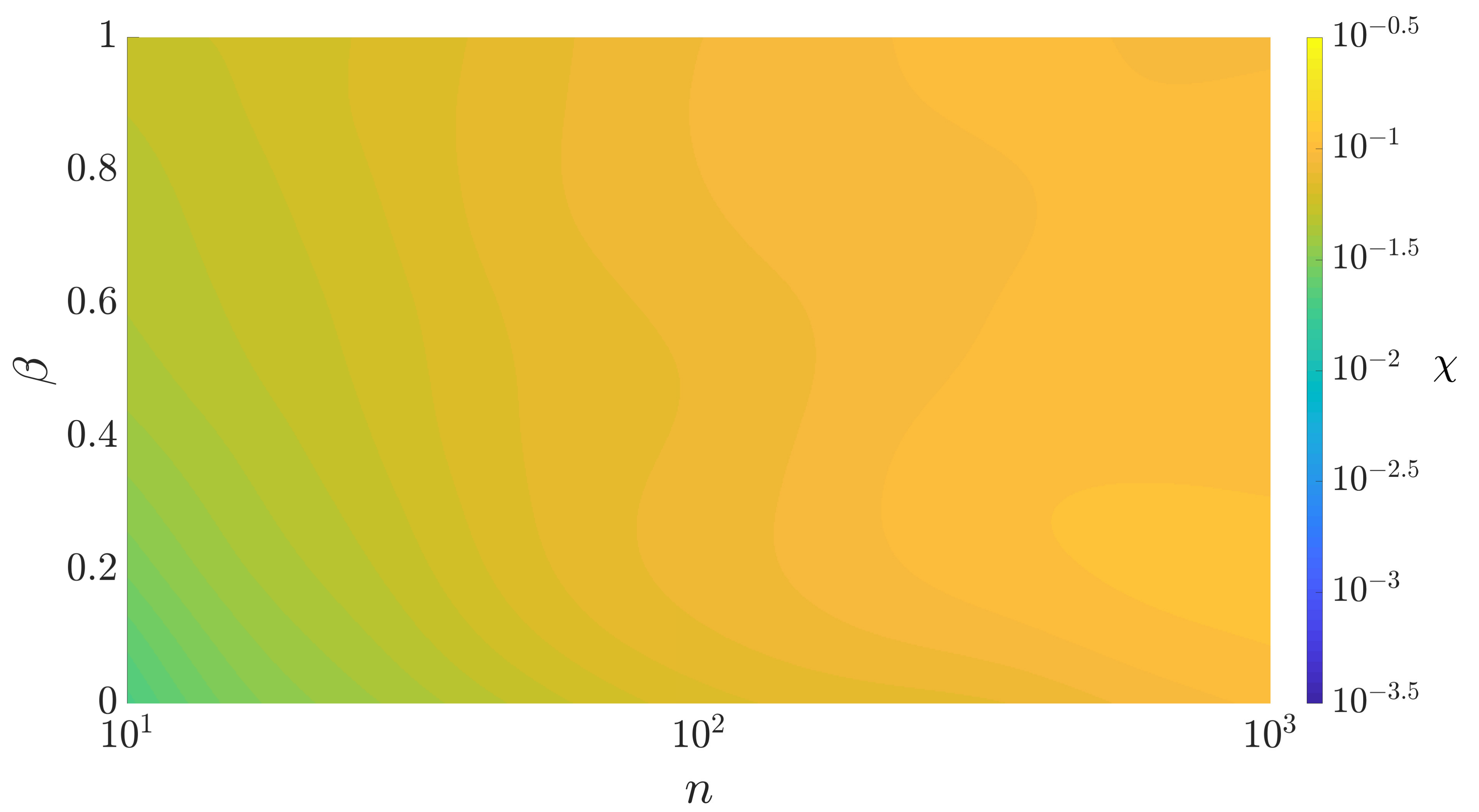}
      \\
      \rotatebox{90}{\hspace*{2ex} \tiny Induction Method}
      &
      \includegraphics[height=0.13\textheight]{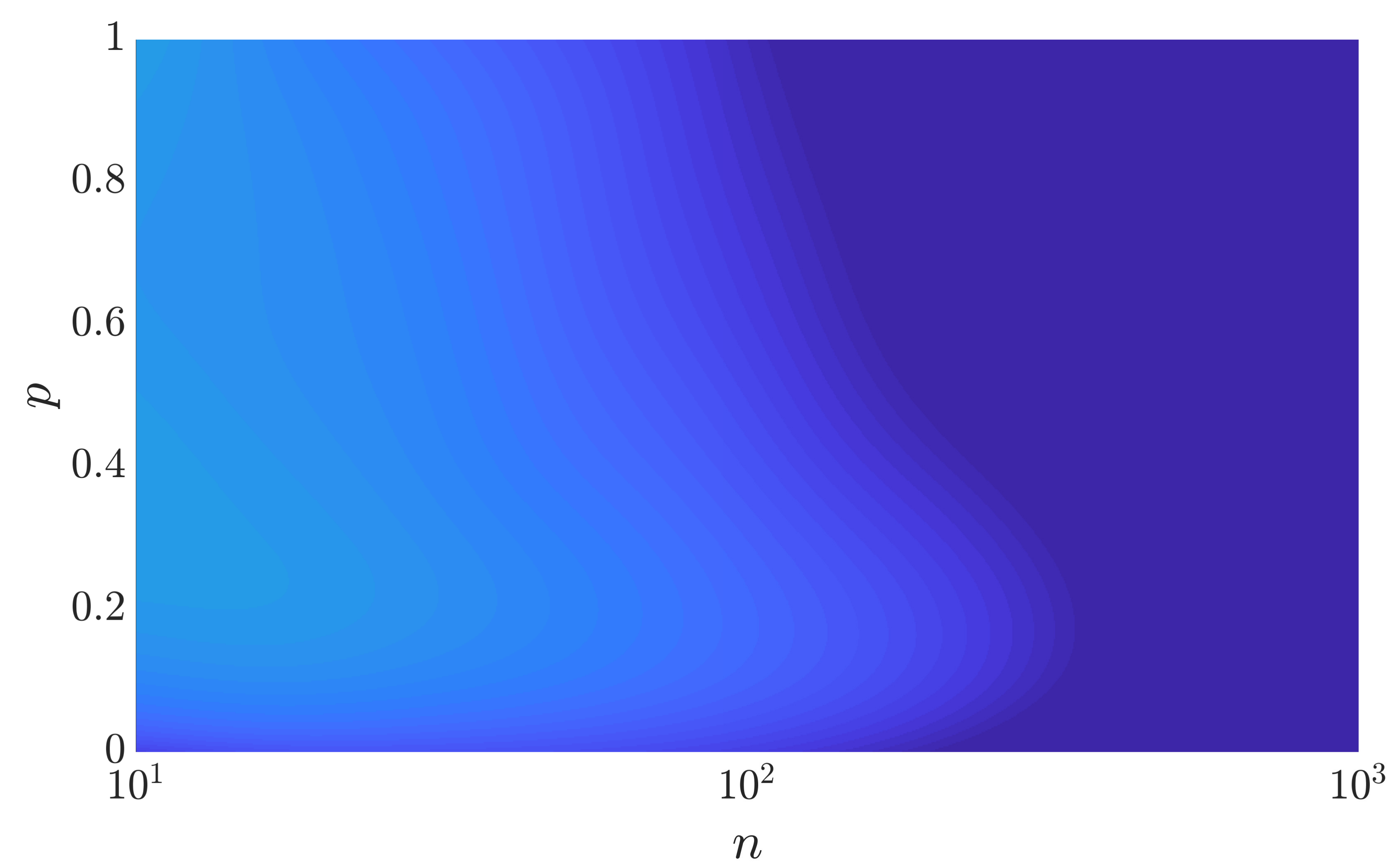}
      &
      \includegraphics[height=0.13\textheight]{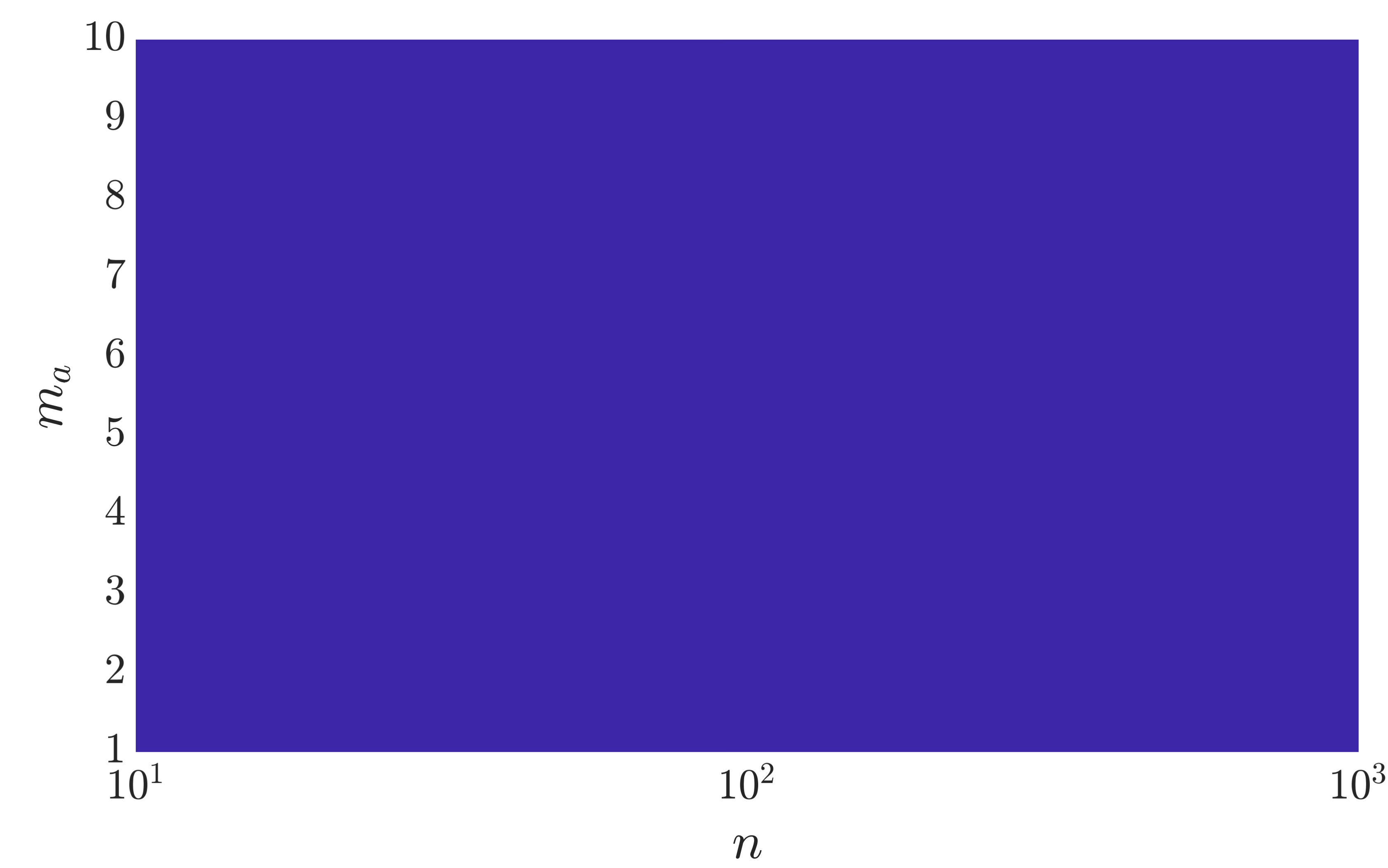}
      &
      \hspace*{-20pt}  \includegraphics[height=0.13\textheight]{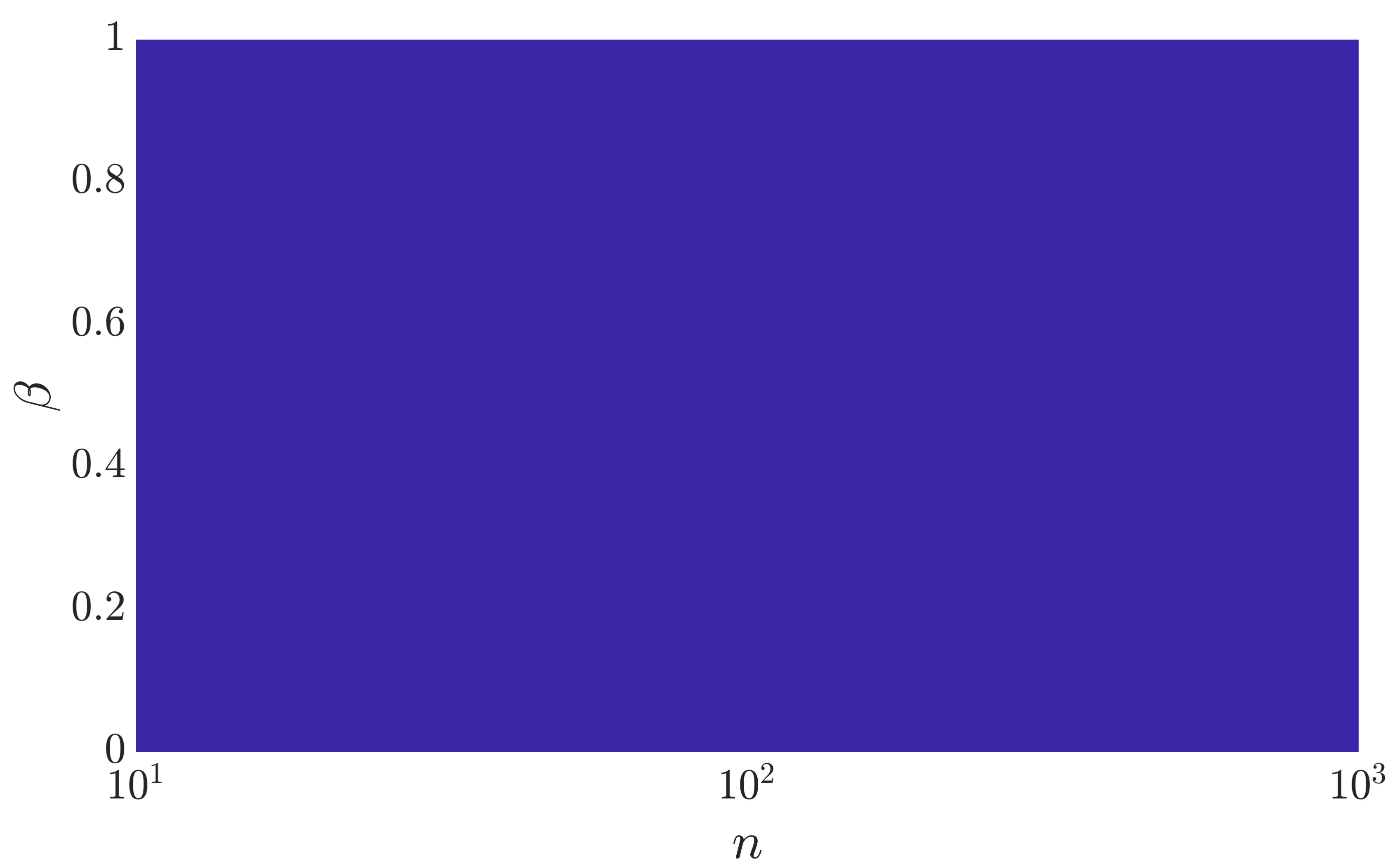}
      \\
      \rotatebox{90}{\hspace*{2ex} \tiny Raw Connectivity}
      &
      \includegraphics[height=0.13\textheight]{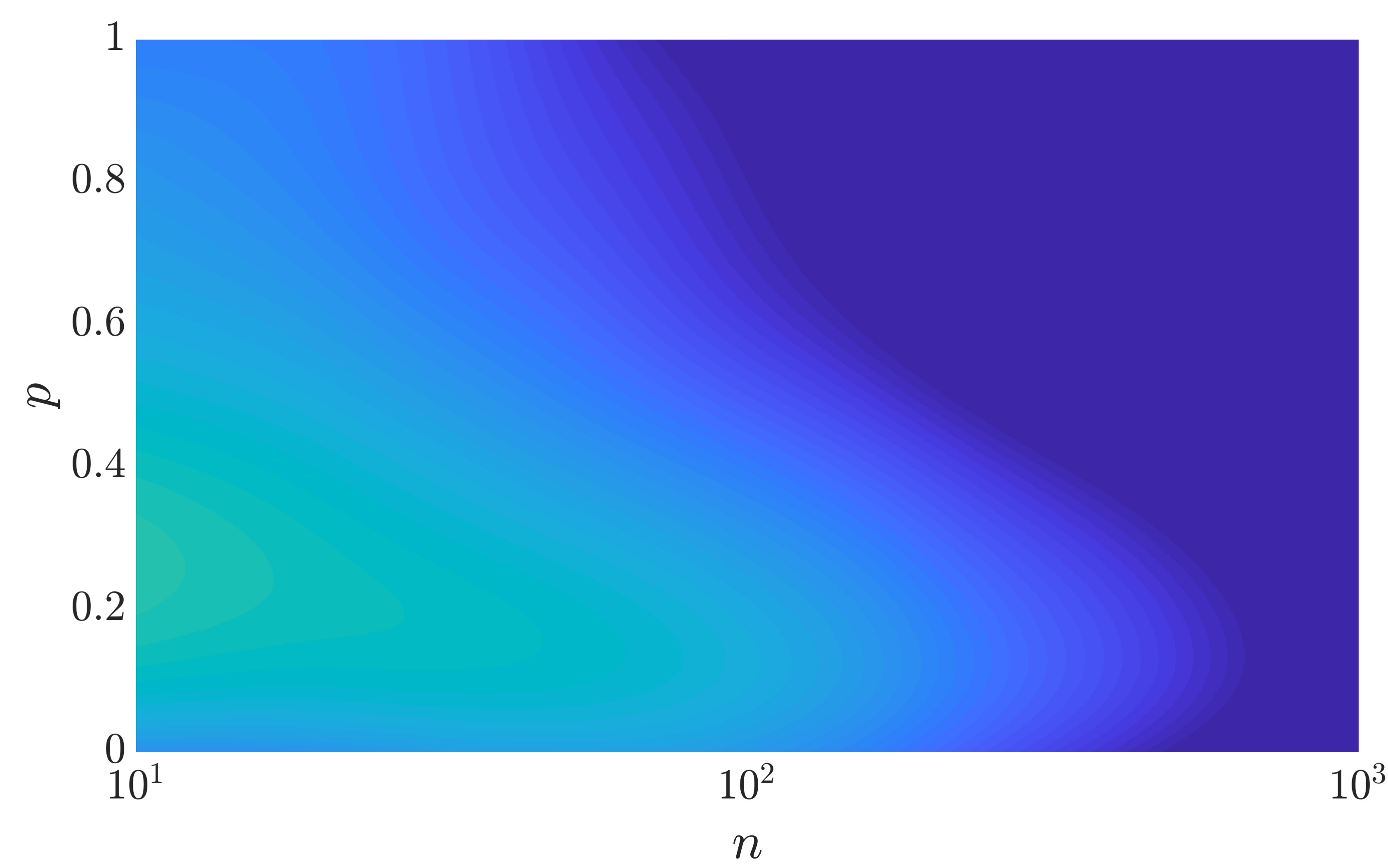}
      &
      \includegraphics[height=0.13\textheight]{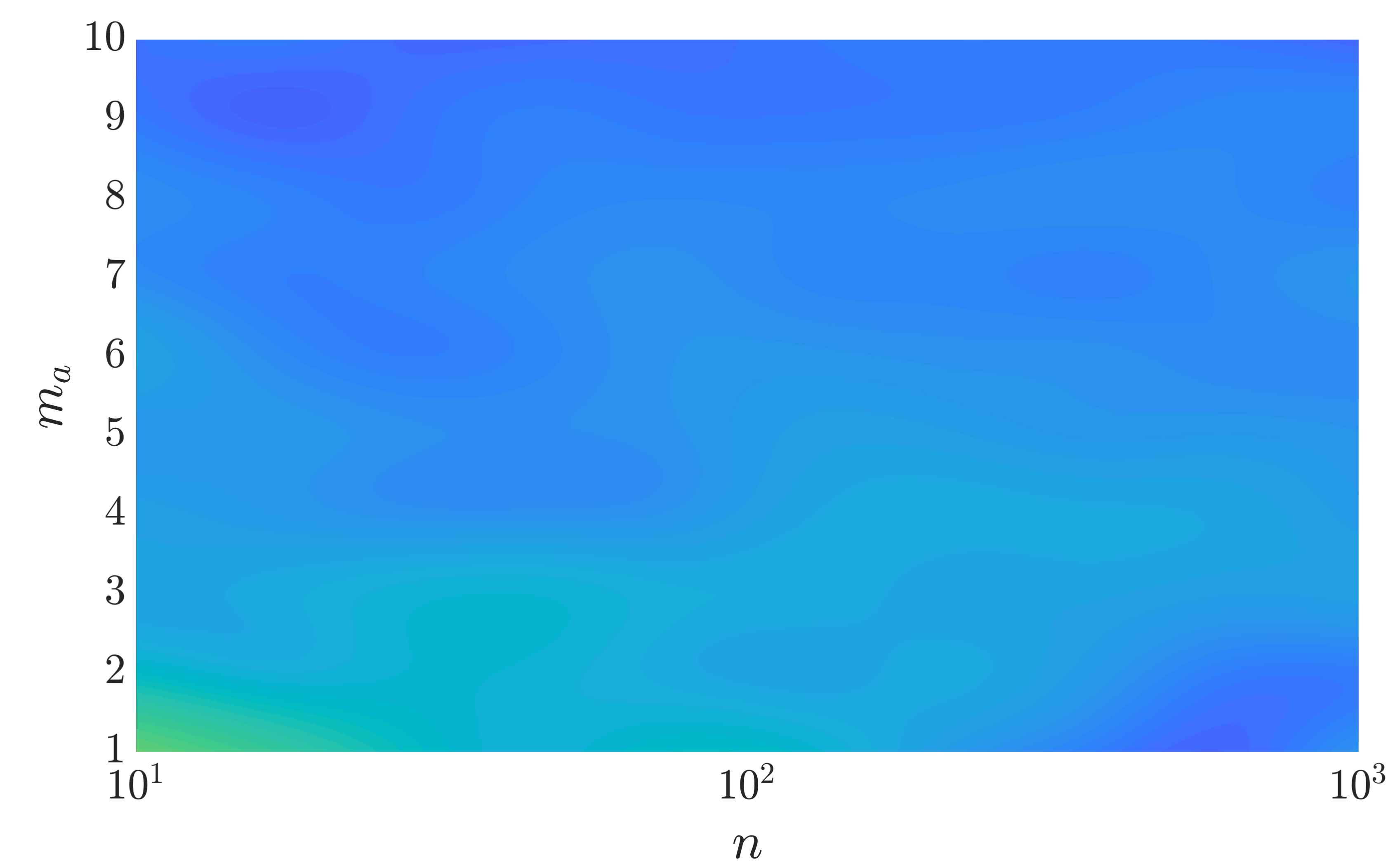}
      &
      \hspace*{-20pt}  \includegraphics[height=0.13\textheight]{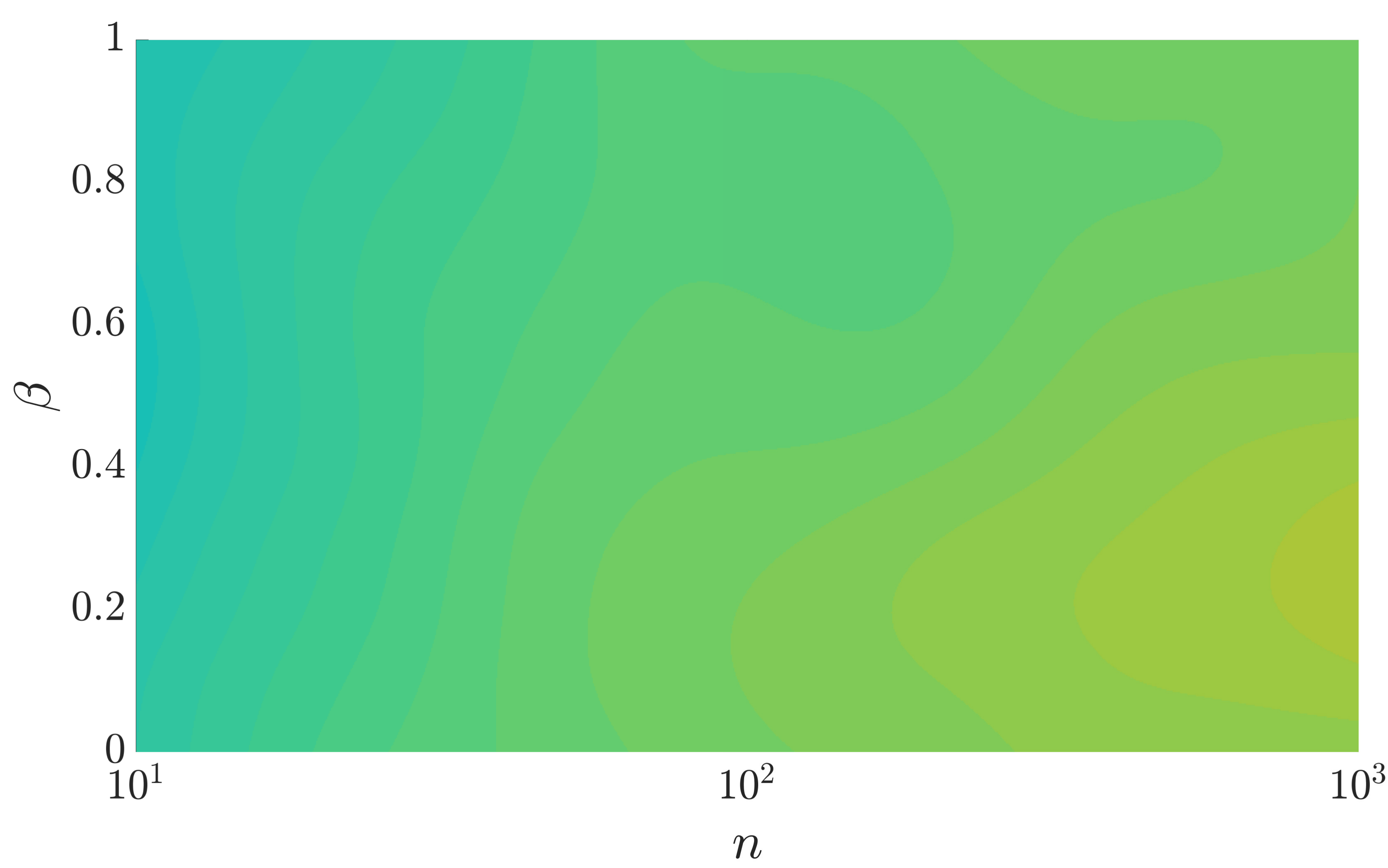}
      \\
      \rotatebox{90}{\hspace*{2ex} \tiny \parbox{20ex}{Undirected\\ Raw Connectivity}}
      &
      \includegraphics[height=0.13\textheight]{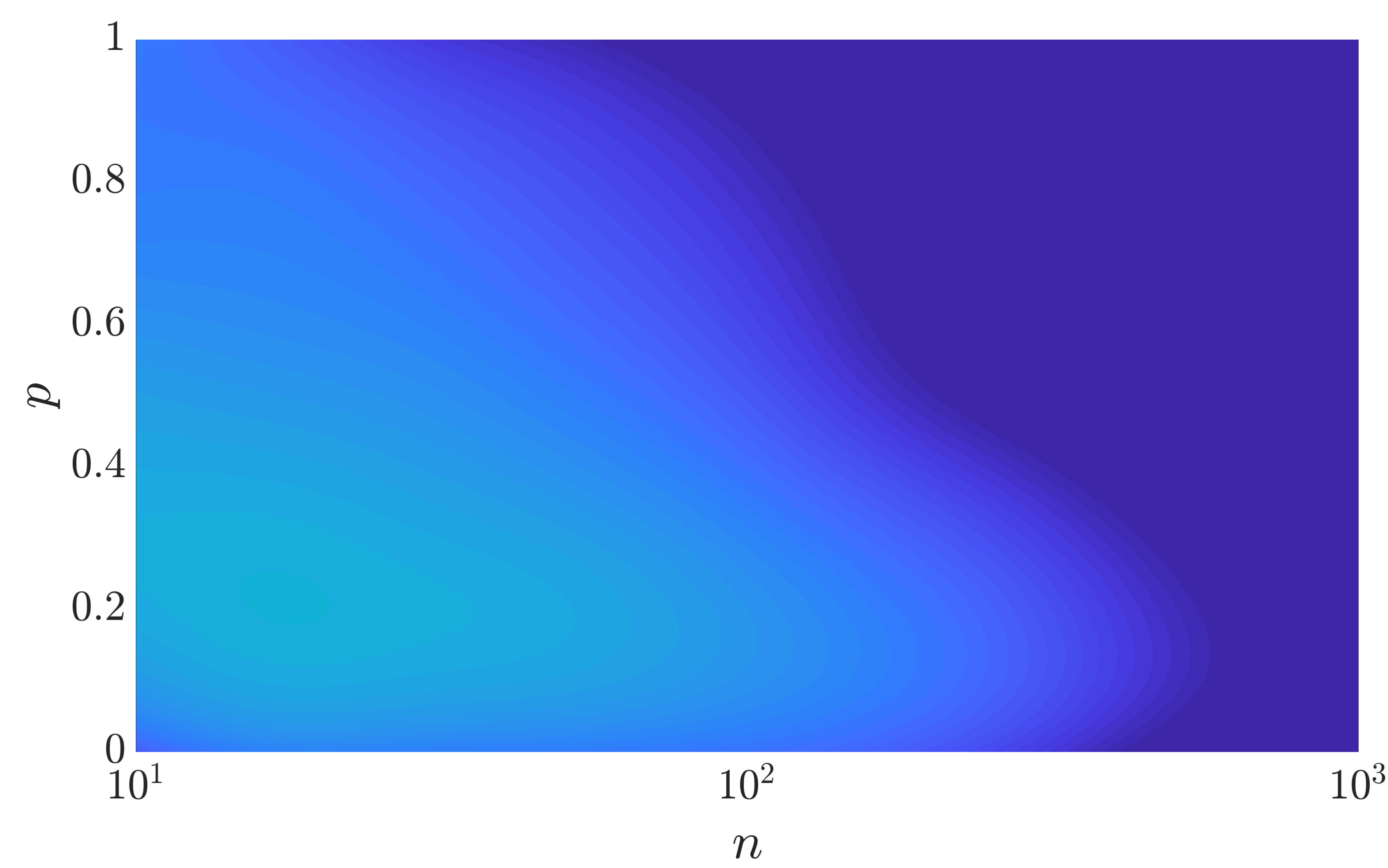}
      &
      \includegraphics[height=0.13\textheight]{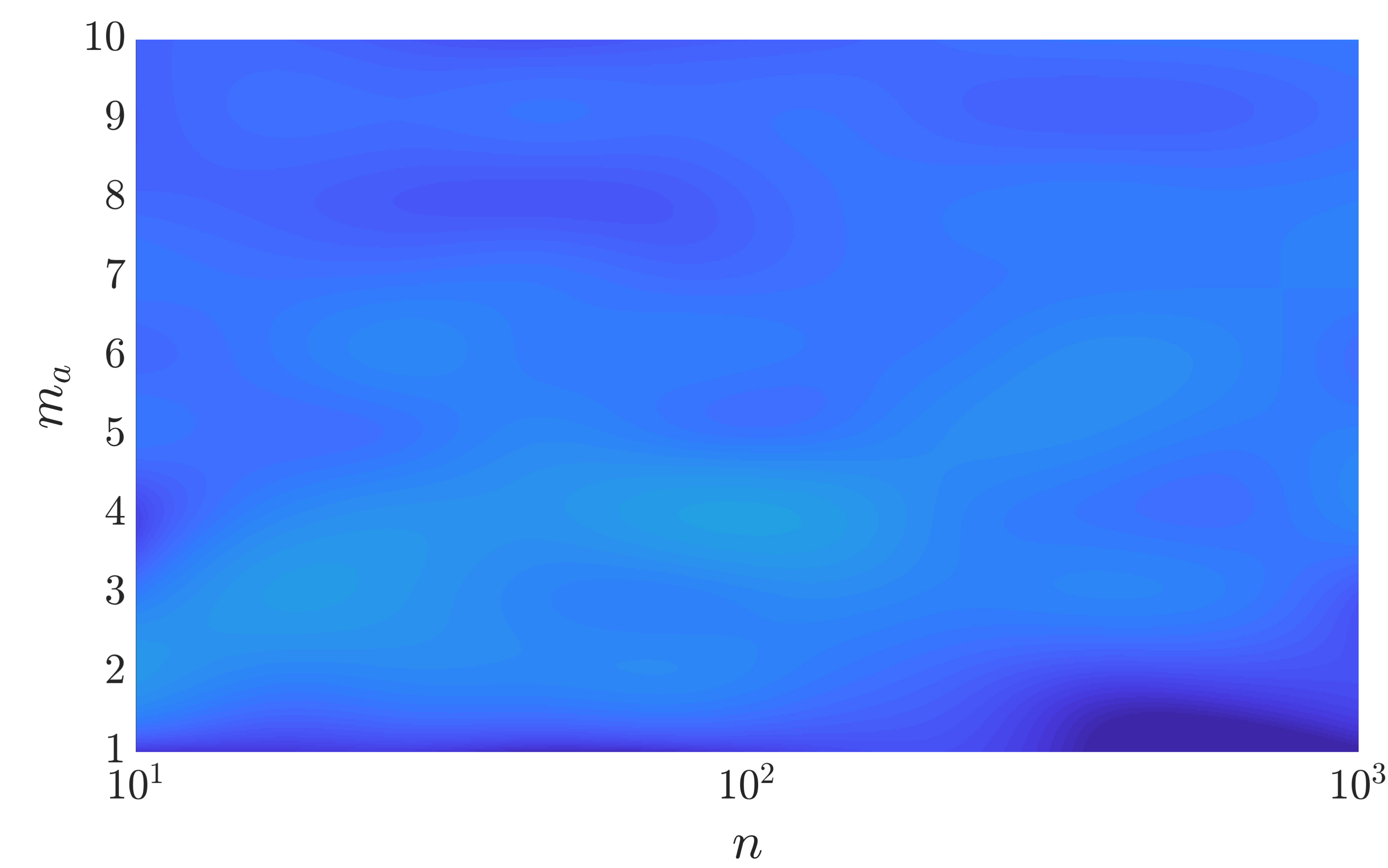}
      &
      \hspace*{-20pt}  
      \includegraphics[height=0.13\textheight]{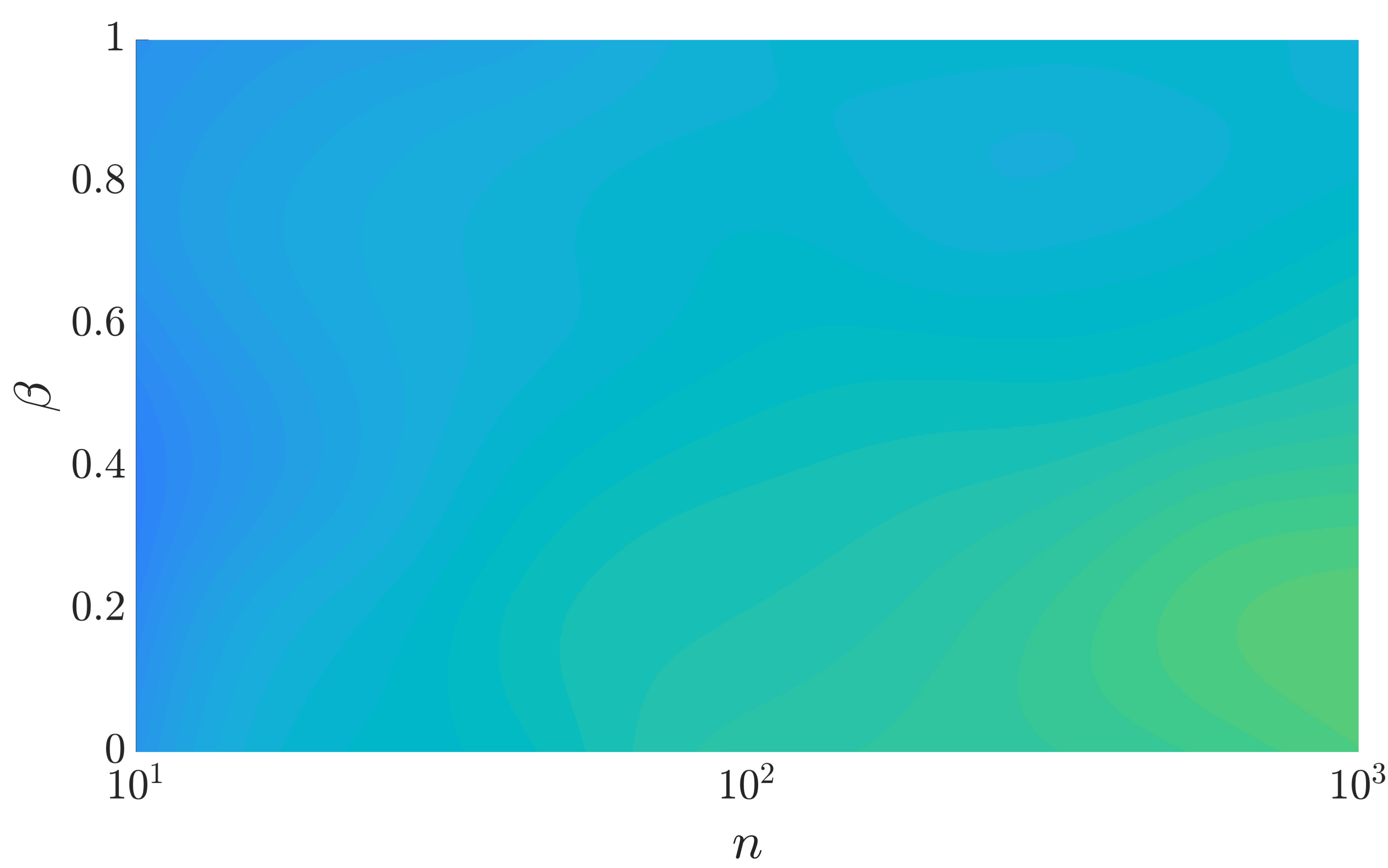}
      \\
      \rotatebox{90}{\hspace*{2ex} \tiny \parbox{20ex}{Weighted Undirected
          Raw Connectivity}}
      &
      \includegraphics[height=0.13\textheight]{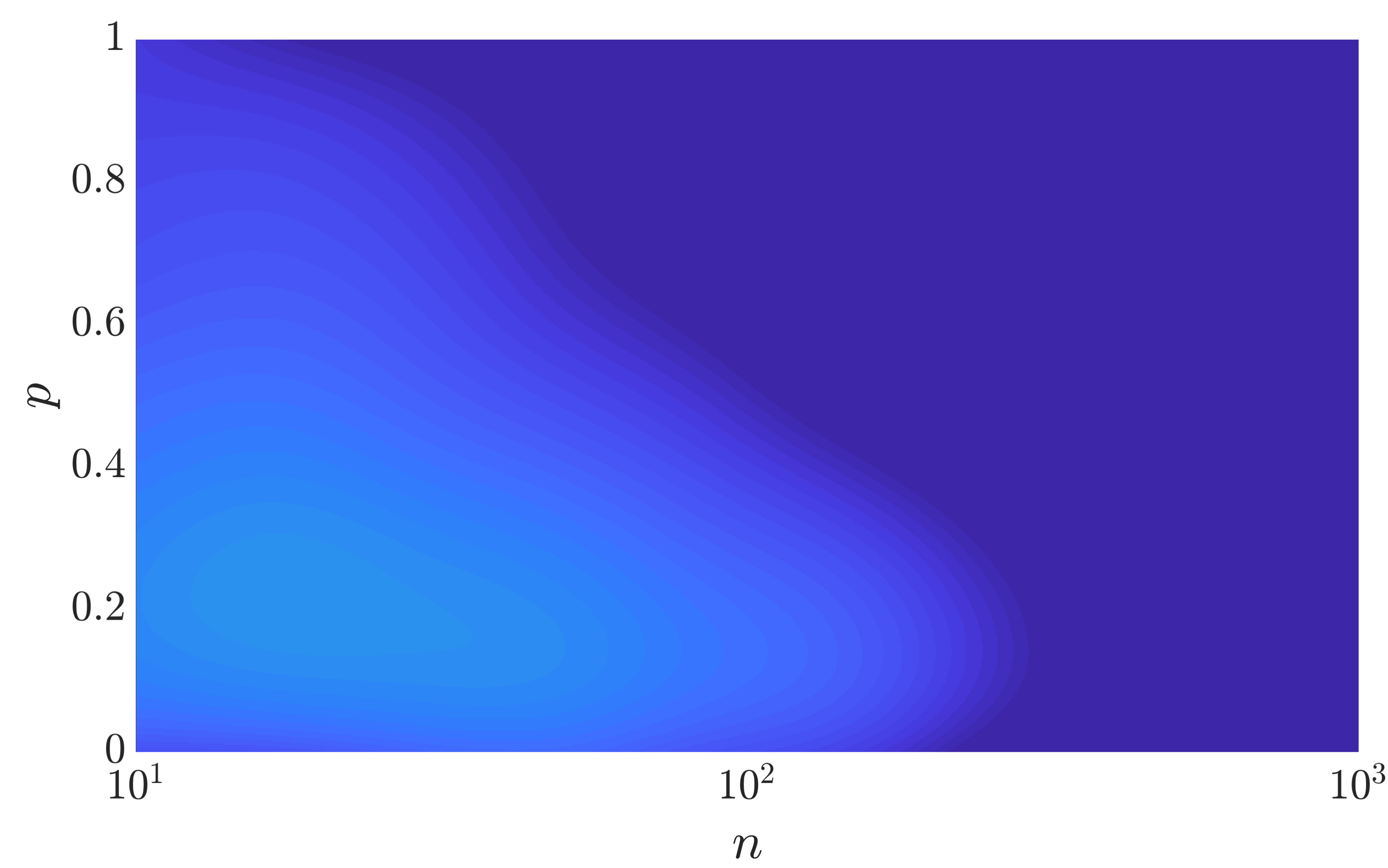}
      &
      \includegraphics[height=0.13\textheight]{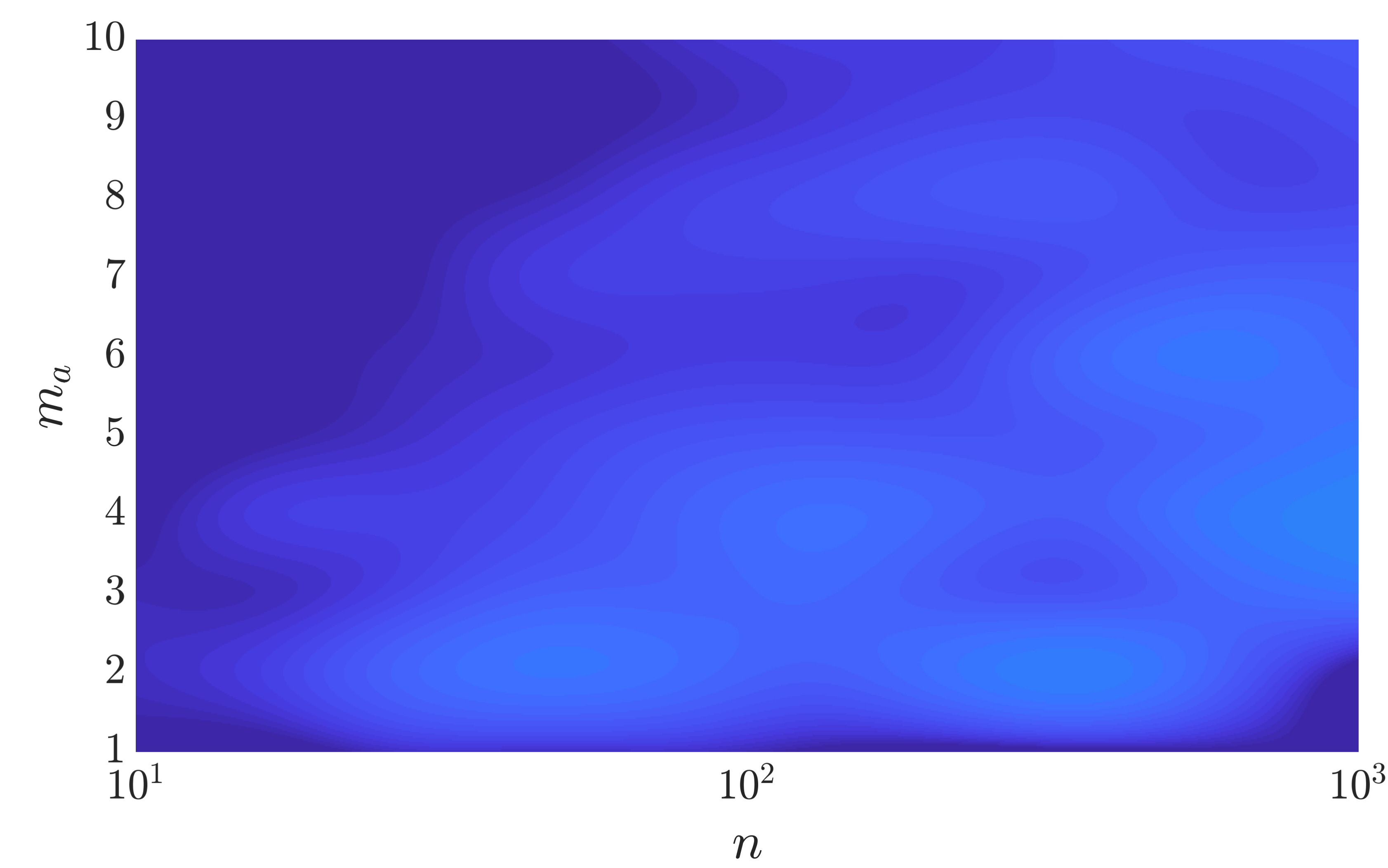}
      &
      \hspace*{-20pt}  
      \includegraphics[height=0.13\textheight]{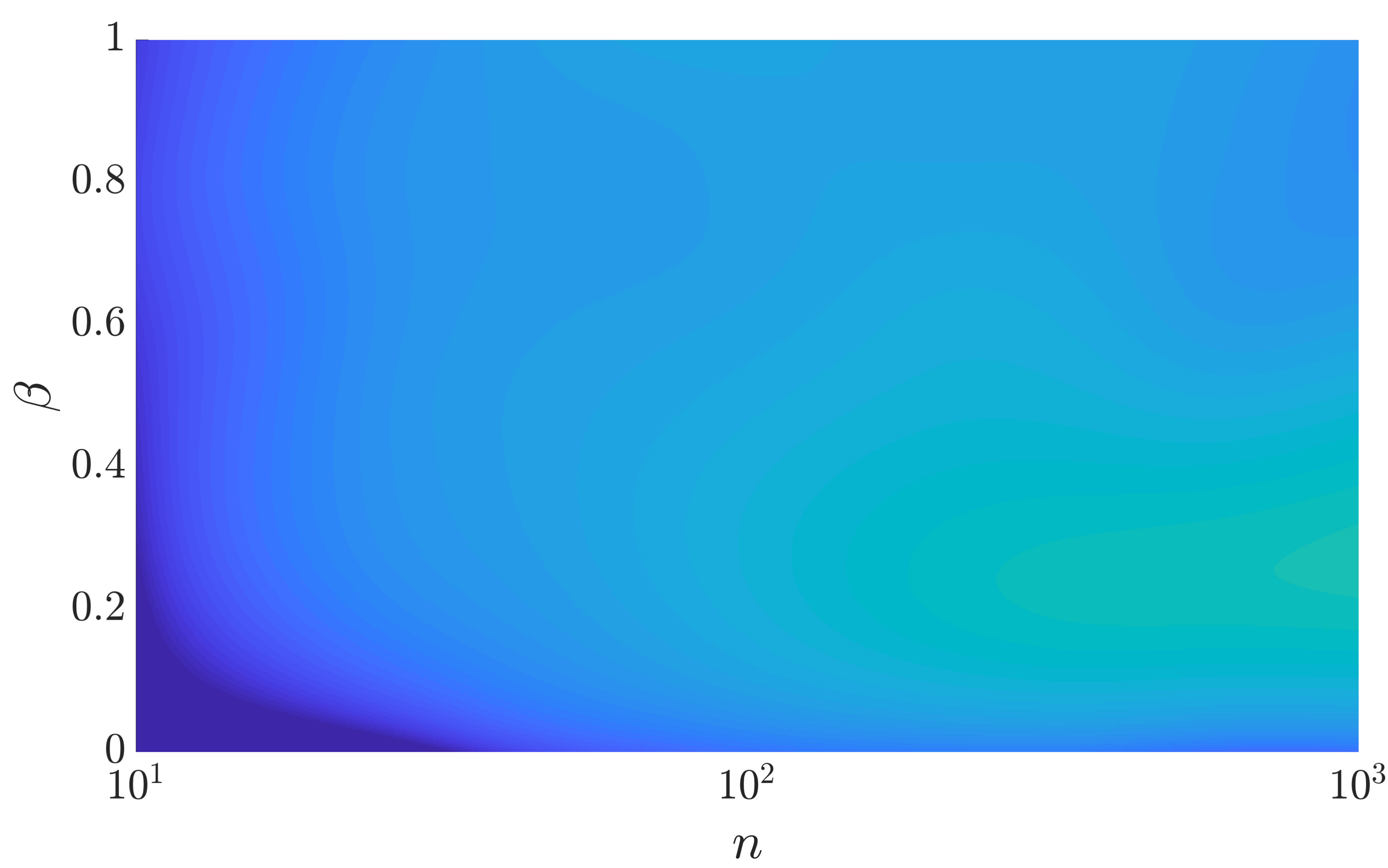}
    \end{tabular}
    \end{center} 
    \caption{Average value of $\chi$ for different  methods of obtaining dynamical
      adjacency matrix $A$ from static connectivity~$C$. The plots
      show the effect of these methods on  TVCS.
      The details on how to obtain the plots are similar to Figure~6
      in the main text. All matrices are normalized by their spectral
      radius for uniformity and comparison. The plots show a sizable
      enhancement (respectively depression) of $\chi$ by the transmission
      (respectively induction) method compared to raw connectivity, except
      for Erd\"os-R\'enyi networks whose $\chi$ maintains a robust
      pattern irrespective of the method of obtaining the dynamic
      adjacency matrix $A$ from the raw static
      connectivity~$C$.}\label{sfig:rand-all}
\end{figure}

Unless otherwise stated, we use the transmission method in this work. 
Nevertheless, it is to be noted that the method used for obtaining $A$ from $C$ can have profound effects on network controllability and should thus be chosen carefully. Figure~\ref{sfig:rand-all} illustrates this concept by showing the mean $\chi$-value of ER, BA, and WS networks for a number of different choices for this transformation. 

\section{Comparison Between Gramian-based Measures of Controllability}\label{sec:measures}

\new{
In this section, we first derive and elaborate on the relationship between the eigenvalues of the Gramian and control energy. Then, we discuss the different Gramian-based measures of controllability and their respective properties. 
}

\new{
Assume that $\Wc_K$ is invertible (the network dynamics~\eqref{eq:dyn} is controllable).
Then, for any $x_f \in \real^n$, among the (usually
infinitely many) choices of $\{u(k)\}_{k = 0}^{K - 1}$ that take the
network from $x(0) = 0$ to $x(K) = x_f$, the one that has the smallest
energy is given by~\cite[Thm 6.1]{CTC:98} 
\begin{align*}
  u^*(k) = b(k)^T (A^T)^{K - 1 - k} \Wc_K^{-1} x_f, \quad k \in \{0,
  \ldots, K - 1\}.
\end{align*}
Similar expression holds for arbitrary $x_0$, but it is
  customary to evaluate control energy starting from the network's
  unforced equilibrium $x = 0$. It is immediate to verify that this gives the minimal energy $\sum_{k
  = 0}^{K - 1} u^{*2}(k) = x_f^T \Wc_K^{-1} x_f$. Therefore, the
unit-energy reachability set is given by
\begin{align*}
\setdef{x_f \in \real^n}{x_f^T \Wc_K^{-1} x_f \le 1}.
\end{align*}
Since $\Wc_K^{-1}$ is positive definite, this is a hyper-ellipsoid in
$\real^n$, with axes aligned with the eigenvectors of $\Wc_K$. Let
$(\lambda_i, v_i)$ be an eigen-pair of $\Wc_K$ and $x_f = c
v_i$. Then,
\begin{align*}
  x_f^T \Wc_K^{-1} x_f \le 1 \Leftrightarrow c^2 \lambda_i^{-1} \le 1
  \Leftrightarrow |c| \le \lambda_i^{1/2},
\end{align*}
showing that the axis lengths of this hyper-ellipsoid are given by the
square roots of the eigenvalues of $\Wc_K$.
Intuitively, the ``larger'' the reachability
hyper-ellipsoid, the ``more controllable" the network dynamics
(equation~\eqref{eq:dyn}) are.
To quantify how large the hyper-ellipsoid is, several measures based
on the eigenvalues of $\Wc_K$ have been proposed in the
literature~\cite{PCM-HIW:72,THS-FLC-JL:16,FP-SZ-FB:14}.
}
Let $\lambda_1 \ge \lambda_2 \ge \cdots \ge \lambda_n \ge 0$
denote the eigenvalues of $\Wc_K$. The most widely used Gramian-based
measures are
\begin{itemize}
\item $\tr(\Wc_K) = \lambda_1 + \lambda_2 + \cdots + \lambda_n$,
\item $\tr(\Wc_K^{-1})^{-1} = (\lambda_1^{-1} + \lambda_2^{-1} +
  \cdots + \lambda_n^{-1})^{-1}$,
\item $\det(\Wc_K) = \lambda_1 \lambda_2 \cdots \lambda_n$,
\item $\lambda_{\min}(\Wc_K) = \lambda_n$.
\end{itemize}
It is clear from these relationships that all these measures, except
for $\tr(\Wc_K)$, approach $0$ if $\lambda_n \to 0$. This property,
i.e., the behavior of a measure as $\lambda_n \to 0$, is the most
critical difference between $\tr(\Wc_K)$ and the other three
measures. For the rest of this discussion, let $f_c(\cdot)$ be any of
$\tr((\cdot)^{-1})^{-1}$, $\det(\cdot)$, or
$\lambda_{\min}(\cdot)$. Since the network is (Kalman-) controllable
if and only if $\lambda_n > 0$, having $f_c(\Wc_K) > 0$ guarantees
network controllability while $\tr(\Wc_K) > 0$ does not. This is a
major disadvantage of $\tr(\Wc_K)$ for small networks, where
controllability in all directions in state space is both achievable
and desirable. As the size of the network grows, however, $\lambda_n$
typically decays exponentially fast to zero~\cite{FP-SZ-FB:14},
irrespective of network structure. This exponential decay of
worst-case controllability is even evident in the example network of
Figure~\ref{fig:tv-adv}(a) comprising of only $n = 5$ nodes.

Computationally, this means that $\lambda_n$ (and in turn
$f_c(\Wc_K)$) can quickly drops below machine precision as $n$ grows. In
fact, for $K = 10$ and double-precision arithmetics, this happens for
$n \sim 15$, making the TVCS (equation~\eqref{eq:opt}) with $f = f_c$ numerically
infeasible (as it involves the comparison of $f_c(\Wc_K)$ for
different $\{b_k\}_{k = 0}^{K - 1}$, which may be zero up to machine
accuracy). Further, notice that the computational complexity of TVCS
for $f = f_c$ grows as $n^K$ due to the NP-hardness of TVCS, 
enforcing the use of sub-optimal greedy algorithms even if machine
precision was not a concern (see~\cite{YZ-FP-JC:16-cdc} and the
references therein for details).

In addition to the computational aspects of TVCS, the exponential
decay of $\lambda_n$ also has theoretical implications for the choice
of $f$. When using $f = f_c$, TVCS seeks to assign the control nodes
$\{\iota_k\}_{k = 0}^{K - 1}$ such that controllability is maintained
in all directions in the state space, with special emphasis on the
hardest-to-reach directions. The use of $\tr(\Wc_K)$, on the other
hand, involves maximizing the average of Gramian eigenvalues, which
usually strengthens the largest eigenvalues and spares the few
smallest ones. In large networks, the latter is in general more
realistic as controllability is hardly needed in all $n$ directions of
the state space. As discussed in detail
in~\cite{SG-FP-MC-QKT-ABY-AEK-JDM-JMV-MBM-STG-DSB:15}, this seems to
be the case in the resting-state structural brain networks: this paper
shows that $\tr(\Wc_K)$ is maximized by controlling specific brain
regions that have long been identified as the structural ``core'' or
``hubs'' of the cerebral cortex, while the Gramian is itself close to
singular.

Further, due to the same strong dependence of $f_c(\Wc_K)$ but not
$\tr(\Wc_K)$ on $\lambda_n$, we often observe that $\tr(\Wc_K)$ is
significantly less sensitive to the choice of the control nodes
$\{\iota_k\}_{k = 0}^{K - 1}$, leading to orders of magnitude smaller
$\chi$ than that of $f_c(\Wc_K)$ (Figure~\ref{fig:tv-adv}(b)). This
means that $\Vc$ is only a small subclass of networks that benefit
from TVCS measured by $f_c$. This also has a clear interpretation,
since maintaining controllability in all directions in the state space
requires a broader distribution of the control nodes that facilitates
the reach of the control action $\{u(k)\}_{k = 0}^{K - 1}$ to all the
nodes in the network.

\new{
Finally, we highlight the need for development and analysis of measures that are neither strongly reliant on the least controllable directions (such as $f_c(\Wc_K)$) nor mainly ignore them (such as $\tr(\Wc_K)$). Two such candidates are:
\begin{itemize}
\item $\tr(C^T \Wc_K C)$ where $C$ is a matrix (or vector) with columns that point towards some particular directions of interest in the state space. This measure is a modular set function similar to $\tr(\Wc_K)$~\cite{THS-FLC-JL:16}, but the extensions of the notion of $2k$-communicability and the relationship between class $\Ic$/$\Vc$ networks and scale-heterogeneity are unclear;
\item appropriate approximations of $\log(f_c(\Wc_K))$. While computing the exact value of $\log(f_c(\Wc_K))$ is subject to the same issues as $f_c(\Wc_K)$ itself, approximations can be used that provide a \emph{mitigation} of the effects of the smallest eigenvalues of $\Wc_K$. In the case of $f_c(\cdot) = \det(\cdot)$, e.g., various algorithms have been proposed to approximate $\log \det$ of large matrices, see, e.g.~\cite{ZB-GHG:96,RPB-RKP:99,AR:02,ICFI-DJL:11,IH-DM-JS:15,CB-PD-PK-EK-AZ:17,JF-KC-MO-SR-MF:17}. These algorithms, however, are predominantly designed with the aim of reducing the computational complexity of determinant calculation and not mitigation of the effects of its high condition number, and often rely on assumptions (such as sparsity or knowledge of lower and upper bounds on matrix eigenvalues) that do not apply to $\Wc_K$. Thus, development of \emph{appropriate} approximations of $\log(f_c(\Wc_K))$ constitutes a warranted direction for future research.
\end{itemize} 
}

\section{Relationships Between $2k$-Communicability, Degree, and Eigenvector Centrality}\label{sec:relations}

The notion of $2k$-communicability introduced in this article has
close connections with the degree and eigenvector centrality in the
limit cases of $k = 1$ and $k \to \infty$, respectively. 
Recall that
the out-degree centrality and $2$-communicability of a node $i \in \Nc$
are defined as, respectively,
\begin{align*}
d^\text{out}_i &= \sum_{j = 1}^n a_{ji},
\\
R_i(1) &= \sum_{j = 1}^n a_{ji}^2.
\end{align*}
Therefore, if the network is unweighted (i.e., all the edges have the
same weight), then $R_i(1) \propto d^\text{out}_i$, so
$2$-communicability and out-degree centrality result in the same
ranking of the nodes (in particular, $r(1)$ is the node with the
largest out-degree). As edge weights become more heterogenous, these
two rankings become less correlated, with $2$-communicability putting
more emphasis on stronger weights.

A similar relation exists between $\infty$-communicability and left
eigenvector centrality, as we show next. Let $v_1, u_1 \in \real^n$ be
the right and left Perron-Frobenius eigenvectors of $A$, respectively,
normalized such that $v_1^T v_1 = u_1^T v_1 = 1$. Since the network is
by assumption strongly connected and aperiodic, we have
\begin{align*}
\lim_{k \to \infty} \Big(\frac{1}{\rho(A)} A\Big)^k = v_1 u_1^T.
\end{align*}
Thus for any $i \in \Nc$,
\begin{align*}
\lim_{k \to \infty} \Big(\frac{1}{\rho(A)} \Big)^{2k} R_i(k) &= \lim_{k \to \infty} \Big(\frac{1}{\rho(A)} \Big)^{2k} \big((A^k)^T A^k\big)_{ii} = (u_1 v_1^T v_1 u_1^T)_{ii} = u_{1, i}^2.
\end{align*} 
Given that dividing $R_i(k)$ by $\rho(A)^{2k}$ for all $i$ does not
change the ranking of nodes, we define $R_i(\infty) = u_{1, i}^2$ for
all $i$. Since squaring non-negative numbers preserves their order,
nodal rankings based on $\infty$-communicability and left eigenvector
centrality are identical.

\section{Nodal Dominance}\label{sec:dom}

Among the networks where the nodes with the greatest $R_i(1)$ and
$R_i(\infty)$ coincide (i.e., $r(0) = r(\infty)$), there is a higher
chance (than in general) that any network belongs to class
$\Ic$. However, about half of these networks still belong to class
$\Vc$, meaning that there exists $1 < k < \infty$ such that $r(k) \neq
r(0)$. To assess the importance of this time-variation of optimal
control nodes, we define the \emph{dominance} of the node $r(0)$ (over
the rest of the network) as follows. Let $r'(0)$ be the index of the
node with the second largest $R_i(1)$ (largest after removing
$r(0)$). Similarly, let $r'(\infty)$ be the index of the second
largest $R_i(\infty)$. We define
\begin{align*}
  \text{Dominance of} \ r(0) = \min\Big\{&\frac{R_{r(0)}(0) -
    R_{r'(0)}(0)}{R_{r(0)}(0)}, \frac{R_{r(0)}(\infty) -
    R_{r'(\infty)}(\infty)}{R_{r(0)}(\infty)}\Big\}.
\end{align*} 
A small dominance indicates that another node has very similar value
$R_i(0)$ \emph{or} $R_i(\infty)$ to $r(0)$, while a large dominance is
an indication of a large gap between $R_{r(0)}(k)$ and the next
largest $R_i(k)$ for both $k = 0$ \emph{and} $k \to \infty$.

\section{Networks with Multiple Inputs}\label{sec:mi}

Consider a multiple-input network, namely, a network in which $m \ge 1$ nodes are controlled at every time step. Let $\iota_k^1, \dots, \iota_k^m \in \Nc$ denote the indices of the control nodes at every time $k$, and $\iotab_k = \{\iota_k^1, \dots, \iota_k^m\}$. Then, the corresponding TICS and TVCS are defined as
\begin{subequations}
\begin{align}\label{eq:opt-TI-MI}
\max_{\iotab_0, \dots, \iotab_{K - 1} \in \Nc} \hspace{10pt} &f(\Wc_K)
\\
\text{s.t.} \hspace{25pt} &\iotab_0 = \cdots = \iotab_{K - 1} \label{eq:opt-TI-MI-a}
\end{align}
\end{subequations}
and\begin{align}\label{eq:opt-MI}
\max_{\iotab_0, \dots, \iotab_{K - 1} \in \Nc} \hspace{10pt} f(\Wc_K),
\end{align}
respectively. Accordingly, a multiple-input network is said to belong to class $\Ic$ if the solution of~\eqref{eq:opt-MI} satisfies~\eqref{eq:opt-TI-MI-a}, and to class $\Vc$ otherwise.

Clearly, for a multiple-input network to belong to class $\Vc$, any of
the first $m$ largest of $\{R_i(k)\}_{i = 1}^n$ should change over
time, which is often implied by (but does not imply) a change in
$r(k)$. Therefore, the condition of Theorem~\ref{thm:main-suf} 
\new{can still be used as a tight proxy for networks in $\Vc$,}
but is too conservative and can be relaxed as follows: assume that $A$ is
irreducible, aperiodic, and diagonalizable. Let $\setdef{r_j^d \in
  \real^n}{j \in \Jc^d}$ be the set of \new{nodes with the $m$ highest}
$2$-communicabilities, where the index set $\Jc^d$ accounts for
different choices of rankings if there are nodes with equal
$2$-communicabilities.  Similarly, let $\setdef{r_j^c \in \real^n}{j
  \in \Jc^c}$ be the set of \new{nodes with the $m$ highest}
  centralities. Then,
if $r_{j_1}^d \neq r_{j_2}^c$ for all $(j_1, j_2) \in \Jc^d \times
\Jc^c$, the network belongs to class $\Vc$ when $K$ is sufficiently
large. The proof of this statement is a straightforward generalization
of the proof of Theorem~\ref{thm:main-suf} and thus omitted.

Similarly, the three conditions in Theorem~\ref{thm:main-nec} can be
generalized to undirected multiple-input networks as follows (with
similar proofs as the proof of Theorem~\ref{thm:main-nec}):
\begin{enumerate}
\item[\emph{(i)}] For all $i \in \until{m}$,
  \begin{align*}
    \frac{1 - w_{i 1}}{\sum_{\ell \le m + 1, \ell \neq i + 1}
      w_{\ell 1}} \le \frac{|\lambda_1| - |\lambda_2|}{|\lambda_1|
      - |\lambda_n|}.
  \end{align*}
  This condition can be simplified, at the expense of being more
  conservative, to
  $\frac{1 - w_{i 1}}{i w_{i 1}} \le \frac{|\lambda_1| -
    |\lambda_2|}{|\lambda_1| - |\lambda_n|}$,
  for all $i \in \until{m}$,
\item[\emph{(ii)}] for all $i \in \until{m}$, $w_{i 2} = 1 -
  w_{i 1}$,
\item[\emph{(iii)}] the network has three or fewer nonzero
  eigenvalues with different absolute values and
  \begin{align*}
    R_1(1) \ge R_2(1) \ge \cdots \ge R_m(1) \ge R_i(1),
  \end{align*}
  for all $i \in \{m + 1, \dots, n\}$.
\end{enumerate}

\begin{figure}
  \begin{center}
    \setlength{\tabcolsep}{2pt}
    \begin{tabular}{cccc}
    &
      {\tiny Erd\"os-R\'enyi} & {\tiny
      Barab\'asi-Albert} &  {\tiny Watts-Strogatz}
      \\
      \rotatebox{90}{\hspace*{10ex} \tiny 1 Control Node} 
      &
      \includegraphics[height=0.18\textheight]{ER_d_w_trans}
      &
      \includegraphics[height=0.18\textheight]{BA_d_w_trans}
      &
      \includegraphics[height=0.18\textheight]{WS_d_w_trans}
      \\
      \rotatebox{90}{\hspace*{10ex} \tiny 5 Control Nodes} 
      &
      \includegraphics[height=0.18\textheight]{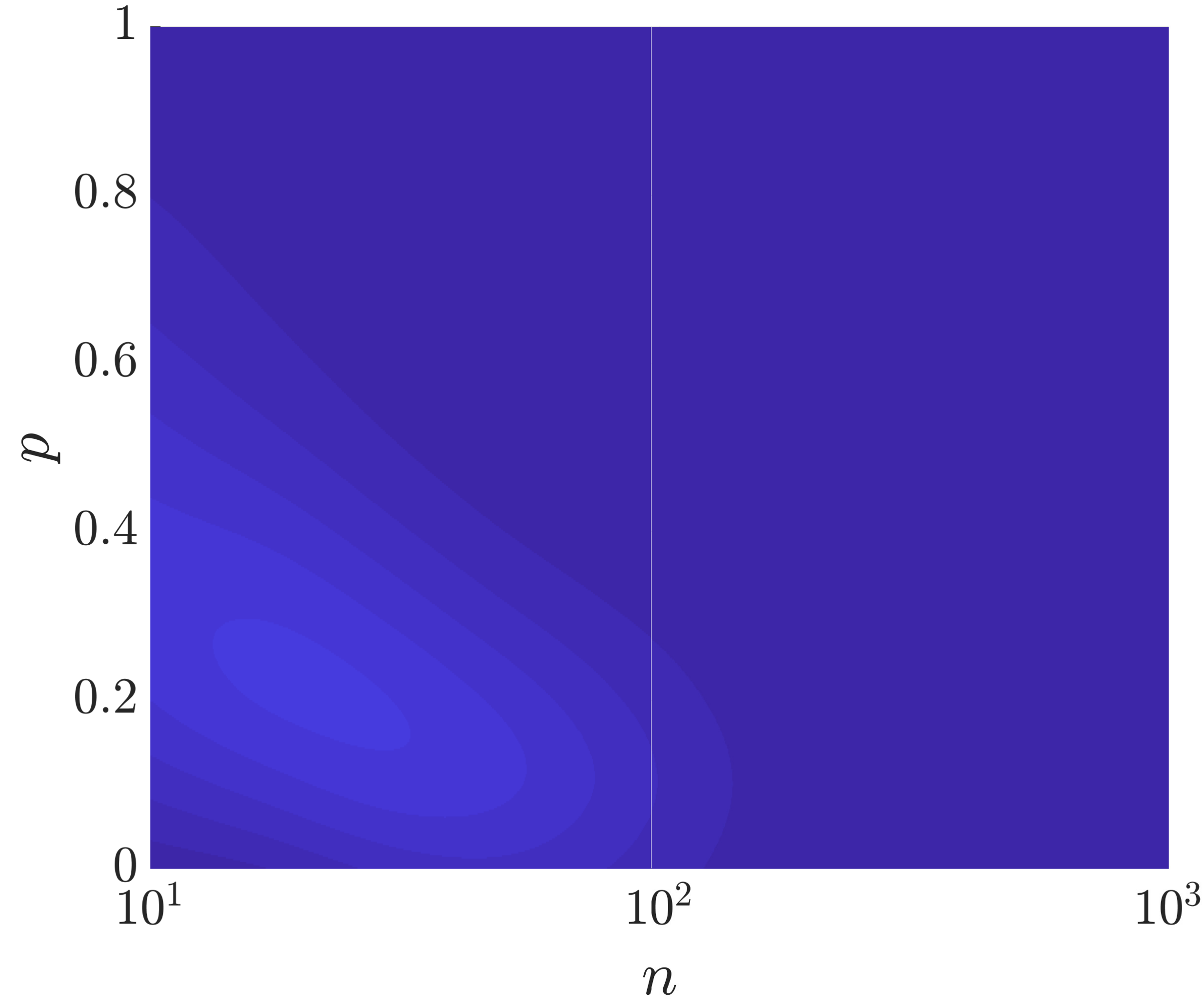}
      &
      \includegraphics[height=0.18\textheight]{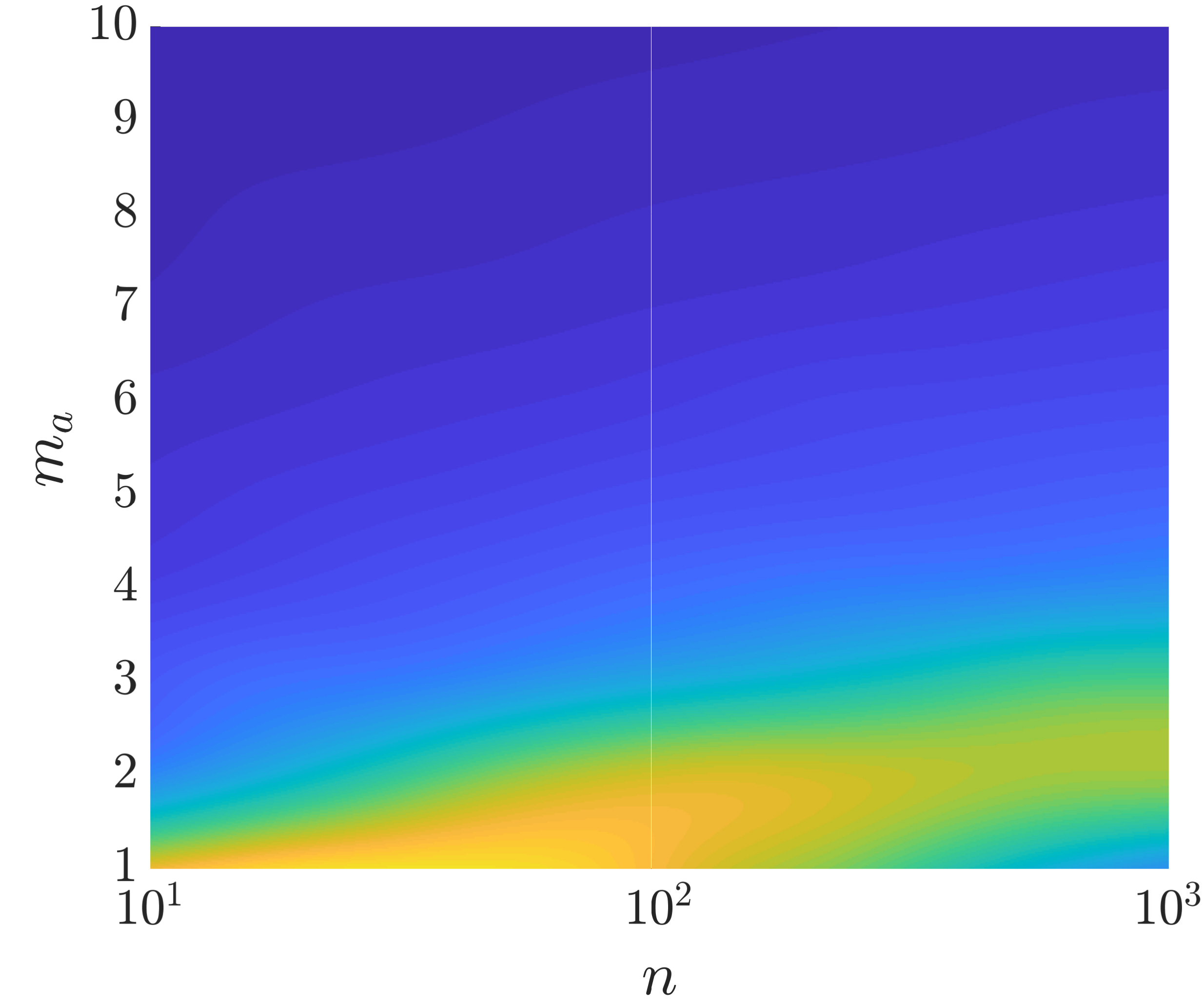}
      &
      \hspace*{-22pt}
      \includegraphics[height=0.18\textheight]{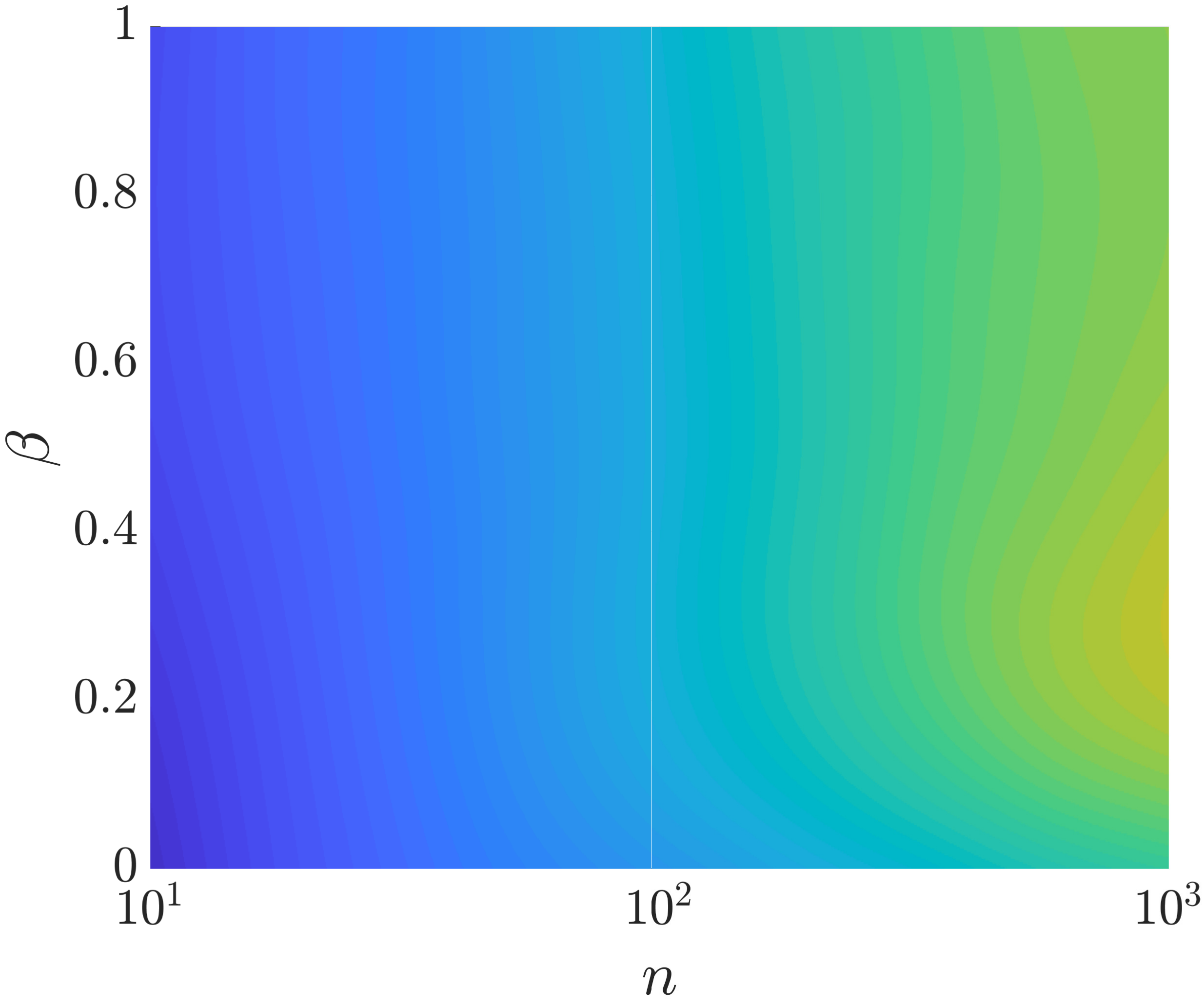}
      \\
      \rotatebox{90}{\hspace*{10ex} \tiny 10 Control Nodes} 
      &
      \includegraphics[height=0.18\textheight]{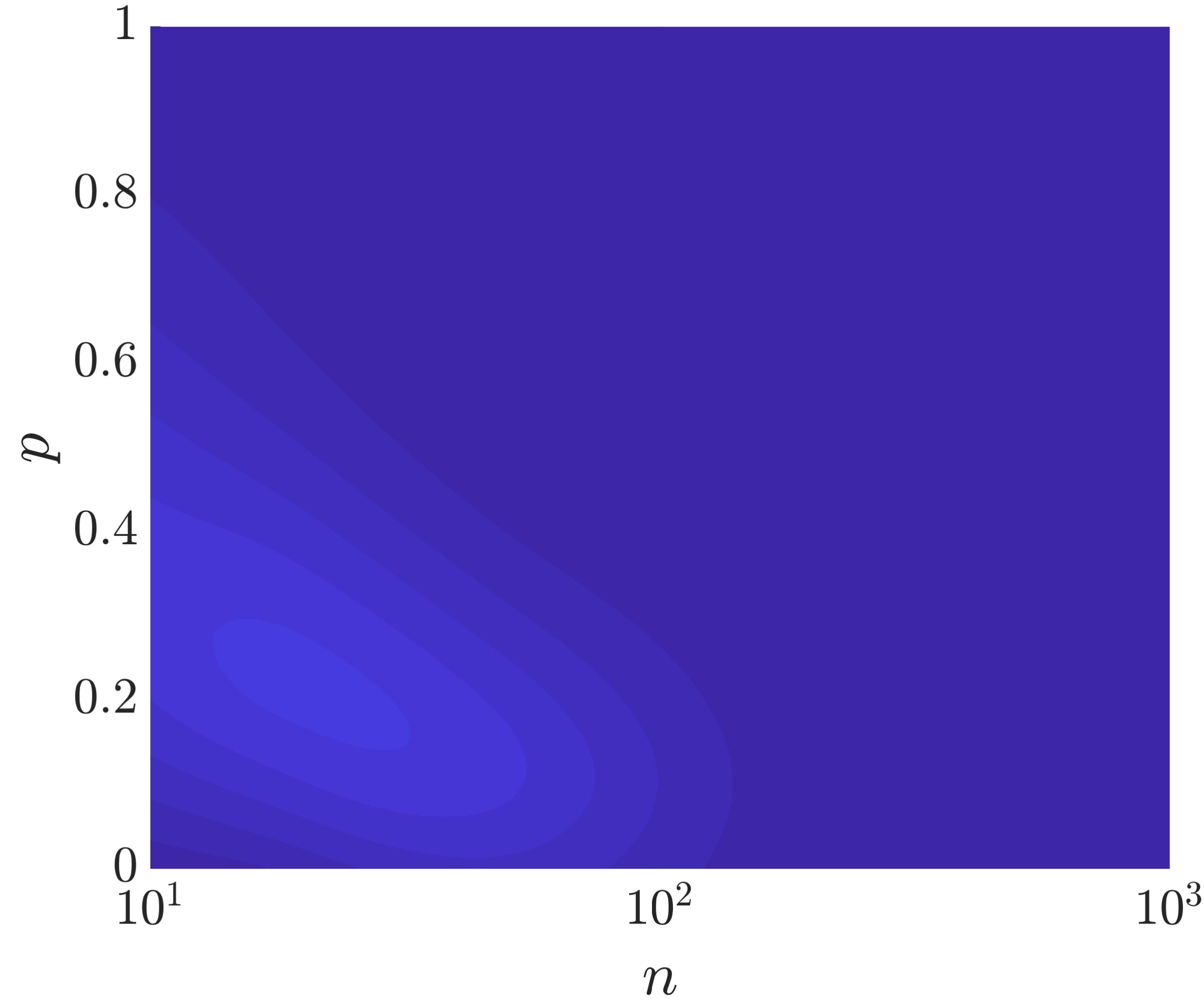}
      &
      \includegraphics[height=0.18\textheight]{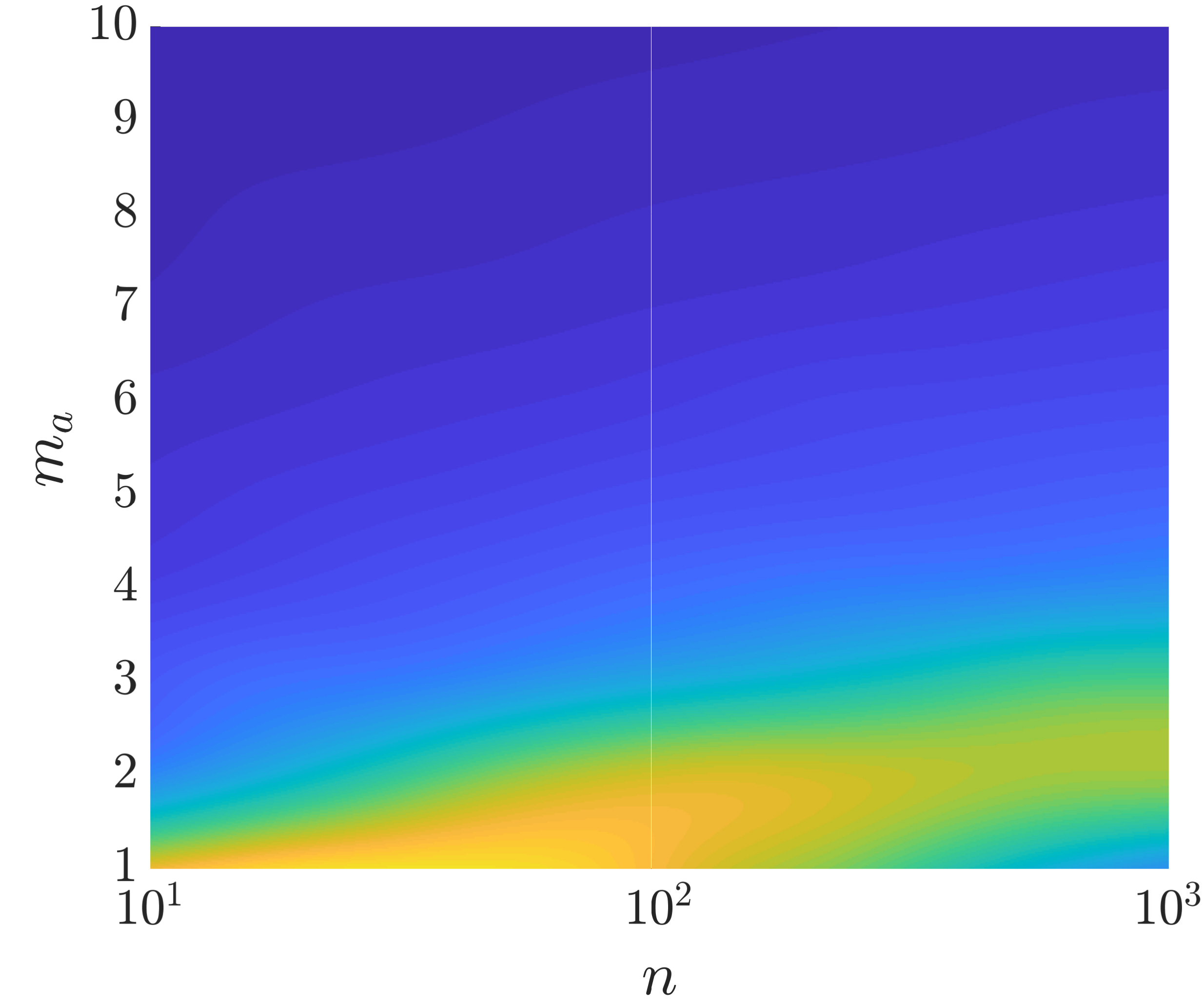}
      &
      \hspace*{-22pt}
      \includegraphics[height=0.18\textheight]{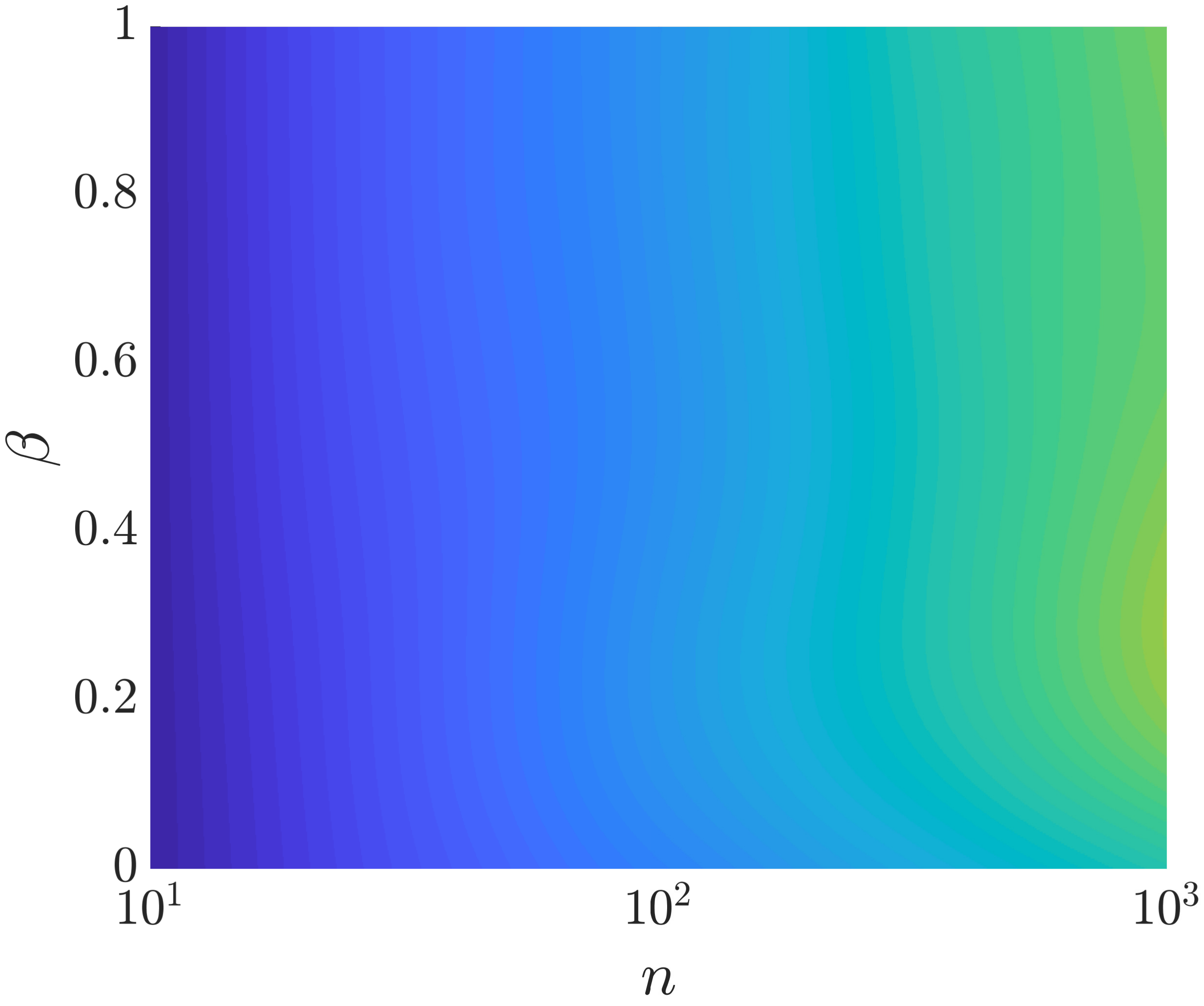}
    \end{tabular}
  \end{center} 
  \caption{Average value of $\chi$ for networks with increasing number
    of inputs.  The plots show the effect of multiple inputs on TVCS.
    The details on how to obtain the plots are similar to Figure~6 in
    the main text. The dynamic adjacency matrix $A$ is obtained from
    the raw static connectivity $C$ using the transmission method in
    all cases. These plots show a slight depression in the benefit of
    TVCS as the number of control nodes grows, despite the fact that
    networks with more control nodes have a higher probability of
    belonging to $\Vc$ (namely, having $\chi > 0$).}
  \label{sfig:mi}
\end{figure}

Finally, Figure~\ref{sfig:mi} illustrates the effect of $m$ on $\chi$-values of ER, BA, and WS networks discussed in the main text.

\section{$2 k$-Communicabilities of Simple
  Networks}\label{sec:simple-networks}

As mentioned in the main text, cf. Figure~5, the simple structure of
homogeneous undirected line, ring, and star networks allows us to
compute their $2k$-communicabilities analytically in closed form, as
derived in the following. \new{Throughout, $\integers$ denotes the set of integers and for $a,
b \in \integers$, $a \divides b$ denotes that $a$ divides $b$.}

\begin{proposition}\longthmtitle{$2 k$-communicabilities of line
    networks}\label{prop:line}
  Consider a line network of $n$ nodes with uniform link weights~$a$
  (and no self-loops). Then, for $i \in \Nc$ and $k \in
  \intpos$,
  \begin{align}\label{eq:rik-trid-full}
    R_i(k) \!=\! a^{2 k}
    \sum_{p \in \Ic} \!\left[\!{{2 k}\choose{k \!+\! p (n \!+\! 1)}}
      \!-\! {{2 k}\choose{k \!+\! p (n \!+\! 1) \!-\! i}}\!\right]\!,
  \end{align}
  where $\Ic = \{-\ceil{\!\frac{k}{n + 1}\!}, \ldots,
  \ceil{\!\frac{k}{n + 1}\!}\}$ and ${{n}\choose{k}} \triangleq 0$ if
  $k \notin \{0, \ldots, n\}$. In particular, if $i \le
  \ceil{\frac{n}{2}}$ and $k \le \ceil{\frac{n}{2}} - 1$,
  \begin{align}\label{eq:rik-trid-simple}
    R_i(k) = a^{2 k} \left[{{2 k}\choose{k}} - {{2 k}\choose{k -
          i}}\right].
  \end{align}
\end{proposition}
   \begin{proof}
     From~\cite[Lemma 1.77]{FB-JC-SM:08cor}, we have
     \begin{align*}
       \lambda_j = 2 a \cos\frac{j \pi}{n + 1} \quad \text{and} \quad
       w_{i j} \propto \sin^2\frac{i j \pi}{n + 1},
   \end{align*} 
   for $i, j \in \until{n}$ where $\propto$ accounts for
   normalization. In order to normalize the eigenvectors, we use the
   identities $\sin^2 \alpha = \frac{1}{2} (1 - \cos 2 \alpha)$ and
   \begin{align}\label{eq:trig-identity}
     \sum_{j = 1}^n \cos\frac{2 s j
       \pi}{n + 1} = -1 \quad \text{for all } s \notdivides n + 1,
   \end{align}
   to get $w_{i j} = \frac{2}{n + 1} \sin^2\frac{i j \pi}{n + 1}$ for
   all $i, j \in \until{n}$ (one can show~\eqref{eq:trig-identity} by
   multiplying and dividing the LHS by $\sin\frac{s \pi}{n + 1}$ and
   using the identity $2 \sin\alpha \cos\beta = \sin(\alpha + \beta) +
   \sin(\alpha - \beta)$ for each term).  Thus, by substitution, we
   have $R_i(k) = \frac{2 a^{2 k}}{n + 1} \sum_{j = 1}^n \tau_{i j
     k}^2$ where
   \begin{align*}
     \tau_{i j k} = 2^k \sin\frac{i
       j \pi}{n + 1} \cos^k\frac{j \pi}{n + 1}.
   \end{align*}
   By using the identity $2 \sin\alpha \cos\beta = \sin(\alpha + \beta)
   + \sin(\alpha - \beta)$, $k$ times and collecting terms, we get
   \begin{align*}
     \tau_{i j k} = \sum_{\ell = 0}^k {{k}\choose{\ell}} \sin\frac{(i +
       k - 2 \ell) j \pi}{n + 1}.
   \end{align*}
   Hence, by squaring $\tau_{i j k}$ and substituting it in $R_i(k)$,
   and using the identity $2 \sin\alpha \sin\beta = \cos(\alpha -
   \beta) - \cos(\alpha + \beta)$, we get
   \begin{align}\label{eq:cancel-sum}
     &R_i(k) = \frac{a^{2 k}}{n + 1} \sum_{\ell, r = 0}^k
     {{k}\choose{\ell}} {{k}\choose{r}} 
     \Big[\sum_{j = 1}^n \cos\frac{2(\ell - r) j \pi}{n + 1} - \sum_{j
       = 1}^n \cos\frac{2(i + k - \ell - r) j \pi}{n + 1}\Big].
   \end{align}
   However, by~\eqref{eq:trig-identity}, the two sums
   in~\eqref{eq:cancel-sum} cancel each other unless $\ell - r \divides
   n + 1$ or $i + k - \ell - r \divides n + 1$ (the cases where both of
   these happen need not be excluded since they automatically
   cancel). Thus,
   \begin{align}\label{eq:Rik-Ic1-Ic2}
     R_i(k) = a^{2 k} \Bigg[\sum_{\Ic_1} {{k}\choose{\ell}}
     {{k}\choose{r}} - \sum_{\Ic_2} {{k}\choose{\ell}}
     {{k}\choose{r}}\Bigg],
   \end{align}
   where
   \begin{align*} 
     \Ic_1 &= \setdef{(\ell, r) \in \{0, \dots, k\}^2}{\ell - r
       \divides n + 1},
     \\
     \Ic_2 &= \setdef{(\ell, r) \in \{0, \dots, k\}^2}{i + k - \ell - r
       \divides n + 1}.
   \end{align*} 
   Defining $p = \frac{n + 1}{\ell - r}$ in the first and $p = \frac{i
     + k - \ell - r}{n + 1}$ in the second sum
   in~\eqref{eq:Rik-Ic1-Ic2}, we get
   \begin{align}\label{eq:Rik-final}
     R_i&(k) = a^{2 k} 
    \sum_{p \in \Ic} \Bigg[\sum_{\ell = 0}^k \!
     {{k}\choose{\ell}} \!  {{k}\choose{\ell \!-\! p (n \!+\! 1)}}
     \!\!-\!\!  \sum_{\ell = 0}^k \!  {{k}\choose{\ell}} \!
     {{k}\choose{\ell \!+\! p (n \!+\!  1) \!-\!  i}}\Bigg],
   \end{align}
   where we have used the identity ${{k}\choose{s}} = {{k}\choose{k -
       s}}$. Equation~\eqref{eq:rik-trid-full} then follows by applying
   the formula $\sum_{\ell = 0}^k {{k}\choose{\ell}} {{k}\choose{\ell
       \pm s}} = {{2 k}\choose{k \pm s}}$~\cite[Eq. 6.69-70]{HWG:10} to
   each of the two sums in~\eqref{eq:Rik-final}. To
   get~\eqref{eq:rik-trid-simple}, note that if $i \le
   \ceil{\frac{n}{2}}$ and $k \le \ceil{\frac{n}{2}} - 1$, then the
   only nonzero term in~\eqref{eq:rik-trid-full} is the one
   corresponding to $p = 0$.
 \end{proof}
 
 According to this result, in the case of no self-loops, the value of
 $R_i(k)$ increases with $i$ until $i = \ceil{\frac{n}{2}}$ (i.e., the
 middle node) for $k \le \ceil{\frac{n}{2}} - 1$ (this can be observed
 from the expression~\eqref{eq:rik-trid-simple}). For general $k$, it
 can be shown that the value of the sum in~\eqref{eq:rik-trid-full}
 for $R_i(k)$ is strongly dominated by the summand corresponding to
 the index $p = 0$, which increases with $i$ until $i =
 \ceil{\frac{n}{2}}$ and decreases afterwards. Thus, \textbf{the
   optimal control node corresponds always to (one of) the
   center node(s)}, i.e., $b^*(k) = e_{\ceil{\frac{n}{2}}}$ for all
 $k$. If nodes have uniform self-loops (i.e., self-loops all with the
 same weight), $R_i(k)$ can no longer be computed analytically but
 simulations show the exact same behavior;

\begin{proposition}\longthmtitle{$2 k$-communicabilities of ring
    networks}\label{prop:ring}
  Consider a ring network of $n$ nodes and uniform link weights $a$
  (with no self-loops). Then, for $i \in \Nc$ and $k \in
  \intpos$,
  \begin{align}\label{eq:rik-ring}
    R_i(k) = \frac{(2 a)^{2 k}}{n} \Bigg[1 + 2 \sum_{j =
      1}^{\ceil{\frac{n}{2}} - 1} \cos^{2 k} \Big(\frac{2 j
      \pi}{n}\Big) + \delta^\text{E}_n\Bigg],
  \end{align}
  where $\delta^\text{E}_n$ equals one if $n$ is even and zero
  otherwise.
\end{proposition}
 \begin{proof}
   From~\cite[Lemma 1.77]{FB-JC-SM:08cor}, we have $\lambda_j = 2 a
   \cos\frac{2 j \pi}{n}$ and (after normalization of eigenvectors),
   \begin{align*}
     w_{i j} = \begin{cases} \frac{2}{n} \cos^2\frac{2 (i - 1) j
         \pi}{n} & \!\!\text{if} \ 1 \le j < \frac{n}{2},
       \\
       \frac{2}{n} \sin^2\frac{2 (i - 1) (n - j) \pi}{n} &
       \!\!\text{if} \ \frac{n}{2} < j < n,
       \\
       \frac{1}{n} & \!\!\text{if} \ j = n, \ \text{or} \ n \in
       \integers  \
       \text{even and} \ j = \frac{n}{2},
     \end{cases}
   \end{align*}
   for $i, j \in \until{n}$. Note that to normalize the eigenvectors,
   we follow a similar procedure to the one described in the proof of
   Lemma~\ref{prop:line} (setting $s = 2 i$ and substituting $n$ by $n
   - 1$ in~\eqref{eq:trig-identity}). The result then follows by
   substituting these expressions in $R_i(k)$.
 \end{proof}
 
 We can infer from the preceding result that without self-loops, the
 value of $R_i(k)$ is independent of $i$ (as shown
 by~\eqref{eq:rik-ring}) for a uniform ring network, so \textbf{the
   optimal control node is arbitrary for all $k$}. Similar
 result can be proved analytically if the nodes have uniform
 self-loops.

\begin{proposition}\longthmtitle{$2 k$-communicabilities of star
    networks}\label{prop:star}
  Consider a star network given by
  \begin{align}\label{eq:star}
    A = \left[\begin{array}{cc} l_c & a^T \\ a & l_p I_{n -
          1} \end{array}\right],
  \end{align}
  where $a \in \real^{n - 1}$ contains the link weights between the
  center node and peripheral nodes. Then 
  \begin{align}\label{eq:rik-star}
    \notag R_1(k) &= \frac{(\lambda_1 - l_p)^2}{(\lambda_1 - l_p)^2 +
                    \|a\|^2} \lambda_1^{2 k} + \frac{(l_p -
                    \lambda_2)^2}{(l_p - \lambda_2)^2 + \|a\|^2}
                    \lambda_2^{2 k}, 
    \\ 
     R_i(k) &= \frac{a_{i - 1}^2}{(\lambda_1 - l_p)^2 + \|a\|^2}
                    \lambda_1^{2 k} + \frac{a_{i - 1}^2}{(l_p -
                    \lambda_2)^2 + \|a\|^2} \lambda_2^{2 k}  
    + \frac{\|a\|^2 - a_{i - 1}^2}{\|a\|^2} l_p^{2 k},
  \end{align}
  for all $k \in \intpos \cup \{0\}$ and $i \in \{2, \dots, n\}$, where
  \begin{align}\label{eq:star-lambda12}
      \lambda_{1, 2} = \frac{l_c + l_p \pm \sqrt{(l_c - l_p)^2 + 4
          \|a\|^2}}{2}.
    \end{align}
\end{proposition}
 \begin{proof}
     Using the formula
     \begin{align*}
 \left|\begin{array}{cc}
         P & Q \\ R & S
 \end{array}\right|
     = (P - 1) |S| + |S - R Q|,
     \end{align*}
     for scalar $P$, row vector $Q$, column vector $R$, and square
     matrix $S$, and some algebra, we get
       $|s I_n - A| = \big(s^2 - (l_c + l_p) s + l_c l_p - \|a\|^2\big)
       (s - l_p)^{n - 2}$,
     so the eigenvalues of $A$ are given by
     \begin{align}
     \lambda_{3, \ldots, n} = l_p,
     \end{align} 
    and $\lambda_{1, 2}$ in~\eqref{eq:star-lambda12}. Note that we may or may not have $|\lambda_1| \ge \cdots \ge
    |\lambda_n|$ as the order depends on the values of the parameter.
     By solving $(A - \lambda_j I_n) v_j = 0$ for $j = 1, 2$, and then
     using the orthogonality of eigenvectors, we get
     \begin{align}\label{eq:star-v}
       v_{1, 2} \propto \left[\begin{array}{c} \lambda_{1, 2} - l_p \\
           a \end{array}\right], \quad (v_j)_1 = 0 \quad \forall j \in
       \{3, \ldots, n\},
     \end{align}
     where $\propto$ accounts for normalization. The result then
     follows by substituting~\eqref{eq:star-lambda12}-\eqref{eq:star-v}
     into $R_i(1) = \sum_j v_{i j}^2 \lambda_j^2$ separately for $i =
     1$ and $i \ge 2$, and simplifying.
 \end{proof}

 Using this result, if all self-loop weights are the same ($l_c = l_p$
 in~\eqref{eq:star}), then $R_1(1) > R_i(1)$ for all $i \ge 2$
 from~\eqref{eq:star-diff1}. Therefore Theorem~\ref{thm:main-nec}(iii)
 implies that \textbf{the center node is the optimal control node at all
 times}.
 
 \section{Additional Lemmas and Proofs}\label{sec:proofs}

In this section, we formulate and prove a number of lemmas that underlie the main results of this paper and also provide the proofs of the main results presented in the main text. \new{Throughout, $\cplx$ denotes the set of complex numbers and for $M \in \cplx^{n \times n}$, $\overline M$ and $M^*$ denote its complex conjugate and complex conjugate transpose, respectively,
and $M^{-*} = (M^*)^{-1}$. Further, for $\lambda \in \real^n$ and $\ell \in \intpos \cup \{0\}$,
$\lambda^\ell \triangleq [\lambda_1^\ell \ \cdots \ \lambda_n^\ell]^T$
and $|\lambda| \triangleq [|\lambda_1| \ \cdots \
  |\lambda_n|]^T$.}

\begin{proof}[\bfseries Proof of Theorem~\ref{thm:main-suf}]
  Define
  \begin{align*}
    U = V^{-\star}.
  \end{align*}
  Notice that the columns of $U$ are the left eigenvectors of $A$,
  with the same order as in $\Lambda$ and $V$. Since for any $k$,
  \begin{align*}
    (A^k)^T A^k = (A^k)^* A^k = (V \Lambda^k U^*)^* V \Lambda^k U^*,
  \end{align*} 
  it follows that for any $i$ and $k$,
  \begin{align*}
    R_i(k) = \big((A^k)^T A^k\big)_{ii} = (V \Lambda^k U_{i, :}^*)^* V
    \Lambda^k U_{i, :}^* = \|V \Lambda^k U_{i, :}^*\|_2^2, 
  \end{align*} 
  where $U_{i, :}$ denotes the $i$th row of $U$. For simplicity,
  define $c^{(i, k)} = V \Lambda^k U_{i, :}^* \in \cplx^n$. It is
  straightforward to check that
  \begin{align*}
    c^{(i, k)}_\ell = \sum_{j = 1}^n v_{\ell j} \overline u_{ij}
    \lambda_j^k, \qquad \ell \in \until{n} 
  \end{align*} 
  so
  \begin{align*}
    R_i(k) = \sum_{\ell = 1}^n \big|c^{(i, k)}_\ell\big|^2 =
    \sum_{\ell = 1}^n c^{(i, k)}_\ell \overline{c^{(i, k)}_\ell}
    &= \sum_{\ell = 1}^n \sum_{j = 1}^n \sum_{m = 1}^n v_{\ell j}
      \overline v_{\ell m} \overline u_{ij} u_{im} \lambda_j^k
      \overline \lambda_m^k 
    \\
    &= \sum_{j, m = 1}^n \underbrace{\left(\sum_{\ell = 1}^n v_{\ell
      j} \overline v_{\ell m} \overline u_{ij}
      u_{im}\right)}_{\beta^{(j, m)}_i} \lambda_j^k \overline
      \lambda_m^k. 
  \end{align*} 
  Dividing both sides by $\lambda_1^{2k}$ and taking the limits as
  $k \to \infty$, we see that for all $i$,
  $\beta^{(1, 1)}_i = u_{i1}^2 = R_i(\infty)$ (notice that
  $u_{i1} \in \real$ for all $i$ since $\lambda_1 \in \realpos$
  according to the Perron-Frobenius Theorem~\cite[Fact
  4.11.4]{DSB:09}). Choose $r(1) \in \argmax_i R_i(1)$ and
  $r(\infty) \in \argmax_i R_i(\infty)$. The network belongs to class
  $\Vc$ if for some $k > 1$,
  \begin{align*}
    R_{r(\infty)}(k) > R_{r(1)}(k) &\Leftrightarrow R_{r(\infty)}(\infty) \lambda_1^{2k} + \sum_{(j, m) \neq (1, 1)} \beta_{r(\infty)}^{(j, m)} \lambda_j^k \overline \lambda_m^k > R_{r(1)}(\infty) \lambda_1^{2k} + \sum_{(j, m) \neq (1, 1)} \beta_{r(1)}^{(j, m)} \lambda_j^k \overline \lambda_m^k
    \\
                                   &\Leftrightarrow \left[R_{r(\infty)}(\infty) - R_{r(1)}(\infty)\right] \lambda_1^{2k} > \sum_{(j, m) \neq (1, 1)} \left[\beta_{r(1)}^{(j, m)} - \beta_{r(\infty)}^{(j, m)}\right] \lambda_j^k \overline \lambda_m^k
    \\
                                   &\stackrel{\text{(a)}}{\Leftarrow} \left[R_{r(\infty)}(\infty) - R_{r(1)}(\infty)\right] \lambda_1^{2k} > \lambda_1^k |\lambda_2|^k \left|\sum_{(j, m) \neq (1, 1)} \beta_{r(1)}^{(j, m)} - \beta_{r(\infty)}^{(j, m)}\right|
    \\
                                   &\Leftrightarrow \left[R_{r(\infty)}(\infty) - R_{r(1)}(\infty)\right] \lambda_1^{2k} > \lambda_1^k |\lambda_2|^k \sum_{(j, m) \neq (1, 1)} \left|\beta_{r(1)}^{(j, m)}\right| + \left|\beta_{r(\infty)}^{(j, m)}\right|
    \\
                                   &\Leftarrow \left[R_{r(\infty)}(\infty) - R_{r(1)}(\infty)\right] \lambda_1^k > |\lambda_2|^k \cdot 2 \max_{i \in \until{n}} \sum_{j, m = 1}^n \left|\beta_i^{(j, m)}\right|,
\end{align*} 
where in (a) we have used the fact that $|\lambda_j \overline
\lambda_m| \le \lambda_1 |\lambda_2|$ for any $(j, m) \neq (1,
1)$. Now, using the definition of $\beta_i^{(j, m)}$,
\begin{align*}
  \sum_{j, m = 1}^n \left|\beta_i^{(j, m)}\right| &\le \sum_{j, m =
    1}^n \sum_{\ell = 1}^n |v_{\ell j}| |v_{\ell m}| |u_{ij}| |u_{im}|
  \\
  &= \sum_{j, m = 1}^n |u_{ij}| |u_{im}| \left(\sum_{\ell = 1}^n
    |v_{\ell j}| |v_{\ell m}|\right)
  \\
  &\stackrel{\text{(b)}}{\le} \sum_{j, m = 1}^n |u_{ij}| |u_{im}|
  \underbrace{\|V_{:, j}\|_2}_{1} \underbrace{\|V_{:, m}\|_2}_{1}
  \\
  &= \|U_{i, :}\|_1^2
  \\
  &\le \|U\|_\infty^2,
\end{align*} 
where (b) follows from the Cauchy-Schwarz inequality. Thus,
\begin{align*}
  R_{r(\infty)}(k) > R_{r(1)}(k) &\Leftarrow
  \left[R_{r(\infty)}(\infty) - R_{r(1)}(\infty)\right] \lambda_1^k >
  |\lambda_2|^k \cdot 2 \|U\|_\infty^2
  \\
  &\Leftrightarrow k > \frac{\log 2 \|U\|_\infty^2 -
    \log\left[R_{r(\infty)}(\infty) -
      R_{r(1)}(\infty)\right]}{\log\lambda_1 - \log|\lambda_2|}.
\end{align*}
Therefore, the result follows by choosing $K > \bar K$, where $\bar K
= \ceil{\frac{\log 2 \|U\|_\infty^2 - \log\left[R_{r(\infty)}(\infty)
      - R_{r(1)}(\infty)\right]}{\log\lambda_1 - \log|\lambda_2|}}$.
\end{proof}

The following lemma will be useful in the proof of Theorem~\ref{thm:main-nec}.
  
\begin{lemma}\label{lem:diag-W}
  Let $W \in \real^{n \times n}$ be a doubly-stochastic matrix and
  $\gamma \in \realnonneg^n$ be such that $\gamma_1 \ge \cdots \ge
  \gamma_n$. If $ \frac{1 - w_{1 1}}{w_{1 1}} \le \frac{\gamma_1 -
    \gamma_2}{\gamma_1 - \gamma_n}, $ then $ 1 \in \argmax_{1 \le i
    \le n} \ (W \gamma)_i$.
\end{lemma}
\begin{proof}
  Note that we have
  \begin{align*}
    \frac{1 - w_{1 1}}{w_{1 1}} \le \frac{\gamma_1 -
      \gamma_2}{\gamma_1 - \gamma_n} &\Leftrightarrow (\gamma_1 -
    \gamma_2) w_{1 1} \ge (\gamma_1 - \gamma_n) (1 - w_{1 1})
    \\
    &\Leftrightarrow \gamma_n + w_{1 1} (\gamma_1 - \gamma_n) \ge
    \gamma_2 + (1 - w_{1 1}) (\gamma_1 - \gamma_2)
    \\
    &\Rightarrow \forall i \ge 2 \quad \gamma_n + w_{1 1} (\gamma_1
    - \gamma_n) \ge \gamma_2 + w_{i 1} (\gamma_1 - \gamma_2),
  \end{align*} 
  where the last implication is because $w_{i 1} \le 1 - w_{1 1}$
  for all $i \in \until{n}$. The last inequality can be equivalently
  expressed, for any $i \in \{2, \dots, n\}$, as
  \begin{align*}
    w_{1 1} \gamma_1 + (1 - w_{1 1}) \gamma_n \ge w_{i 1} \gamma_1 +
    (1 - w_{i 1}) \gamma_2,
  \end{align*}
  which, given that $\gamma_n \le \gamma_j \le \gamma_2$ for all $j
  \in \{2, \dots, n\}$, implies
  \begin{align*} 
    w_{11} \gamma_1 + \sum_{j = 2}^n w_{1 j} \gamma_j \ge w_{i 1}
    \gamma_1 + \sum_{j = 2}^n w_{i j} \gamma_j,
  \end{align*}
  for any $i \in \{2, \dots, n\}$.  This can be equivalently written
  as
  \begin{align*}
    \sum_{j=1}^n w_{1 j} \gamma_j \ge \sum_{j=1}^n w_{i j} \gamma_j
    \Leftrightarrow (W \gamma)_1 \ge (W \gamma)_i,
  \end{align*}
  completing the proof.
\end{proof}

\begin{proof}[\bfseries Proof of Theorem~\ref{thm:main-nec}]
\new{
For convenience, let $\lambda =
[\lambda_1 \ \cdots \ \lambda_n]^T$.  After some algebraic manipulations, one can
show that
\begin{align}\label{eq:W-lambda2ki}
 R_i(k) = (A^{2 k})_{i i} = \sum_{j=1}^n v_{ij}^2 \lambda_j^{2k} = (W \lambda^{2 k})_i.
\end{align}
The assumption that node $1$ has the largest eigenvector centrality is equivalent to the largest element of the first column of $W$ being $w_{1 1}$, i.e.,
\begin{align}\label{eq:w11}
  w_{1 1} = \max_{1 \le i \le n} w_{i 1},
\end{align}
or, also equivalently, $r(\infty) = 1$. This can always be realized by a permutation of the rows of~$W$
achieved by relabeling the node with the largest centrality as
node~$1$ (note that relabeling the nodes only permutes the rows of $W$ and not its columns.
The order of its columns is arbitrary and corresponds to the
order of the diagonal elements of $\Lambda$).
}

  \new{The claim of the theorem is trivial in all cases for $k=0$.} Under condition \emph{(i)},
  for $k= 1$, we have
  \begin{align*}
    \frac{\lambda_1^2 - \lambda_2^2}{\lambda_1^2 - \lambda_n^2} \!=\!
    \frac{|\lambda_1| - |\lambda_2|}{|\lambda_1| - |\lambda_n|}
    \frac{|\lambda_1| + |\lambda_2|}{|\lambda_1| + |\lambda_n|}
    \!\ge\! \frac{|\lambda_1| - |\lambda_2|}{|\lambda_1| -
      |\lambda_n|} \!\ge\! \frac{1 - w_{1 1}}{w_{1 1}}.
  \end{align*}
  For $k \ge 2$, using the above inequality, we have
  \begin{align*}
    \frac{\lambda_1^{2 k} - \lambda_2^{2 k}}{\lambda_1^{2 k} -
      \lambda_n^{2 k}} &= \frac{\lambda_1^2 - \lambda_2^2}{\lambda_1^2
      - \lambda_n^2} \frac{\lambda_1^{2 k - 2} + \cdots + \lambda_2^{2
        k - 2}}{\lambda_1^{2 k - 2} + \cdots + \lambda_n^{2 k - 2}}
    \ge \frac{1 - w_{1 1}}{w_{1 1}}.
  \end{align*}
  Thus, the result follows from Lemma~\ref{lem:diag-W}.

  Under condition \emph{(ii)}, for any $k \ge 1$,
  \begin{align*}
    1 \in \argmax_{1 \le i \le n} R_i(k) &\Leftrightarrow \sum_{j = 1}^n w_{1
      j} \lambda_j^{2 k} \ge \sum_{j = 1}^n w_{i j} \lambda_j^{2 k}
    \\
    &\Leftrightarrow w_{1 1} \lambda_1^{2 k} + (1 - w_{1 1})
    \lambda_2^{2 k} \ge \sum_{j = 1}^n w_{i j} \lambda_j^{2 k}
    \\
    &\Leftarrow w_{1 1} \lambda_1^{2 k} + (1 - w_{1 1}) \lambda_2^{2
      k} \ge w_{i 1} \lambda_1^{2 k} + (1 - w_{i 1}) \lambda_2^{2 k}
    \\
    &\Leftrightarrow (w_{1 1} - w_{i 1}) (\lambda_1^{2 k} -
    \lambda_2^{2 k}) \ge 0,
  \end{align*}
  where the last inequality is always true (cf. equation~\eqref{eq:w11}).

  Finally, under condition \emph{(iii)}, first consider the case when
  $|\lambda_1| > |\lambda_2|$. By contradiction, assume $R_i(k) >
  R_1(k)$ for some $i \in \{2, \ldots, n\}$ and $k \ge 2$. Since
  $|\lambda_1| > |\lambda_i|$ for all $i \in \{2, \ldots, n\}$, there
  exists a sufficiently large $\overline k$ where $R_1(\overline k) >
  R_i(\overline k)$ (recall our node labeling convention
  in~\eqref{eq:w11}). Note that it is not necessary for $\overline k$
  to be less than $K$. Thus, $R_1$ and $R_i$ swap orders at least $2$
  times. However, since $A$ has (at most) three distinct nonzero
  eigenvalues, \cite[Theorem 1]{TT:06} implies that $R_1$ and $R_i$
  can swap orders at most once, which is a contradiction. On the other
  hand, if $|\lambda_1| = |\lambda_2|$, then each $R_i$ is essentially
  the sum of at most two distinct exponential functions and thus,
  using~\cite[Theorem 1]{TT:06} again, the order of all $R_i$'s
  remains unchanged for all $k$, yielding the result.
\end{proof}

\begin{proof}[\bfseries Proof of Theorem~\ref{thm:sn}]
  We first prove the first part of the theorem for general (not
  necessarily symmetric) $A_0$ and $E$. Recall that for $k \in \{0,
  \dots, K - 1\}$
\begin{align*}
  r(k) = \argmax_{i \in \Nc} R_i(k) &= \argmax_{i \in \Nc} \big(((A +
  \alpha E)^k)^T (A + \alpha E)^k\big)_{ii}
  \\
  &\stackrel{\text{(a)}}{=} \argmax_{i \in \Nc} \big(((\alpha^{-1} A +
  E)^k)^T (\alpha^{-1} A + E)^k\big)_{ii},
\end{align*} 
where (a) holds because the maximizer of a set is invariant to the
scaling of all the elements of the set by a constant. Using
$\lim_{\alpha \to \infty} \alpha^{-1} A + E = E$ and the continuity of
polynomials, we get
\begin{align*}
\lim_{\alpha \to \infty} R_i(k) = \tilde R_i(k),
\end{align*}
where $\tilde R_i$ denotes the $2k$-communicabilities of a node $i$ in
the additive network $E$. Since $E$ is not acyclic, powers of $E$
never vanish, and thus
\begin{align*}
\forall k \in \{0, \dots, K - 1\} \ \exists i \in \until{n_1} \qquad \tilde R_i(k) > 0,
\end{align*}
while $\tilde R_i(k) = 0$ for $i \in \{n_1 + 1, \dots, n\}$ and all $k$. Therefore, for any $k \in \{0, \dots, K - 1\}$, there exists $\bar \alpha_k > 0$ such that
\begin{align*}
r(k) \in \until{n_1},
\end{align*}
for $\alpha > \bar \alpha_k$. The claim follows by taking $\bar \alpha = \max_{k \in \{0, \dots, K - 1\}} \bar \alpha_k$.

Now, assume $A_0$ and $E$ are symmetric. As before, let $\lambda = [\lambda_1 \ \cdots \ \lambda_n]^T \in
    \real^n$ and $V \in \real^{n \times n}$ be the vector of
    eigenvalues (with $|\lambda_1| \ge \cdots \ge
    |\lambda_n|$) and the matrix of eigenvectors of $A$, respectively,
    and $W$ be the element-wise square of $V$. Recall that this gives
\begin{align*}
 R_i(k) = (A^{2 k})_{i i} = \sum_{j=1}^n v_{ij}^2 \lambda_j^{2k} = (W \lambda^{2 k})_i.
\end{align*}
Let $i^* \in \until{n_1}$ be the node with the greatest eigenvector centrality in $E$ and $\gamma \in \real^n$ be any vector such that $\gamma_1 \ge \cdots \ge
  \gamma_n \ge 0$. Fix $i \in \{n_1 + 1, \dots, n\}$ arbitrarily and let $r
  \le n_1$ be the rank of $E$. Using the inequalities
  \begin{align*}
    \sum_{j = 1}^n w_{i^* j} \gamma_j &\ge w_{i^* 1} \gamma_1, \\
    \sum_{j = 1}^r w_{i j} \gamma_j &\le \gamma_1 \sum_{j = 1}^r w_{i
      j},
    \\
    \sum_{j = r + 1}^n w_{i j} \gamma_j &\le \gamma_{r + 1},
  \end{align*} 
  it follows that $(W
  \gamma)_{i^*} > (W \gamma)_i$ if
  \begin{align}\label{eq:sn-proxy}
    w_{i^* 1} \gamma_1 > \gamma_1 \sum_{j = 1}^r w_{i j} + \gamma_{r +
      1}.
  \end{align}
  Note that if~\eqref{eq:sn-proxy} holds for $\gamma = |\lambda|$,
  then it holds for $\gamma = \lambda^{2 k}$ for any $k \ge 1$. This is
  because
  \begin{align*}
    &w_{i^* 1} \lambda_1^{2 k} = |\lambda_1|^{2 k - 1} \cdot w_{i^* 1}
    |\lambda_1|
    > |\lambda_1|^{2 k - 1} \Big(|\lambda_1| \sum_{j = 1}^r w_{i j} +
    |\lambda|_{r + 1}\Big) > \lambda_1^{2 k} \sum_{j = 1}^r w_{i j} +
    \lambda_{r + 1}^{2 k}.
  \end{align*}
  Therefore, our proof strategy is to find $\overline \alpha$ such
  that~\eqref{eq:sn-proxy} holds for $\gamma = |\lambda|$ if
  $\alpha > \overline \alpha$. To this end, let $\tilde
    \lambda = [\tilde \lambda_1 \ \cdots \ \tilde \lambda_n]^T \in
    \real^n$ and $\tilde V \in \real^{n \times n}$ be the vector of
    eigenvalues (with $|\tilde \lambda_1| \ge \cdots \ge
    |\tilde \lambda_n|$) and the matrix of eigenvectors of $E$, respectively,
    and $\tilde W$ be the element-wise square of $\tilde V$. Note that
    $\tilde W$ has the structure 
  \begin{center}
  \begin{tikzpicture}
  \node[draw=none] at (0, 0) (matrix) {
  \parbox{0.48\textwidth}{
  \begin{align}\label{eq:tilde-W}
  \arraycolsep=8pt
  \def\arraystretch{1.4}
  \def\bracketHeight{18pt}
  \tilde W = 
  \left[\rule{0cm}{\bracketHeight}\right.
  \hspace{-5pt}
  \text{\raisebox{1pt}{$\begin{array}{cc} \star & 0 \\ 0 & \star \end{array}$}}
  \hspace{-5pt}
  \left]\rule{0cm}{\bracketHeight}\right.
  \end{align}
  }
  };
  \node at (37pt, 8pt) (right-top1) {$\big\}$};
  \node[right of=right-top1, xshift=-21pt] (right-top2) {$\scriptstyle n_1$};
  \node at (37pt, -9pt) (right-bottom1) {$\big\}$};
  \node[right of=right-bottom1, xshift=-15pt] (right-bottom2) {$\scriptstyle n - n_1$};
  \node[rotate=270] at (0pt, 18pt) (top-left1) {$\big\{$};
  \node[above of=top-left1, yshift=-23pt] (top-left2) {$\scriptstyle n_1$};
  \node[rotate=270] at (21pt, 18pt) (top-right1) {$\big\{$};
  \node[above of=top-right1, yshift=-22pt] (top-right2) {$\scriptstyle n - n_1$};
  \node[right of=matrix, xshift=35pt] (dot) {$.$};
  \end{tikzpicture}
  \end{center}
  \vskip -13pt
  \noindent In the following, we bound $\lambda$ and $V$ using
  perturbation theory of eigenvalues and eigenvectors.  For simplicity
  of exposition, we only deal with the case where the $r$ nonzero
  eigenvalues of $E$ are all distinct (the proof for the general case
  proceeds along the same lines but is more involved).

  To bound the eigenvalues in $\lambda$, let $\pi_A:\until{n} \to
  \until{n}$ be a permutation that re-orders the eigenvalues of $A$
  based on their \emph{signed} value, i.e., $\lambda_{\pi_A(1)} \ge
  \lambda_{\pi_A(2)} \ge \cdots \ge \lambda_{\pi_A(n)}$.  Define
  $\pi_E$ similarly for $E$ (i.e., such that $\tilde
    \lambda_{\pi_E(1)} \ge \tilde \lambda_{\pi_E(2)} \ge \cdots \ge
    \tilde \lambda_{\pi_E(n)}$). By Weyl's Theorem~\cite[Thm
  4.3.1]{RAH-CRJ:85},
  \begin{align}\label{eq:weyl}
    |\lambda_{\pi_A(j)} - \alpha \tilde \lambda_{\pi_E(j)}| \le
    \rho(A_0),
  \end{align}
  for all $j \in \until{n}$. We know from the Perron-Frobenius
    theorem~\cite[Fact 4.11.4]{DSB:09} for nonnegative matrices that
    $\pi_A(1) = \pi_E(1) = 1$. Therefore, \eqref{eq:weyl} implies that
    \begin{subequations}\label{eq:lambda-pert}
      \begin{align}
        \alpha \rho(E) - \rho(A_0) \le \lambda_1 \le \alpha \rho(E) +
        \rho(A_0).
      \end{align}
      Moreover, since $E$ has $n - r$ zero eigenvalues, \eqref{eq:weyl}
      implies that $A$ has \emph{at least} $n - r$ eigenvalues with
      absolute value less than or equal to $\rho(A_0)$, i.e.,
      \begin{align}
        |\lambda_{r + 1}| \le \rho(A_0).
      \end{align}
    \end{subequations}
    Next, we bound the eigenvectors in~$V$.  Define
    \begin{align*}
      \delta_E = \min\{\tilde \lambda_{\pi_E(j)} - \tilde
      \lambda_{\pi_E(j + 1)} \; | \; \tilde \lambda_{\pi_E(j)} - &
      \tilde \lambda_{\pi_E(j + 1)} > 0, \
      j \in \until{n - 1}\}.
    \end{align*}
    Using~\cite[Cor. 1]{YY-TW-RJS:15}, we have
    \begin{align}\label{eq:dk}
      \|v_{\pi_A(j)} - \tilde v_{\pi_E(j)}\| \le \frac{2^{3/2}
        \|A_0\|}{\alpha \delta_E},
    \end{align}
    for $j \in \pi_E^{-1}(\until{r})$. To see this, set
      $\Sigma = \alpha E$ and $\hat \Sigma = A_0$
      in~\cite[Cor. 1]{YY-TW-RJS:15}. This is the only place where we
      need the nonzero eigenvalues of $E$ to be distinct. If $E$ has a
      repeated nonzero eigenvalue, then the corresponding eigenvectors
      are no longer unique, i.e., one has to study the perturbation of
      eigenspaces rather than eigenvectors. Therefore, one can no
      longer use the simplified variant~\cite[Cor. 1]{YY-TW-RJS:15} of
      the Davis-Kahan Theorem but the original result itself, which
      provides essentially the same result but is more technically
      involved.
      
      Using $\pi_A(1) = \pi_E(1) = 1$ and~\eqref{eq:dk},
    we get
    \begin{align}\label{eq:w-bounding}
      |w_{i^* 1} - \tilde w_{i^* 1}| &= |v_{i^* 1}^2 - \tilde v_{i^*
        1}^2| \le 2 | |v_{i^* 1}| - |\tilde v_{i^* 1}| |
      \\
      \notag &\le 2 |v_{i^* 1} - \tilde v_{i^* 1}| \le 2 \|v_1 - \tilde
      v_1\| \le \frac{2^{5/2} \|A_0\|}{\alpha \delta_E},
    \end{align}
    which together with $\tilde w_{i^* 1} \ge \frac{1}{n_1}$ gives
   \begin{subequations}\label{eq:V-pert}
     \begin{align}
       w_{i^* 1} \ge \frac{1}{n_1} - \frac{2^{5/2} \|A_0\|}{\alpha
         \delta_E}.
     \end{align}
     To derive similar bounds on $w_{i j}, j \in \until{r}$ (recall
     that we fixed $i \in \{n_1 + 1, \dots, n\}$ arbitrarily at the beginning of
     the proof), we need to choose $\alpha > \frac{2
       \rho(A_0)}{|\tilde \lambda_r|}$. This choice of $\alpha$
     guarantees that $\pi_A(j) \in \until{r}$ for all $j \in
     \pi_E^{-1}(\until{r})$. Therefore, using~\eqref{eq:dk}
     and~\eqref{eq:tilde-W} and following the same steps as
     in~\eqref{eq:w-bounding}, we get
     \begin{align}
       w_{i j} \le \frac{2^{5/2} \|A_0\|}{\alpha \delta_E}, \quad j
       \in \until{r}.
     \end{align}
   \end{subequations}
   Now, using~\eqref{eq:lambda-pert} and~\eqref{eq:V-pert},
   \eqref{eq:sn-proxy} holds with $\gamma = |\lambda|$ if
   \begin{align*}
     \left(\frac{1}{n_1} - \frac{2^{5/2} \|A_0\|}{\alpha
         \delta_E}\right) &\left(\alpha \tilde \lambda_1 -
       \rho(A_0)\right)
     > \left(\alpha \tilde \lambda_1 + \rho(A_0)\right)
     \frac{r 2^{5/2} \|A_0\|}{\alpha \delta_E} + \rho(A_0),
   \end{align*}
   which itself holds if $\alpha > \overline \alpha$, where
   \begin{align*}
     \overline \alpha \triangleq \max&\Big\{\!1, \frac{2
       \rho(A_0)}{\tilde \lambda_r}, \frac{8 \|A_0\|}{\delta_E} \Big(1
     \!+ \frac{\rho(A_0)}{\rho(E)}\Big) n_1^2 + 2
     \frac{\rho(A_0)}{\rho(E)} n_1\!\Big\} ,
  \end{align*}
  completing the proof.
\end{proof}

\section{Description of the Analyzed Real Networks}\label{sec:real}

The real networks studied in this work have been acquired from a multitude of sources, which we list here for easier reproduction of our results. All the databases are freely and publicly available.

\begin{itemize}
\item \textbf{BCTNet fMRI~\cite{MR-OS:10}:} This is a human
  whole-brain functional network. Nodes represent brain areas and
  edges represent fMRI co-activations. The dataset is available online
  at \url{https://sites.google.com/site/bctnet/datasets}.
\item \textbf{Cocomac~\cite{RB-TW-MD:12}:} This is a macaque
  whole-brain structural network based on the Felleman and Van Essen
  atlas. Nodes represent brain areas and edges represent axonal
  projections (nerve tracts) between them. The dataset is retrieved
  from \url{http://cocomac.g-node.org/services/axonal_projections.php}
  by entering the specifications detailed in
  \url{http://cocomac.g-node.org/main/faq.php#connectivity matrix}.
\item \textbf{BCTNet Cat~\cite{MR-OS:10}:} This represents the cat
  structural thalamocortical network. Nodes represent thalamocortical
  areas and edges represent nerve tracts between them. The dataset is
  available online at
  \url{https://sites.google.com/site/bctnet/datasets}.
\item \textbf{C. elegans~\cite{DJW-SHS:98}:} This dataset contains the
  neural network of Caenorhabditis elegans worm (C. elegans). Nodes
  represent individual neurons and edges represent the total number of
  synapses and gap junctions between any pair of neurons. The dataset of
  available online at \url{https://toreopsahl.com/datasets/#celegans}.
\item \textbf{air500~\cite{VC-RPS-AV:07}:} This is the network of the
  500 busiest commercial airports in the United States in 2002. Nodes
  represent airports and edges represent flights between them. The dataset
  is available online at
  \url{https://toreopsahl.com/datasets/#usairports}.
\item \textbf{airUS~\cite{TO:11}:} This is the complete US airport
  network in 2010. Nodes and edges represent airports and flights
  between them, respectively. The dataset is available online at \url{https://toreopsahl.com/datasets/#usairports}.
\item \textbf{airGlobal~\cite{TO:11}:} This dataset contains the
  global airport network according to \url{OpenFlights.org}. Nodes and
  edges represent airports and flights between them, respectively. The
  dataset is available online at  \url{https://toreopsahl.com/datasets/#usairports}.
\item \textbf{Chicago~\cite{RWE-KSC-YJL-DEB:83,DEB-KSC-MEF-YJL-KTL-RWE:85}:} This dataset represents the road transportation network of the Chicago region, USA. Nodes are transport nodes while edges represent connections between them. The dataset is available online as~\cite{KONECT-Chicago:17}.
\item \textbf{E. coli~\cite{HK-JES-JS-IL:15}:} This is the probabilistic functional gene network of E. coli. Nodes represent genes and edges represent interactions between them. The dataset is available online at \url{http://www.inetbio.org/ecolinet/downloadnetwork.php} (The integrated network).
\item \textbf{Yeast~\cite{DB-YZ-LC-HX-XZ-HL-JZ-SS-LL-NZ-GL-RC:03}:} This network represents the protein-protein interaction (PPI) network in the budding yeast. Nodes and edges represent proteins and the interactions among them, respectively. The dataset is available online at \url{http://vlado.fmf.uni-lj.si/pub/networks/data/bio/Yeast/Yeast.htm}.
\item \textbf{Stelzl~\cite{US-UW-ML-CH-FHB-HG-MS-MZ-AS-SK-JT-SM-CA-NB-SK-AG-ET-AD-SK-BK-WB-HL-EEW:05}:} This is a protein-protein interaction network in humans. Nodes and edges represent proteins and the interactions among them, respectively. The dataset is available online as~\cite{KONECT-Stelzl:17}.
\item \textbf{Figeys~\cite{RME-PC-FE-HL-PT-SC-LM-MDR-LO-ML-RT-MD-YH-AH-LM-SZ-OO-YVB-ME-YS-JV-MA-JPPL-HSD-IIS-BK-KH-KC-KG-BM-RK-SLLA-MFM-GBM-TT-DF:07}:} Similar to above, this is a protein-protein interaction network in humans where nodes and edges represent proteins and the interactions among them, respectively. The dataset is available online as~\cite{KONECT-Figeys:17}.
\item \textbf{Vidal~\cite{JR-KV-TH-THK-AD-NL-GFB-FDG-MD-NAG:05}:} Similar to above, this is a protein-protein interaction network in humans where nodes and edges represent proteins and the interactions among them, respectively. The dataset is available online as~\cite{KONECT-Vidal:17}.
\item \textbf{westernUS~\cite{DJW-SHS:98}:} This dataset describes the high voltage power grid in the Western States of the US. Nodes represent transformers, substations, and generators, and the edges represent high-voltage transmission lines. The dataset is available online at \url{https://toreopsahl.com/datasets/#uspowergrid}.
\item \textbf{Florida~\cite{REU-JJH-MSE:00}:} This network describes
  the food web in the cypress wetlands of South Florida during the wet
  season. Nodes represent taxa and an an edge denotes that a taxon
  uses another taxon as food. The dataset is available online
  as~\cite{KONECT-Florida:17}.
\item \textbf{LRL~\cite{NDM-JJM-TK-MS:91}:} The networks describes the
  food web of Little Rock Lake, Wisconsin, USA. Nodes represent
  autotrophs, herbivores, carnivores and decomposers while links
  represent food sources. The dataset is available online
  as~\cite{KONECT-LRL:17}.
\item \textbf{Facebook group~\cite{JL-JJM:12}:} This dataset describes the
  social interactions among a group of Facebook users. Nodes and edges
  represent profiles and the connections between them,
  respectively. The dataset is available online at
  \url{http://snap.stanford.edu/data/egonets-Facebook.html}.
\item \textbf{E-main~\cite{HY-ARB-JL-DFG:17,JL-JK-CF:07}:} This
  datasets contains E-main communications in a research
  institution. Nodes represent institution members and edges exist
  between any ordered pair of members if one has sent at least one
  E-main to the other. The dataset is available online at
  \url{http://snap.stanford.edu/data/email-Eu-core.html}.
\item \textbf{Southern Women~\cite{AD-BBG-MRG:41}:} This is a social
  network of 18 Southern women. Nodes are individuals and edges
  represent mutual attendance at one of the 14 events recorded. The
  dataset is available online at
  \url{https://toreopsahl.com/datasets/#southernwomen}.
\item \textbf{UCI P2P~\cite{TO-PP:09}:} This dataset describes an
  online community among the students of the University of California,
  Irvine. Nodes represent individuals and edges represent at least one
  message sent between any pair of them. The dataset is available
  online at
  \url{https://toreopsahl.com/datasets/#online_social_network}.
\item \textbf{UCI Forum~\cite{TO:13}:} This network is based on the
  same online community as in UCR P2P, but an edge exists between two
  individuals if they posted on the same topic in a forum. This
  dataset is also available online at
  \url{https://toreopsahl.com/datasets/#online_social_network}.
\item \textbf{Freeman's EIES~\cite{SCF-LCF:79}:} This is a network of
  researchers working on social network analysis. Nodes represent
  researchers and edges represent personal relationships between
  them. The dataset is available online at
  \url{https://toreopsahl.com/datasets/#FreemansEIES} (the second
  dataset in the list).
\item \textbf{Dolphins~\cite{DL-KS-OJB-PH-ES-SMD:03}:} This is a
  social network of bottlenose dolphins observed between 1994 and
  2001. The nodes are the bottlenose dolphins and edges indicate a
  frequent association between them. The dataset is available online
  as~\cite{KONECT-dolphins:17}.
\item \textbf{Physicians~\cite{JC-EK-HM:57}:} This network captures
  innovation spread among 246 physicians in four towns in Illinois,
  USA. A node represents a physician and an edge represents that one
  physician recognizes the other as theor friend or that they turn to
  them if they need advice or are interested in a discussion. The
  dataset is available online as~\cite{KONECT-physicians:17}.
\item \textbf{Org. Consult Advice \& Value~\cite{RLC-AP:04}:} These
  are intra-organizational networks between employees of a consulting
  company. The nodes are individuals, and the edges represent
  frequency of information or advice requests (Org. Consult Advice)
  and the value placed on the information or advice received
  (Org. Consult Value). The datasets are available online at
  \url{https://toreopsahl.com/datasets/#Cross_Parker}.
\item \textbf{Org. R\&D Advice \& Aware~\cite{RLC-AP:04}:} Similar to
  the networks above, these describe intra-organizational interactions
  among the members of a research team in a manufacturing
  company. Nodes represent individuals, and edges represent the extent
  to which individuals received advice from their peers to accomplish
  their work (Org. R\&D Advice) and employees' awareness of each
  others' knowledge and skills (Org. R\&D Aware). The datasets are
  available online at
  \url{https://toreopsahl.com/datasets/#Cross_Parker}.
\end{itemize} 

\bibliographystyle{IEEEtran}

\end{document}